%% file: Transitions_SPDEs_arXiv.tex
\documentclass[11pt,reqno]{amsart}

\pdfoutput=1
\numberwithin{equation}{section} \oddsidemargin=-.0cm
\usepackage{amscd}
\usepackage{amsmath}
    \usepackage{amssymb}
\usepackage{graphicx}%
\usepackage{bm} 
\usepackage{needspace}
\usepackage{newtxtext}
\usepackage{booktabs,caption} 

\usepackage[colorlinks=true, pdfstartview=FitV, linkcolor=blue,
            citecolor=blue, urlcolor=blue]{hyperref} 
\usepackage{subcaption}
\usepackage{sidecap}
\usepackage{grffile}
\usepackage{xcolor}

\usepackage{scalerel}

\graphicspath{{./Figures/}}

\providecommand{\U}[1]{\protect\rule{.1in}{.1in}}

\newcommand{\deltaR}{\delta \! \Ra}
\newcommand{\mkr}{\color{black}}

\def\bea{\begin{equation} \begin{aligned}}
\def\eea{\end{aligned} \end{equation}}
\def\beas{\begin{equation*} \begin{aligned}}
\def\eeas{\end{aligned} \end{equation*}}
\def\bes{\begin{equation*}}
\def\ees{\end{equation*}}
\def\Forall{\text{ } \forall \:}
\def\d{\, \mathrm{d}}

\def\be{\begin{equation}}
\def\ee{\end{equation}}
\def\adots{
  \mathinner{\mkern1mu\raise1pt\hbox{.}\mkern2mu\raise4pt\hbox{.}
  \mkern2mu\raise7pt\vbox{\kern7pt\hbox{.}}\mkern1mu}}

\def\vecu{\bm{u}}
\def\Ra{R}  
\def\Pr{\text{\rm Pr}}
\def\D{\mathcal{D}} 

\def\diffusion{\sigma}

\def\s{\mathfrak{s}}
\def \c{\mathfrak{c}}

\def\Id{{\rm Id}}
\def\Re{{\rm Re \,}}
\def\Im{{\rm Im \,}}

\def\cH{H}
\def\adots{
  \mathinner{\mkern1mu\raise1pt\hbox{.}\mkern2mu\raise4pt\hbox{.}
  \mkern2mu\raise7pt\vbox{\kern7pt\hbox{.}}\mkern1mu}}
  
\newtheorem{thm}{Theorem}[section]
\newtheorem{lem}{Lemma}[section]
\newtheorem{defi}{Definition}[section]

\newtheorem{rem}{Remark}[section]
\newtheorem{cor}{Corollary}[section]

\usepackage{thmtools}
\declaretheoremstyle[
spaceabove=6pt, spacebelow=6pt,
headfont=\normalfont\bfseries,
notefont=\mdseries, notebraces={(}{)},
bodyfont=\it,
postheadspace=0.4em,
headpunct=.
]{hyp_style}

\declaretheorem[style=hyp_style, name=Condition, preheadhook={}]{hypA}
\declaretheorem[style=hyp_style, name=Condition, preheadhook={}]{hypP}
\declaretheorem[style=hyp_style, name=Condition, preheadhook={}]{hypL}
\def\bh{\begin{hypA}}
\def\eh{\end{hypA}} 
\def\bhp{\begin{hypP}}
\def\ehp{\end{hypP}} 
\def\bhl{\begin{hypL}}
\def\ehl{\end{hypL}} 

\def\bt{\begin{thm}}
\def\et{\end{thm}}
\def\bl{\begin{lem}}
\def\el{\end{lem}}
\def\bd{\begin{defi}}
\def\ed{\end{defi}}
\def\bc{\begin{cor}}
\def\ec{\end{cor}}
\def\bp{\begin{proof}}
\def\ep{\end{proof}}
\def\br{\begin{rem}}
\def\er{\end{rem}}
\def\bi{\begin{itemize}}
\def\ei{\end{itemize}}

\title[Transitions in Stochastic Non-equilibrium Systems]{Transitions in Stochastic Non-equilibrium Systems: \\
Efficient Reduction and Analysis}
\author[M.~D.~Chekroun]{Micka\"el D.~Chekroun}
\address[MDC]{Department of Atmospheric and Oceanic Sciences, University of California, Los Angeles, CA 90095-1565, USA, and\\
Department of Earth and Planetary Sciences, Weizmann Institute,
Rehovot 76100, Israel
}
\email{mchekroun@atmos.ucla.edu}  

\author[H.~Liu]{Honghu Liu}
\address[HL]{Department of Mathematics, Virginia Tech, Blacksburg, VA 24061, USA}
\email{hhliu@vt.edu}

\author[J.~C.~McWilliams]{James C.~McWilliams}
\address[JCM]{Department of Atmospheric and Oceanic Sciences, University of California, Los Angeles, CA 90095-1565, USA}
\email{jcm@atmos.ucla.edu}

\author[S.~Wang]{Shouhong Wang}
\address[SW]{Department of Mathematics,
Indiana University, Bloomington, IN 47405}
\email{showang@indiana.edu}

\date{\today}

\begin{document}
\maketitle

\begin{abstract}
A central challenge in physics is to describe non-equilibrium systems driven by randomness, such as a randomly growing interface, or fluids subject to random fluctuations that account e.g. for local stresses and heat fluxes in the fluid which are not related to the velocity and temperature gradients. For deterministic systems with infinitely many degrees of freedom, normal form and center manifold theory have shown a prodigious efficiency to often completely characterize how the onset of linear instability translates into the emergence of nonlinear patterns, associated with genuine physical regimes. However, in presence of random fluctuations, the underlying reduction principle to the center manifold is seriously challenged due to large excursions caused by the noise, and the approach needs to be revisited. 

In this study, we present an alternative framework to cope with these difficulties exploiting the approximation theory of stochastic invariant manifolds, on one hand, and energy estimates measuring the defect of parameterization of the high-modes, on the other. To operate for fluid problems subject to stochastic stirring forces, these error estimates are derived under assumptions regarding dissipation effects brought by the high-modes in order to
suitably counterbalance the loss of regularity due to the nonlinear terms. 
As a result, the approach enables us to predict, from reduced equations of the stochastic fluid problem, the occurrence in large probability of a stochastic analogue to the pitchfork bifurcation, as long as the noise's intensity and the eigenvalue's magnitude of the mildly unstable mode scale accordingly.

In the case of SPDEs forced by a multiplicative noise in the orthogonal subspace of e.g. its mildly unstable mode, our parameterization formulas show that the noise gets transmitted to this mode via non-Markovian coefficients, and that the reduced equation is only stochastically driven by the latter.  These coefficients depend explicitly on the noise path's history, and their memory content is self-consistently determined by the intensity of the random force and its interaction through the SPDE's nonlinear terms. Applications to a stochastic Rayleigh-B\'enard problem  are detailed, for which conditions for a stochastic pitchfork bifurcation (in large probability) to occur, are clarified.
\end{abstract}

{\small
\tableofcontents
}

\section{Introduction} \label{Sec_intro}

The development of instability in a physical system may often be described in terms of the temporal evolution of the amplitudes of certain normal modes, namely those that are mildly unstable and those that are only slightly damped in linear theory. When the number of these nearly marginal modes is finite, their amplitudes are governed by ordinary differential equations (ODEs) in which the growth rates of the linear theory have been renormalized by nonlinear terms \cite{crawford1991introduction,crawford1991symmetry}.
Intuitively, the reason for this reduction is a simple separation of time scales. Modes that have just crossed the imaginary axis have a small real part and are evolving
on long time scales, all the other fast modes rapidly adapting themselves to the slow modes. 

 For deterministic dynamical physical systems with infinitely many degrees of freedom (d.o.f.), normal form and center manifold theory \cite{Hen81,MW05} along with their recent extension \cite{MW14} have shown a prodigious efficiency to often completely characterize how the onset of linear instability translates into the emergence of nonlinear patterns, associated with genuine physical regimes. To do so, the theory along with its recent advances identify, for a given problem, a nonlinear reduction mapping that permits to reduce the many d.o.f.~to a few essential variables (governed by reduced ODEs), able to predict in turn the dynamical transitions \cite{MW14}.

However, in presence of random fluctuations, the reduced ODEs produced by these deterministic theories are no longer valid objects to account for phenomena pertaining to the stochastic realm such as noise-induced transitions. In the stochastic context, the reduction mapping needs indeed to account for these random fluctuations. This mapping consists of e.g.~a parameterization of the stable modes  which becomes then stochastic due to noise. The existence of such reduction mappings has been abundantly studied under various settings for stochastic partial differential equations (SPDEs) \cite{BF95,DPD96,DLS03,DLS04}, but the question of their efficient calculation or approximation has been much less addressed, mainly in some special cases \cite{BHP07} or by means of deterministic (and thus incomplete) formulas \cite{blomker2010qualitative}.

 Only in the recent years, the derivation of general, stochastic, approximation formulas to the leading-order (in the nonlinear terms) of stochastic invariant manifolds (IMs) appeared in \cite{CLW15_vol1} for a broad class of SPDEs subject to a multiplicative noise. Yet, in spite of their relevance demonstrated through numerical examples  \cite{CLW15_vol2}, 
the usage of such formulas to derive reduced stochastic differential equations (SDEs) aimed at predicting dynamical transitions in SPDEs, calls for more understanding.\footnote{In fact, the bifurcation analysis for SPDEs is available only for a few particular cases \cite{CLR01,wang2014existence}. Even for stochastic ODEs, the question of how to describe a  stochastic bifurcation is not completely settled; see e.g.~\cite[Chap.~9]{Arnold98} and \cite{CF98}.}

 Indeed, what makes the success of center/unstable manifold theory to describe bifurcations for PDEs---the so-called reduction principle \cite{crawford1991introduction,SY02}---from low-dimensional reduced systems fails to operate in the stochastic realm, calling for a new approach to apprehend bifurcation analysis in presence of noise. The reason is that a breakdown of the reduction principle to the center/unstable manifold now takes place in the stochastic setting, and an exponential slaving of the stable modes onto the unstable ones  is no longer valid, due to large excursions caused by the (white) noise, even when its intensity is small. At a technical level, the exponential attraction by invariant manifolds that holds in presence of a sufficiently large spectral gap
and that extends to SPDEs with globally Lipschitz nonlinearity \cite[Corollary 4.3]{CLW15_vol1}, is violated in the stochastic setting for SPDEs with nonlinearities that are only locally Lipschitz,  such as encountered in stochastic fluid problems. In these cases, the solution is not guaranteed to stay within the neighborhood over which the local stochastic IMs exist, due to noise-induced deviations.

In this article we bring  justification regarding the usage of reduced SDEs built upon stochastic IMs and their approximations formulas based on \cite{CLW15_vol1} to predict dynamical transitions in SPDEs,  via estimates measuring  the  error between the solutions to the reduced SDEs and those of the original SPDE. 
More generally, the goal of this article is twofold: (i) To advance the fundamental understanding for the efficient derivation of reduced 
models able to describe spatiotemporal pattern formations and their transitions subject to random fluctuations in non-equilibrium systems described by SPDEs, and (ii) to derive error estimates relevant for the analysis of the spatio-temporal transitions, in a random environment.

Other approaches have shown their usefulness to derive low-dimensional reduced systems able to capture the essential macroscopic dynamics of SPDEs. Examples include the averaging method \cite{cerrai2009averaging,kifer2004some,pavliotis2008multiscale,wang2013macroscopic}, the amplitude equation approach \cite{blomker2004multiscale,Blomker07,blomker2020impact,BHP07}, and approaches rooted in the singular perturbation theory of Markov processes \cite{MTV01,kurtz1973limit}. 

However, it remains challenging to derive reduced systems and rigorous error estimates to predict dynamical transitions in e.g.~stochastic fluid problems involving typically quadratic nonlinearities. For such problems, the main issue lies in the fact that the nonlinearity does not act only on the dominant modes but also influences them through non-dominant modes \cite{blomker2005approximation,BHP07}, and thus makes central the way the noise interacts with the nonlinear dynamics, either forcing the dominant or the neglected modes, or both.   This is where, as explained below, the stochastic reduction mappings provided by the approximation theory of stochastic IMs  \cite{CLW15_vol1} show their usefulness.
As we will see, these mappings parameterize the neglected variables beyond Ornstein-Uhlenbeck approximations (cf.~\cite{blomker2005center,MTV01,blomker2009amplitude})  by producing essential variables able to predict dynamical transitions (with large probability) for a broad class of SPDEs with quadratic nonlinearities (such as arising in fluid problems), regardless of whether the noise acts only on the resolved or unresolved scales, or both (cf.~\cite{BHP07}).

The main difficulty to perform an efficient stochastic reduction for SPDEs is thus to track properly the noise-path dependence which in many situations requires to be properly parameterized to render account for noise-induced phenomena and other interactions between the noise and nonlinear dynamics. 
This is the purpose of Theorem~\ref{Thm_CM_approx} recalled below  that summarizes  from \cite{CLW15_vol1} general formulas of stochastic reduction mappings based on the approximation results of \cite[Thm.~6.1 \& Cor.~6.1]{CLW15_vol1}. These results entail that the  corresponding stochastic reduction mappings not only provide the leading-order approximation of the underlying stochastic IM but also the tracking of its time-stochastic dependence  via explicit  random coefficients; see  \eqref{Phi_n_general} below.  As already pointed out  in \cite[Chap.~5]{CLW15_vol2} and more substantiated here,  these coefficients are non-Markovian as dependent  explicitly on the history of the noise path, conveying thus exogeneous memory effects \cite{HO07, Hai09} rooted in the (very) small-scale fluctuations (noise).\footnote{These exogeneous memory effects are to be contrasted with the endogeneous ones encountered in the reduction of nonlinear autonomous systems as predicted by the Mori-Zwanzig (MZ) theory; see e.g.~\cite{Chorin2002,GKS04,zhu2018estimation,falkena2019derivation,santos2021reduced}. The former are functionals of the ``past'' of the noise and emerge in the reduction of stochastic systems by means of stochastic IMs \cite{CLW15_vol2} whereas the latter are functionals of the past of the resolved variables that arise in the reduction of autonomous systems when the validity of (deterministic) IMs or, more generally, of the conditional expectation alone in the MZ expansion, breaks down; see \cite{Stinis06,CLM16_Lorenz9D,CLM19_closure,Chekroun2021c}.}

Mathematically, these coefficients are shown to  correspond to stationary solutions of auxiliary SDEs whose drift part depends on the distance to certain resonances between the critical eigenvalues (losing stability)  and the stable part of the spectrum; see \eqref{M_eq} below.  Such a distance to resonance controls in turns the decay of temporal correlations of these coefficients, as  well as their statistics, allowing for non-Gaussian ones---unlike Ornstein-Uhlenbeck processes encountered in other reduction methods  \cite{blomker2005center,MTV01,blomker2009amplitude}---with tails all the more pronounced than their ratio with the noise amplitude approaches unity; see Appendix~\ref{Appendix_A}.

As already shown through numerical examples  in \cite{CLW15_vol2}, these non-Markovian coefficients constitute in fact the key ingredients for the reduced SDEs to achieve good performance, and can be seen as part of the essential variables produced by the underlying stochastic reduction mapping.  For instance, in the case when the noise acts only on the (unresolved) stable modes, such a reduction mapping shows that the noise gets transmitted to the (resolved) unstable modes via these non-Markovian coefficients, and the reduced equation is only stochastically driven by the latter. In contrast, any deterministic parameterization of the unresolved modes will lead in this case  to a deterministic ODE system, which is insufficient to capture any noise-induced phenomena in the original system. In Sec.~\ref{Sec_23} below, we illustrate on a simple two-dimensional system, the ability to capture its stochastic dynamics and its large excursions from such non-Markovian reduced equations.

 In this article, we go beyond these numerical examples and show that these non-Markovian coefficients play furthermore a key role to conduct a rigorous error analysis thanks to the auxiliary SDEs they satisfy. 
 Our focus is on SPDEs involving energy-preserving quadratic nonlinearities motivated by physical problems such as the prediction and characterization of pattern formation arising in fluid problems after the onset of instability and subject to noise disturbance.   In order to circumscribe the difficulties while conveying the main ideas, we place ourselves in the case of a supercritical pitchfork bifurcation, assuming that exactly one eigenmode loses its stability at the critical parameter value. An attentive reader will notice that the framework laid out here is actually not limited to the pitchfork scenario  and adaptable to many other bifurcation/transition scenarios, including the Hopf bifurcation and other attractor bifurcations.

 At the basis of our error estimates is a ``Markovianisation'' of the one-dimensional non-Markovian reduced SDE governing the local dynamics in the case of a pitchfork bifurcation scenario; cf.~\eqref{Eq_reduced} and \eqref{Eq_X_Phi} below. In this reformulation, the  dynamical equation satisfied by the stochastic parameterization $\Phi(t)$ is obtained by It\^o's formula, which is made explicit thanks to the aforementioned auxiliary SDEs determining the random coefficients involved in $\Phi$; see again \eqref{M_eq}.  By recasting the original SPDE  into a system of two coupled equations  satisfied by the unstable and stable components of the SPDE solution (see \eqref{Eq_Quadratic_v2}), the estimates become facilitated, exploiting in particular the dynamical equation satisfied by $\Phi(t)$ (in \eqref{Eq_X_Phi}) aimed at parameterizing the stable component of the SPDE solution.

{\mkr The error estimates assessing the quality of our non-Markovian reduced SDEs to predict dynamical transitions in SPDEs, are then organized through Secns.~\ref{Sec_modeling_error} and \ref{Sec_err_estimates} as follows.} 
After establishing in Sec.~\ref{Sec_apriori_bounds} a suitable a priori bound for the solution $(X(t),\Phi(t))$ of the  Markovianized reduced system \eqref{Eq_X_Phi}, we derive in Sec.~\ref{Sec_residual_estimates} a residual error estimate for the non-Markovian reduced SDE \eqref{Eq_reduced}; see Lemma~\ref{Lem_apriori_est} and Theorem~\ref{Thm_residual} below. This residual error estimate gives---with large probability and over long time intervals---the error made by $(X(t),\Phi(t))$ in satisfying the SPDE model \eqref{Eq_Quadratic_v2}. Our result shows that this error is all the smaller as the eigenvalue $\epsilon$  associated with the unstable mode  is small and the noise intensity scales as $\sqrt{\epsilon}$.

We then derive in Sec.~\ref{Sec_err_estimates} rigorous error estimates between the solutions of the reduced systems and those of the original SPDE emanating from small initial data. The result presented in Theorem~\ref{Thm_error_estimates} shows that the same bounds controlling the residual error in Theorem~\ref{Thm_residual} control also in large probability and over long time intervals,  the error between the solutions of the reduced systems and those of the SPDE.

The derivation of this error estimate requires suitable a priori bounds of the SPDE solution itself. 
As explained in Sec.~\ref{Sec_dicuss_bounds} below, the main difficulty in deriving these a priori bounds compared to those for $(X(t),\Phi(t))$ 
lies in the infinite-dimensional character of the projection onto the stable subspace of the SPDE solution to which applies the  
 loss of regularity via the nonlinearity $B$, requiring {\it de facto} new estimates.  
 In contrast, $\Phi(t)$ is of finite-dimensional range due to our working assumptions; see condition \eqref{cond_interaction_table} below. This latter condition is  typically encountered  for fluid problems defined e.g.~on rectangular spatial domains and is made here to avoid unnecessary technicalities. 
We show that with the addition of mild dissipative assumptions on the linear operator $L_\lambda$, one can deal with more general situations and counterbalance the  loss of regularity induced by the nonlinear term $B$ to derive thereof the required a priori bounds on the SPDE solution; see Conditions~\ref{Cond_L1}-\ref{Cond_L2} and Remark~\ref{Rmk_operator_A} below. 

Equipped with these error estimates, the non-Markovian reduced SDE \eqref{Eq_reduced} provides thus a natural normal form to describe in large probability the notion of a stochastic pitchfork bifurcation for SPDEs. 
Indeed, thanks to the error estimates of Theorem~\ref{Thm_error_estimates}, the non-Markovian reduced SDEs derived from the approximation formulas of stochastic IMs given in Theorem~\ref{Thm_CM_approx}, provide an effective way to predict stochastic pitchfork bifurcations (in large probability) for a broad class of SPDEs issued from fluid dynamics, as long as the noise intensity and the eigenvalue's magnitude of the mildly unstable mode, scale accordingly.

As an application, we illustrate in Sec.~\ref{RBC_sec} below that our framework is able to predict a stochastic pitchfork bifurcation in a stochastic Rayleigh-B\'enard model, in terms of natural conditions involving the distance of the Rayleigh number to its first critical value; see Theorem \ref{RBC_thm}.  The article is then concluded by Sec.~\ref{Sec_conclusion} with a  discussion about future developments of the framework introduced here.


\section{Reduced-order dynamics of stochastic PDEs near the onset of instability} \label{Sec_2}
\subsection{SPDEs driven by multiplicative noise}\label{Sec_SEE}
In \cite{CLW15_vol1,CLW15_vol2} was undertaken the endeavor of deriving general, analytic formulas for the approximation of (local) stochastic invariant manifolds of the following class of nonlinear stochastic evolution equations posed on a Hilbert space $H$ and driven by a multiplicative noise according to
\begin{equation} \label{SEE}
\mathrm{d} u = \big( L_\lambda u + F(u) \big) \mathrm{d} t + \sigma u \circ \mathrm{d}W_t, \; u \in H,
\end{equation}
where $\sigma u \circ \mathrm{d}W_t$ means that this stochastic term (of magnitude $\sigma >0$) is considered in the Stratonovich sense.
Here $W_t$ denotes a standard Brownian motion whose probability space is denoted by $\Omega$,  and whose probability law is denoted by $\mathbb{P}$.  A noise's realization, also called a noise path will be thus labeled by an element $\omega$ in $\Omega$. This multiplicative term is such that the perturbation ``scales'' with the solution, such that the total energy is e.g.~(formally) $\mathbb{P}$-almost surely preserved in the case of a nonlinearity $F$ that is energy preserving.

Such equations results from recasting an SPDE with its boundary conditions, in which the unknown $u$ evolves in a functional space, taken to be here a Hilbert space $H$.  The operator $ L_\lambda$ represents a linear differential operator (parameterized by a scalar $\lambda$) while  $F$ is a nonlinear operator that accounts for the nonlinear terms. Both of these operators imply loss of regularity when applied to a function $u$ in $H$. To have a consistent existence theory of solutions and their stochastic invariant manifolds for SPDEs requires to take into account these loss of regularity effects. 

To do so, we place ourselves in the functional setting of sectorial operators and analytic semigroups \cite{Hen81,Pazy83}. We assume that 
\bes
L_\lambda=-A + P_\lambda, 
\ees
where  $A$ is sectorial with domain $D(A)\subset \cH$ which is compactly and densely embedded in $H$. We assume also that $-A$ is stable, while $P_\lambda$ is a low-order perturbation of $A$, i.e.~$P_\lambda: D(A^{\alpha}) \rightarrow \cH$ is a family of bounded linear operators depending continuously on $\lambda$, for some $\alpha$ in $[0,1)$; see \cite[Sec.~1.4]{Hen81} for definitions of fractional power of an operator and the function space $D(A^{\alpha})$. 
In practice, the choice of $\alpha$ should match the loss of regularity effects caused by the nonlinear terms so that 
\bes
F \colon  D(A^{\alpha})  \rightarrow H,
\ees
is a well-defined $C^p$-smooth mapping that satisfies $F(0)=0$, $DF(0)=0$ (tangency condition), and  
\be\label{F_Taylor}
F(u)= F_k(u, \cdots, u) + O(\|u\|^{k+1}_\alpha), 
\ee
with $p > k \ge 2$, and $F_k$ denoting the leading-order operator (on $D(A^{\alpha})$) in the Taylor expansion of $F$, near the origin. 
A broad class of stochastic equations that arise in many branches of physics can be recasted into this abstract functional framework.
Such equations driven by (possibly spatially inhomogeneous) multiplicative noise are met in various contexts such as in turbulence theory (intermittency phenomena) \cite{Birnir13}, ocean dynamics \cite{SFL01}, climate dynamics \cite{kug2008state,majda2009normal,CSG11,ahmed2019explaining,berner2017stochastic},  non-equilibrium phase transitions \cite{CH93,Hin00,Mun04,SSG07}, statistical physics \cite{schenzle1979multiplicative} or population dynamics \cite{HSZ02,PMax16}.

For later usage, we denote by $V$, the function space $D(A^{\alpha})$. In particular $D(A)\subset V \subset H$.    
The norm of $H$ is denoted by $\|\cdot \|$, while the norm of $V$ is denoted by  $\|\cdot \|_V.$ 
Under mild conditions often met in applications about the spectrum of $L_\lambda$ (in particular near the instability onset, \cite[Secns.~3.3 and 6.2]{CLW15_vol1}), one can then prove that for any initial condition $u_0$ in $V$ (sufficiently small) \footnote{The existence and uniqueness result proved in \cite{CLW15_vol1} for globally Lipschitz nonlinearity, applies to the case of 
a standard cutoff argument of the nonlinearity ensures existence and uniqueness for sufficiently small initial data.}, emanates  
a (unique) solution to Eq.~\eqref{SEE} in the sense that   $v = e^{-Z(t,\omega)} u$ is the unique classical solution of 
\begin{equation} \label{REE} 
\frac{\mathrm{d} v}{\mathrm{d} t} = L_\lambda v + Z(t,\omega) v + e^{-Z(t,\omega)} F(e^{Z(t,\omega)} v),  \; v(0)=e^{-Z(0,\omega)} u_0,
\end{equation}
that lies, $\mathbb{P}$-a.s., in
\bea \label{Eq_regularity}
C((0, \tau); D(A)) \cap C([0, \tau); V) \cap C^1((0,\tau); H),
\eea
where $\tau$ is a stopping time\footnote{Defined, loosely speaking, as the (random) minimum time after which $v(t,\omega)$ exits the neighborhood of the origin over which the required estimates apply.} and $Z(t,\omega)$ denotes the Ornstein-Uhlenbeck (OU) process, stationary solution of the scalar Langevin equation
\be\label{Eq_Langevin}
\d r = - r \d t + \sigma \d W_t(\omega). 
\ee
Furthermore, $u$ depends continuously on $\lambda$ and on $u_0$, and is a stochastic  $\cH$-valued adapted process. 

Note that to derive \eqref{REE} from \eqref{SEE}, we need to note that, by using the It\^o's formula, 
\beas
\d e^{-Z(t,\omega)}&=(Z(t,\omega) e^{-Z(t,\omega)} +\frac{\sigma^2}{2} e^{-Z(t,\omega)}) \d t-\sigma e^{-Z(t,\omega)} \d W_t\\
&=Z(t,\omega) e^{-Z(t,\omega)} \d t -\sigma e^{-Z(t,\omega)} \circ \d W_t,
\eeas
where the last equality follows from the conversion between It\^o and Stratonovich integrals. 
On the other hand, 
\bes
\d v =\d (e^{-Z(t,\omega)} u) =u \circ \d e^{-Z(t,\omega)}+e^{-Z(t,\omega)} \circ \d u,
\ees
and \eqref{REE} results from \eqref{SEE}, after simplification. 

Note that the usage of the OU process $Z(t,\omega)$ through the transformation $v = e^{-Z(t,\omega)} u$ is here motivated 
by the fact that the PDE with random coefficients \eqref{REE} enjoys a better temporal regularity than the SPDE \eqref{SEE} due to  the H\"older regularity of the OU process  \cite[Lemma 3.1]{CLW15_vol1}, which allows in turn for applying standard existence and uniqueness results from the theory of non-autonomous PDEs; see  \cite[Theorem 3.3.3 and Corollary 3.5.3]{Hen81}. 
We refer to \cite[Proposition 3.1]{CLW15_vol1} and  \cite[Appendix A]{CLW15_vol1}, for more details.

In this article, we are concerned with the study of bifurcations for the class of problems that can be recast into the abstract formulation of Eq.~\eqref{SEE}.
In other words, we are interested in describing how
linear instability translates to nonlinear dynamics, in presence of multiplicative noise. We are in particular concerned with the dynamical reduction problem near the onset of instability, i.e.~to describe by a low-dimensional object the SPDE dynamics in case of e.g.~the existence of a critical value $\lambda_c$ at which $m$ eigenvalues $\beta_j(\lambda)$ of $ L_\lambda$ (counting algebraic multiplicity) cross the imaginary axis.
In this case, it is known, under the assumptions recalled above, that a finite-dimensional $C^{p}$-smooth stochastic (local) invariant manifold, $\mathfrak{M}_{\lambda}(\omega)$, exists near the origin \cite[Corollary 5.1 and Prop.~6.1]{CLW15_vol1}. 
This manifold is obtained as the graph of a random mapping, $h^\lambda$ over the space $\cH_\c$ of eigenmodes $\boldsymbol{e}_j$ losing stability.  
Inspired by the deterministic theory, such manifolds are natural objects to consider for reduction of SPDEs near the instability onset. 
In \cite{CLW15_vol1,CLW15_vol2}, the approximation problem of such manifolds was thus investigated.  
Theorem \ref{Thm_CM_approx} below gives a summary of the main approximation results derived in \cite[Theorem 6.1]{CLW15_vol1} (see also 
 \cite[Corollary 6.1]{CLW15_vol1}) which allows, in practice, for deriving explicit reduced systems. 
 
 However, due to noise, intrinsic difficulties arise and the prediction of transitions for Stochastic PDEs, based on such reduced systems, requires a particular care. 
Indeed, what makes the success of center/unstable manifold theory to describe bifurcations for PDEs from low-dimensional reduced systems fails to operate in the stochastic realm, calling for a new approach to apprehend bifurcation analysis in presence of noise. The reason is that a breakdown of the reduction principle to the center/unstable manifold now takes place, and exponential slaving of the stable modes onto the unstable ones  is no longer valid, due to large excursions caused by the (white) noise, even when its intensity is small.    
We come back to this issue in Sec.~\ref{Sec_23} below and the later sections (Secns.~\ref{Sec_modeling_error} and \ref{Sec_err_estimates}) about the error estimates. Next we recall  
the approximation formulas that play a key role in the derivation of these error estimates.

\subsection{Leading-order approximation of stochastic invariant manifolds, and non-Markovian terms} \label{Sec_stoch_inv_man}

Note that $L_\lambda$ has a compact resolvent by recalling that $D(A)$ is compactly and densely embedded in $H$ \cite[Prop.~II.4.25]{EN00}. 
As a consequence, since  $L_{\lambda} \colon D(A) \rightarrow H$ is a closed operator due to the sectorial property of $-L_{\lambda}$, we have that for each $\lambda$, the  spectrum of $L_{\lambda}$, $\sigma(L_\lambda)$, consists only of isolated eigenvalues with finite algebraic multiplicities; see \cite[Thm.~III-6.29]{Kato95} (see also \cite[Corollary IV.1.19]{EN00}). This spectral property combined with the sectorial property of $-L_{\lambda}$ implies that  there are at most finitely many eigenvalues with a given real part. The sectorial property of $-L_{\lambda}$ also implies that the real part of the spectrum, $\Re \sigma(L_\lambda)$, is bounded above; see also \cite[Thm.~II.4.18]{EN00}.  

These two properties about $\Re \sigma(L_\lambda)$  allow us in turn to label elements in $\sigma(L_\lambda)$ according to the lexicographical order. According to this rearrangement, we can label the eigenvalues by a single positive integer $n$, so that 
\bea  \label{eq:ordering-1}
\sigma(L_\lambda) = \{\beta_n(\lambda) \:|\:  n \in \mathbb{N}^\ast\}, 
\eea
with, for any $1\le n < n'$, either 
\bea
\Re \beta_{n}(\lambda) > \Re \beta_{n'}(\lambda), 
\eea
or
\bea  \label{eq:ordering-3}
\Re \beta_{n}(\lambda) = \Re \beta_{n'}(\lambda), \; \text{ and } \; \Im \beta_{n}(\lambda) \geq \Im \beta_{n'}(\lambda). 
\eea
In this convention, an eigenvalue of algebraic multiplicity $m$, is repeated $m$ times. Hereafter, this rearrangement is mainly used for simplifying the notations of some theoretical developments, whereas the labeling 
with wavenumbers is often restored when dealing with applications; see Sec.~\ref{RBC_sec} below.

As already mentioned, we are concerned with describing how linear instabilities translate to the nonlinear dynamics, in presence of multiplicative noise.
To do so, the onset of instability is described in terms of the principle of exchange of stabilities (PES) \cite{MW05}, concerned 
with the loss of stability of the basic steady state.
More precisely,  the PES describes situations for which the spectrum of $L_\lambda$ experiences the following change at a critical parameter $\lambda_c$:
\begin{equation} \label{PES}
\begin{aligned}
& \Re \beta_j(\lambda)
\begin{cases} <0 & \mbox{if } \lambda < \lambda_c, \\ =0 & \mbox{if } \lambda = \lambda_c,\\ >0 & \mbox{if } \lambda > \lambda_c,
\end{cases} &&   1 \leq j \leq  m, \\
& \Re \beta_j(\lambda_c) < 0, &&  j\ge m+1,
\end{aligned}
\end{equation}
for some $m>0$, and for $\lambda$ in some neighborhood $\Lambda$ of $\lambda_c$. To this PES condition \eqref{PES} is associated the following decomposition of $\sigma(L_\lambda)$:
\bea \label{splitting}
 & \sigma_{\c}(L_\lambda) = \{\beta_j(\lambda) \:|\: j = 1, \ 2,\ \cdots, m\},\\  
 &\sigma_{\s}(L_\lambda) = \{\beta_j(\lambda) \:|\: j = m+1, \ m+2,\ \cdots\}.
\eea

The PES condition prevents eigenvalues from $\sigma_{\s}(L_\lambda)$ to cross the imaginary axis as $\lambda$ varies  in $\Lambda$. Hence, no eigenvalues other than those of $\sigma_\c(L_\lambda)$ change sign in $\Lambda$. Furthermore, the PES condition implies the following 
uniform spectral gap by reducing $\Lambda$ accordingly \cite[Lemma 6.1]{CLW15_vol1},
\be\label{general gap}
0 > 2 k \eta_{\c} > \eta_{\s}, 
\ee
where $k$ is the leading order of $F$ (see \eqref{F_Taylor}), and  
\beas \label{general etac etas}
& \eta_{\c}= \inf_{\lambda \in \Lambda} \inf_{j = 1, \cdots, m} \{\mathrm{Re} (\beta_j(\lambda))\} , & \eta_{\s}= \sup_{\lambda \in \Lambda} \sup_{j \ge m+1} \{\mathrm{Re}(\beta_j(\lambda))\}.
\eeas
Such a uniform spectral gap condition is required for the approximation formulas of stochastic invariant manifolds recalled in Theorem \ref{Thm_CM_approx} below; see \cite[Theorem 6.1]{CLW15_vol1} and \cite[Corollary 6.1]{CLW15_vol1}. It also implies (exponential) dichotomy estimates (see \cite[Eqns.~(3.24a,b,c)]{CLW15_vol1} ) to be satisfied by the semigroup generated by $L_\lambda$ and that are key to ensure the existence of stochastic invariant manifolds (see \cite[Appendix B]{CLW15_vol1}) associated with the $m$ modes losing their stability according to \eqref{PES}.

To these modes, we associate the reduced state  space, $\cH_{\c}$,  given by 
\be\label{Eq_Hc}
\cH_{\c} = \mathrm{span}\{\boldsymbol{e}_1, \cdots, \boldsymbol{e}_m\},
\ee
while a mode $\boldsymbol{e}_n$ with $n\geq m+1$ denotes a stable mode. 
In what follows, the indices $j_1,\cdots,j_k$ correspond to the $m$ critical wavenumbers losing stability (allowing repetition)  at $\lambda=\lambda_c$. 
We finally introduce the notation 
\be
(\boldsymbol{k},\boldsymbol{\beta}_c(\lambda))=\sum_{\ell=1}^{k}\beta_{j_\ell}(\lambda),
\ee
with $\boldsymbol{k}$ denoting the $k$-tuple, $(1,\cdots,1)$.   The projector onto the subspace $\cH_\s$ (resp.~$H_\c$) spanned by the stable modes (resp.~given by \eqref{Eq_Hc}) is denoted by $\Pi_\s$ (resp.~$\Pi_\c$).  The inner product in $\cH$ is denoted by $\langle \cdot, \cdot \rangle$. We have then the following stochastic invariant manifold  approximation theorem which  lies at the core of our reduction approach near the 
onset of instability for non-equilibrium systems subject to random fluctuations.

\bt \label{Thm_CM_approx}
Assume that $F$ and $L_\lambda$ satisfy the assumptions recalled above, and that the PES condition \eqref{PES} is satisfied.   Then for each $\lambda$ in a neighborhood $\Lambda$ of $\lambda_c$,  Eq.~\eqref{SEE} admits a stochastic (local) invariant manifold, $\mathfrak{M}_\omega^{\lambda}=\mbox{graph}(h_\omega^\lambda)$, with $h_\omega^\lambda$ that maps $\cH_\c$ into the regular space $D(A^\alpha)$, for every $\omega$ in $\Omega$.

The mapping $h_\omega^\lambda$ characterizing $\mathfrak{M}_\omega^{\lambda}$ is a time-dependent stochastic mapping which is approximated by  
another such mapping $\Phi_\omega^\lambda$ satisfying for any $\epsilon >0$,
\be \label{SEE App CMF formula est}
 \|h_\omega^\lambda(X,t) - \Phi_\omega^\lambda(X, t)\|_V \leq \epsilon \|X\|^k, \; \omega \in \Omega, \; t \in \mathbb{R},
\ee
as long as $X$ lies in a neighborhood $\mathcal{N}_\epsilon$ of the origin in $\cH_\c$ spanned  by the $m$ eigenmodes $\boldsymbol{e}_j$ losing their stability (see \eqref{Eq_Hc}), as lambda crosses $\lambda_c.$

The mapping $\Phi_\omega^\lambda$ in \eqref{SEE App CMF formula est} is given explicitly by the following random Lyapunov-Perron integral:
\bea \label{LP integral}
\Phi_\omega^\lambda(X, t)=\int_{-\infty}^0   e^{\sigma (k-1)W_{t,s}(\omega) \Id}e^{- s L_\lambda } P_{\s} F_k(
e^{s L_\lambda} X) \, \mathrm{d}s,
\eea
with $W_{t,s}(\omega)=W_{t+s}(\omega)-W_t(\omega).$

This mapping is well defined if $(\boldsymbol{k},\boldsymbol{\beta}_c(\lambda)) >\beta_n(\lambda)$ for $n\geq m+1$ and possesses the following expansion:
\be\label{Eq_phi}
 \Phi^{\lambda}_{\omega}(X, t) =\sum_{n\geq m+1}  \Phi^{n,\lambda}_{\omega}(X,t) \boldsymbol{e}_n, \; X\in \cH_\c,  \; \omega \in \Omega.
\ee
Here 
\be \label{Phi_n_general}
 \Phi^{n,\lambda}_{\omega}(X, t) = \sum_{(j_1, \cdots, j_k )}  F_{j_1\cdots j_k}^n M_{j_1\cdots j_k}^{n, \lambda}(t,\omega)  X_{j_1} \cdots X_{j_k}, 
\ee
 where  $1 \le j_1, \cdots, j_k \le m$, $X_j=  \langle X, e_j \rangle$, and the $F_{j_1\cdots j_k}^n$ are  coefficients accounting for the $n$-th component of the nonlinear interactions (through the leading-order term $F_k$) between the low modes $\boldsymbol{e}_{j_1}$, $\cdots$, $\boldsymbol{e}_{j_k}$ (in $\cH_c$), namely:
\be
F_{j_1\cdots j_k}^n=\langle F_k(\boldsymbol{e}_{j_1}, \cdots, \boldsymbol{e}_{j_k}), \boldsymbol{e}_n^\ast \rangle, \; 1\leq j_1, \cdots, j_k \leq m,
\ee
where the $ \boldsymbol{e}_n^\ast $ denote the eigenmodes of the adjoint operator of $L_\lambda$.  

The $M_{j_1\cdots j_k}^{n, \lambda}$-terms are 
path-dependent coefficients making explicit the $(t,\omega)$-dependence of $\Phi_\omega^\lambda$, which are obtained as stationary solutions of the auxiliary scalar SDE:
\be \label{M_eq} 
\mathrm{d} M   = \left(1 -  \Big((\boldsymbol{k},\boldsymbol{\beta}_c(\lambda)) - \beta_n(\lambda)\Big) M \right )\mathrm{d} t   - \sigma  (k-1) M \circ \mathrm{d} W_t.
\ee
\et

This theorem is a reformulation of \cite[Proposition 6.1 and Theorem 6.1]{CLW15_vol1} whose proofs can be found in \cite[Section 6.4]{CLW15_vol1}.  Note that although proved for $L_\lambda$ self-adjoint, the conclusions of these results extend to the case
$L_\lambda$ diagonalizable in $\mathbb{C}$. We emphasize  that such approximation results have been  also extended in \cite[Corollary 7.1]{CLW15_vol1} to situations in which stable modes are included in the reduced state space which has important consequences for certain applications requiring higher-dimensional reduced state spaces to ensure better approximation of the SPDE dynamics as previously illustrated in \cite[Chapter 6]{CLW15_vol2} in the context of capturing  by reduced equations noise-induced large excursions.

The above theorem provides via the coefficients of $\Phi_\omega^\lambda (X,t)$ an explicit tracking of the  time-stochastic dependence of the (unknown) random mapping $h_\omega^\lambda (X,t)$. Prior works  proposed deterministic approximation formulas for $h_\omega^\lambda (X,t)$ in the case $F(u)=B(u,u)$ \cite{blomker2010qualitative}, or  $F(u)=u^p$ (both with error $O(\|X\|)$) \cite{CDZ11}; see also  \cite[Section 4.1]{CLW15_vol2}. However, such deterministic approximation formulas are of limited virtues for applications; see e.g.~example \eqref{AB_syst}.

Instead,  the path-dependent coefficients $M_{j_1\cdots j_k}^{n, \lambda}$ (simply referred to as $M$-terms hereafter), depend on the past of the noise, when, for $\lambda$ in a neighborhood of $\lambda_c$,  the following non-resonance conditions are satisfied
\be\label{NR} \tag{NR}
\Bigl (F_{j_1\cdots j_k}^n \neq 0 \Bigr) \Longrightarrow \biggl ( \Re(\boldsymbol{k},\boldsymbol{\beta}_c(\lambda)) > \Re \beta_n(\lambda) \biggr),  \; 1\leq j_1,\cdots,j_k \leq m, \; n\geq m+1.
\ee
Under such an (NR)-condition, the $M$-terms are indeed given by 
 \be \label{Eq_Mn_formula_v2}
M_{j_1\cdots j_k}^{n, \lambda}(t,\omega) = \int_{-\infty}^0  e^{\big( (\boldsymbol{k},\boldsymbol{\beta}_c(\lambda)) - \beta_n(\lambda)\big) s + \sigma (k-1)W_{t,s}(\omega)} \mathrm{d}s,
\ee
with $W_{t,s}(\omega)=W_{t+s}(\omega)-W_t(\omega)$. Note that due to the NR-condition, the integral in \eqref{Eq_Mn_formula_v2} is well defined almost surely, since the Brownian motion $W_t$ satisfies the sublinear growth condition $\lim_{t\rightarrow \pm \infty}W_t/t = 0,$ almost surely (see \cite[Lemma 3.1]{CLW15_vol1}). Finally, note that the PES condition implies that the NR-condition
 is satisfied for $\lambda$ in  neighborhood of $\lambda_c$, as  $\sum_{\ell=1}^{k}\Re \beta_{j_\ell}(\lambda) >k \eta_\c > 2k \eta_\c >\eta_\s \geq \Re \beta_n(\lambda)$ due to  \eqref{general gap}. 
 
The exponential decaying integrand in \eqref{Eq_Mn_formula_v2} depends on the distance, $(\boldsymbol{k},\boldsymbol{\beta}_c(\lambda)) - \beta_n(\lambda)$, between the $\overbrace{(1:\cdots:1)}^{(k)}$-resonance made up from critical eigenvalues,  and the stable eigenvalue $\beta_n(\lambda)$. Such a distance to resonance controls in turns the decay of temporal correlations of the $M$-terms, provided that the noise intensity, $\sigma$, lies in some admissible range; see Lemma ~\ref{Mn recall2} in Appendix \ref{Appendix_A}. As a result, these  terms exhibit  decay of correlations of a ``reddish'' nature, although allowing for non-Gaussian statistics with heavy tails all the more pronounced than $\sigma/\sqrt{(\boldsymbol{k},\boldsymbol{\beta}_c(\lambda)) - \beta_n(\lambda)}$ is close to $1$ from below (in the case of real eigenvalues); see Appendix \ref{Sec_lognormal}.

The $M$-coefficients are thus non-Markovian and convey exogenous memory effects in the sense of  \cite{HO07, Hai09}. It is worthwhile noting that  the ``noise bath" is essential for these coefficients  to exhibit decay of correlations, the $M$-terms being reduced to simply the constant term $((\boldsymbol{k},\boldsymbol{\beta}_c(\lambda)) - \beta_n(\lambda))^{-1}$,  when $\sigma=0$, recovering standard approximations formulas of (deterministic) invariant manifolds; see \cite[Theorem 2]{CLM19_closure}.

These non-Markovian terms are new ingredients 
produced by the  interactions between the noise with the nonlinear effects. The reduced systems built from the corresponding non-Markovian parameterization,  $\Phi^{\lambda}_{\omega}(X, t)$ given by Theorem \ref{Thm_CM_approx},  are thus, in general,  non-Markovian SDEs.
In the case where $F(u)$ is a quadratic nonlinearity, $B(u,u)$, such as arising in fluid problems (see Secns.~\ref{Sec_modeling_error} and \ref{Sec_err_estimates} below), such a non-Markovian reduced SDE takes the form 
\bea \label{Eq_reduced_abstract}
\mathrm{d} X = \bigg(  \Pi_\c L_\lambda X  + \Pi_\c \Big(B(X,X) + &B(\Phi(X,t),X) +  B(X,\Phi(X,t))\\
& + B(\Phi(X,t),\Phi(X,t))  \Big) \bigg) \mathrm{d} t + \sigma X \circ \d W_t, \; X\in \mathbb{R}^m,
\eea
which in coordinate forms reads
 \bea\label{Reduced_SDE_expansion}
 \mathrm{d}X_\ell  & = \bigg (\beta_{\ell}(\lambda) X_\ell   + \overbrace{\sum_{i,j= 1}^m  B^{\ell}_{i j} X_{i}X_{j}}^{(a)}  + \overbrace{\sum_{p, i, j=1}^m \sum_{n = m+1}^\infty ( B^{\ell}_{p n}   + B^{\ell}_{n p})B^n_{i j} M_{i j}^{n,\lambda}(t,\omega) X_{i}X_{j} X_{p}}^{(b)} \\ 
&+ \underbrace{\sum_{n, n' = m+1}^\infty \sum_{i,j = 1}^m \sum_{p, q = 1}^m B^{\ell}_{n n'}B_{i j}^{n} B_{pq}^{n'} M_{i j}^{n,\lambda}(t,\omega) M_{p q }^{n',\lambda}(t,\omega)X_{i}X_{j}X_{p}X_{q}
}_{(c)} \bigg) \, \mathrm{d} t  \\
& \hspace{3cm}+ \sigma X_{\ell} \circ \mathrm{d}W_t, \quad \ 1 \leq \ell \leq m,  \;\; \text{with $B^{\ell}_{ \nu \mu} = \langle B(\boldsymbol{e}_{\nu},\boldsymbol{e}_{\mu}), \boldsymbol{e}_{\ell}^\ast \rangle$}.
\eea
Here the $X_\ell$ are aimed at approximating the amplitudes of the SPDE solution onto the $m$ modes $\boldsymbol{e}_{\ell}$ ($1 \leq \ell \leq m$) losing stability, as $\lambda$ crosses $\lambda_c$.

In the expansion \eqref{Reduced_SDE_expansion}, the $(a)$-terms account for self-interactions between the unstable modes, the  $(b)$-terms account for cross-interactions between unstable and stable modes, and the $(c)$-terms for self-interactions between the stable modes. The non-Markovian terms appear only out of the two latter types of interactions.

As we will see, the $M$-terms are key, under certain circumstances, to appropriately track the fluctuations and large-excursions  caused by the (white) noise, although the original SPDE is Markovian.  This is particularly true for
SPDEs forced by a multiplicative noise in the orthogonal subspace of e.g. its mildly unstable mode. Our parameterization formulas show that the noise gets transmitted to this mode via non-Markovian coefficients, and that the reduced equation is only stochastically driven by the latter.  This situation is a special case of more general spatially  inhomogeneous multiplicative noises, for which the remark below points out the required modifications of Theorem \ref{Thm_CM_approx}. 
\br\label{Rmk_Generalization_App_formulas}
Theorem \ref{Thm_CM_approx} above can be extended to the case of SPDEs driven by spatially  inhomogeneous multiplicative noise, of the form
\begin{equation} \label{SEE2}
\mathrm{d} u = \big( L_\lambda u + F(u) \big) \mathrm{d} t + \sigma(u) \circ \mathrm{d}W_t,
\end{equation}
where  $L_\lambda$ and $F$ satisfy the assumptions of Theorem \ref{Thm_CM_approx} and
\be \label{Eq_noise2}
 \sigma(u) = \sum_{j=1}^{N_f} \sigma_j \langle u, \boldsymbol{e}_j^\ast \rangle \boldsymbol{e}_j,
\ee
with $\sigma_j\geq 0$ and $N_f\geq 1$. While the proof of \cite[Theorem 6.1]{CLW15_vol1} requires a 
certain care to handle such multiplicative noises, the same approach applies and leads to 
approximation formulas of stochastic invariant manifolds for Eq.~\eqref{SEE2} which are still given by \eqref{Eq_phi}--\eqref{Phi_n_general}, with the only difference that the $M$-terms are now replaced by
\bea \label{Eq_Mn_formula_general}
M_{j_1\cdots j_k}^{n, \lambda}(t,\omega) = \int_{-\infty}^0  e^{\big( (\boldsymbol{k},\boldsymbol{\beta}_c(\lambda)) - \beta_n(\lambda)\big) s + (\boldsymbol{k}, \boldsymbol{\sigma}) W_{t,s}(\omega)} \mathrm{d}s, 
\eea
where 
\be
(\boldsymbol{k}, \boldsymbol{\sigma}) = \sum_{\ell=1}^k \sigma_{j_\ell} - \sigma_n.
\ee
These $M$-terms are obtained as stationary solutions of the auxiliary scalar SDE:
\bea \label{M_eq_general} 
\mathrm{d} M = \left(1 -  \Big((\boldsymbol{k},\boldsymbol{\beta}_c(\lambda)) - \beta_n(\lambda)\Big) M \right )\mathrm{d} t -   (\boldsymbol{k}, \boldsymbol{\sigma}) M \circ \mathrm{d} W_t.
\eea
\er

\br  \label{Rmk_Lyap_Exp}
It is also interesting to point out that when Theorem~\ref{Thm_CM_approx} is applied to (finite-dimensional) SDE systems of the form $\mathrm{d} y = \big( A y + F(y) \big) \mathrm{d} t + \sigma y \circ \mathrm{d}W_t$ with $A$ being a $d\times d$ matrix and $F : \mathbb{R}^d \rightarrow \mathbb{R}^d$ a sufficiently smooth function satisfying $F(0) = 0$ and $DF(0) = 0$, our approximation formula \eqref{LP integral} coincides with the one determined by the cohomological equation given by \cite[Eq.~(8.4.9)]{Arnold98} for the lowest-order non-trivial $H^s_{p0}$ therein. Indeed, for such SDEs, after performing the change of variables $x = e^{-Z(t,\omega)}y$, we obtain the following RDE system 
\be \label{Eq_RDE_x}
\frac{\d x}{\mathrm{d} t} = A x + Z(t,\omega) x + e^{-Z(t,\omega)}F(e^{Z(t,\omega)} x);
\ee 
see again~\eqref{REE}. Denoting the analogue of $\Phi$ for this transformed RDE by $\Psi$ and adopting the invariance equation approach such as reviewed in \cite[Sec.~2.2]{CLM19_closure}, the homological equation satisfied by $\Psi(t,\omega)$ reads
\be \label{Eq_homological}
\frac{\d \Psi}{\d t} + D_{x_{\c}} \Psi(x_{\c},t) (A_{\c} x_{\c}  + Z(t,\omega) x_{\c}) - (A_{\s} \Psi + Z(t,\omega) \Psi) = \Pi_{\c} e^{-Z(t,\omega)}F_k(e^{Z(t,\omega)} x_{\c}),
\ee
where $k \ge 2$ denotes the lowest-order term in the Taylor expansion of $F$. Note that compared with \cite[Eq.~(2.27)]{CLM19_closure} for the deterministic case, the above equation involves also a time-derivative term $\frac{\d \Psi}{\d t}$ arising from the time-dependence (or noise-path dependence) nature of $\Psi$. 

Note that \eqref{Eq_homological} is the same as \cite[Eq.~(8.4.9)]{Arnold98} when the RDE system \eqref{Eq_RDE_x} is considered and when $p$ therein equals the lowest order $k$ in the Taylor expansion of $F$ here, since at this lowest order the term $R_{p0}^c(x_c)$ in \cite[Eq.~(8.4.9)]{Arnold98} simply equals $\Pi_{\c} e^{-Z(t,\omega)}F_k(e^{Z(t,\omega)} x_{\c})$.
 
{\mkr One one hand}, the approach {\mkr adopted in} \cite{Arnold98} {\mkr for finite-dimensional SDEs} aims to handle more general SDE settings when {\mkr e.g.~the basic reference state is} more complicated than a steady state, {\mkr by relying on} the multiplicative ergodic theorem and the associated Lyapunov spectra and Oseledets subspaces to {\mkr perform} the decomposition of the stochastic flow. {\mkr Extension to SPDEs using the Oseledets ergodic theorem has been pursued in \cite{MZZ08}.}

{\mkr On the other hand, the practical aspects in the general setting are} still too intricate {\mkr to derive reduced systems with easily computable coefficients} even for high-dimensional SDEs (let alone SPDEs). Indeed, as pointed out in \cite[Sec.~8.4.3]{Arnold98} {\mkr already} for the application to a two-dimensional SDE system associated with the Duffing-van der Pol oscillator in the Pitchfork bifurcation scenario, ``the computational effort for these results is enormous and could only be accomplished by using the computer algebra program MAPLE ....''

{\mkr Instead, by adopting a standard normal modes framework, although we lose the characterization of perturbations near complicated states, we gain insights about (i) the derivation of explicit reduced equations based on rigorous approximation formulas (as reviewed in Theorem~\ref{Thm_CM_approx}), and (ii) error estimates (see Theorem \ref{Thm_error_estimates} below) for SPDEs that cover the challenging case of fluid problems.}

\er

\begin{figure*}[hbtp]
\centering
\includegraphics[height=0.24\textwidth, width=.9\textwidth]{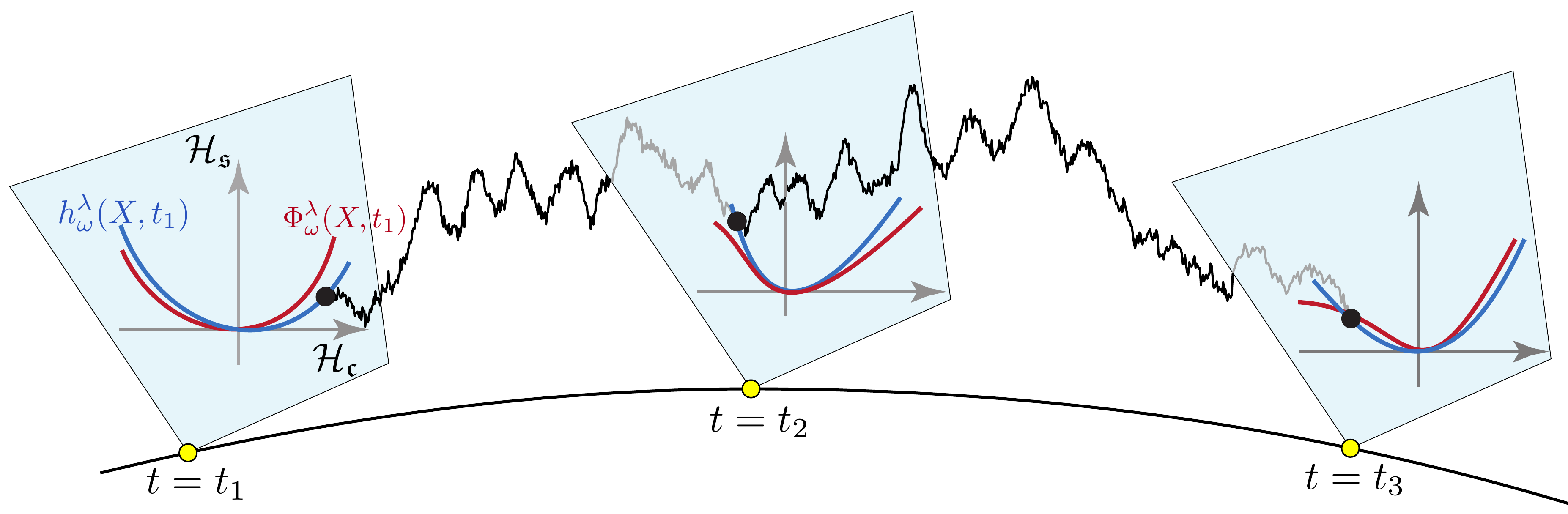}
\caption{{\bf Schematic of the stochastic invariant manifold approximation summarized by Theorem \ref{Thm_CM_approx}.}  The time-dependent stochastic manifold, $\Phi^{\lambda}_{\omega}(X,t)$ (red curve) approximates the stochastic invariant manifold $h^{\lambda}_{\omega}(X,t)$ (blue curve). For a given noise's realization $\omega$, the shape of the latter changes in the course of time. The approximation $\Phi^{\lambda}_{\omega}(X,t)$ is able to track these changes due to the path-dependent, non-Markovian, coefficients, $M_{j_1\cdots j_k}^{n, \lambda}(t,\omega)$, solving the SDE \eqref{M_eq}.  In this schematic, the black curve shows a solution path to the SPDE that evolves on the stochastic invariant manifold.}
\label{Fig_intro2}
\end{figure*}

\subsection{Tracking the large-excursions through the non-Markovian $M$-terms}\label{Sec_23}
As already mentioned, although Theorem \ref{Thm_CM_approx} provides an analytic approximation of the (unknown) slaving function $h_\omega^\lambda$, and in particular 
regarding its dependence with respect to time and the noise path, it suffers from a technical restriction, namely that the approximation  
\eqref{SEE App CMF formula est} is valid over a deterministic neighborhood  $\mathcal{N}_\epsilon$ of the origin.

From an analysis viewpoint, this theorem concedes thus some taste of dissatisfaction since as soon as $X$ becomes a time-dependent 
random variable\footnote{Such as when $X$ is a solution to non-Markovian reduced equation \eqref{Eq_reduced_abstract}.}, large excursions caused by the (white) noise are expected to take place even for small noise intensity, pushing $X$ outside of $\mathcal{N}_\epsilon$ in the course of time and questioning {\it de facto}  the validity of the approximation  
\eqref{SEE App CMF formula est}. This restriction is technical (from the proof), but even if the neighborhood would be allowed to fluctuate with time, still noise could drive the solutions outside of it.

Similarly, exponential attraction results of invariant manifolds that occur in presence of a sufficiently large spectral gap
and that extend to SPDEs  with globally Lipschitz nonlinearity (\cite[Corollary 4.3]{CLW15_vol1}), encounter however similar technical restrictions for SPDEs with nonlinearities that are only locally Lipschitz, such as for fluid problems.  The current proofs do not allow indeed for establishing a pathwise exponential attraction of local stochastic invariant manifolds, preventing 
to rely on a  reduction principle as in the deterministic theory \cite{crawford1991introduction}, to analyze bifurcations. This important but subtle point, misled some authors in their conclusions regarding the analogue of a reduction principle for SPDEs; see \cite[Remark 4.3]{CLW15_vol1}.
We provide in Secns.~\ref{Sec_modeling_error} and \ref{Sec_err_estimates}  below error estimates accounting for large-excursions and that are furthermore relevant for stochastic fluid problems, allowing in turn to make precise the use of stochastic invariant manifolds in bifurcation analysis of such SPDEs.

Before embarking in these error estimates, we provide below a simple example to illustrate these issues, and to show that the non-Markovian $M$-terms allow for tracking the large excursions caused by the noise. For this purpose, we consider the stochastic system:
\bea\label{AB_syst}
&\d A =(\lambda A+AB)\d t,  \; \lambda>0 \\
&\d B=(\kappa B-A^2) \d t + \sigma B\circ \d W_t,\; \kappa<0, \sigma>0,
\eea
in which the noise forces only the equation of the stable mode, here the $B$-equation. 
This system arises in diverse fields, with or without stochastic forces. For instance it arises in the study of nonlinear crystals \cite{schenzle1979multiplicative}, 
and more recently in the study of the finite time blow-up problem for the 3-D Navier Stokes equations  \cite{tao2016finite}, where such systems are used as elemental ``pump gates'' to execute transition of energy from one mode to another.

Here, we are looking for a closed equation  describing the evolution of $A$ without resolving $B$. If one adopts ideas of deterministic unstable/center-manifold reduction to seek for such a closure  \cite{knobloch1983bifurcations}, we arrive using this theory at the parameterization
 $B=-(2\lambda -\kappa)^{-1} A^2$, which leads to the closure
\be\label{Eq_closure1}
\dot{X}=\lambda X-(2\lambda -\kappa)^{-1} X^3.
\ee
This closure is here only deterministic and thus cannot capture the stochastic nature of the dynamics of \eqref{AB_syst}.

We would like instead to use an explicit stochastic parameterization of $B$. To be of practical interest, this stochastic parameterization should be able to track the large excursions caused  by the white noise. Application of our parameterization formulas to this particular system (see Remark \ref{Rmk_Generalization_App_formulas}), provides such a stochastic parameterization. The latter is given here by the non-Markovian parameterization $B=-M_{11}^{2,\lambda} (t,\omega)A^2$, with 
$M_{11}^{2,\lambda}$ being the stationary solution of (see \eqref{M_eq_general})
\be\label{MEq}
\d M=(1-(2\lambda-\kappa)M )\d t +\sigma M \circ \d W_t,
\ee
namely
\be\label{M_absyst}
M_{11}^{2,\lambda}(t,\omega)=\int_{-\infty}^0 e^{(2\lambda-\kappa) s -\sigma (W_{t+s}(\omega)-W_s(\omega))} \d s.
\ee

The closure becomes then the following ODE with a path-dependent coefficient:
\be\label{Eq_closure2}
\dot{X}=\lambda X-M_{11}^{2,\lambda} (t,\omega) X^3. 
\ee
Recall that the AB-system is forced by a noise in the orthogonal subspace of its mildly unstable mode. In this case,  our parameterization formulas show that the noise gets transmitted to this mode via the  non-Markovian coefficient $M_{11}^{2,\lambda}$, and that the reduced equation is only stochastically driven by the latter.  

\begin{figure}
\centering
\includegraphics[width=1\linewidth, height=.5\textwidth]{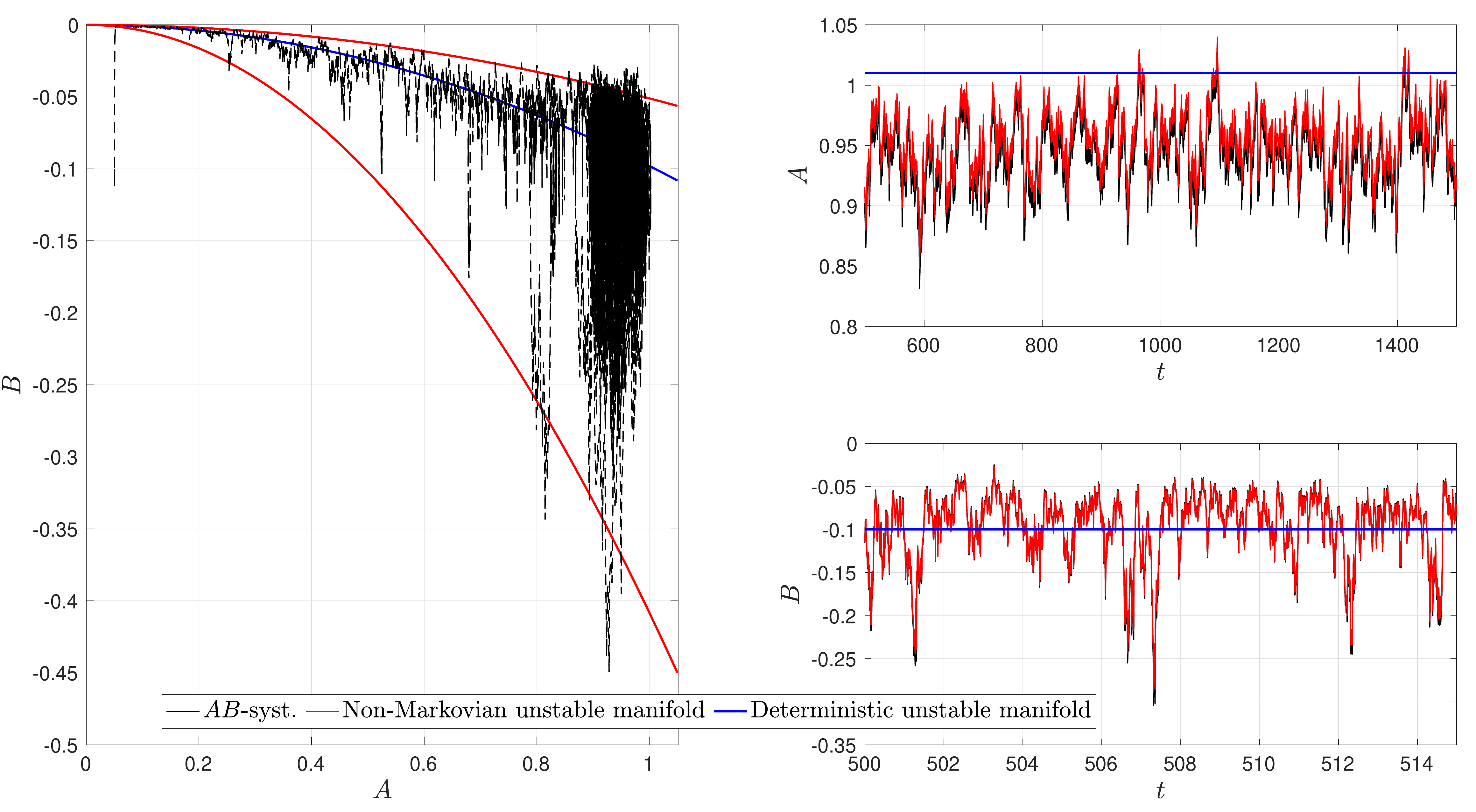}
\vspace{-1ex}
\caption{\footnotesize {\bf Left panel:} A solution path of \eqref{AB_syst} shown in the (A,B)-plane (black dash curve), and two time instances of the non-Markovian unstable manifold approximation, $B=-M_{11}^{2,\lambda} (t,\omega)A^2$ (red curves). The latter is able to track the large excursions of the solution path, while the deterministic unstable manifold gives only a poor approximation to its average motion. {\bf Right, top-to-bottom panels:} $A$- and $B$-time series as simulated from  \eqref{AB_syst} (black curves), the deterministic closure \eqref{Eq_closure1} (blue curves), and non-Markovian closure  \eqref{Eq_closure2} (red curves). }\label{Fig_AB}
\end{figure}

The practical efficiency of such parameterizations to track the large excursions induced by the white noise  is illustrated in Fig.~\ref{Fig_AB}. The left panel shows a solution path of \eqref{AB_syst} displayed in the (A,B)-plane (black dash curve), and two time instances of the non-Markovian unstable manifold approximation, $B=-M_{11}^{2,\lambda} (t,\omega)A^2$ (red curves). In contrast this figures shows that the deterministic unstable manifold (blue curve) gives only a poor approximation of the solution path's average motion. The consequence is that, given a noise realization $\omega$, the $A$- and $B$-trajectories as simulated from  \eqref{AB_syst} (black curves), are approximated to a good precision by the non-Markovian closure  \eqref{Eq_closure2} (red curves), whereas the closure \eqref{Eq_closure1} based on the deterministic unstable manifold 
is only proposing a deterministic constant approximation.

The reasons for such a successful closure by  Eq.~\eqref{Eq_closure2}, can be worked it out on this example as system  \eqref{AB_syst} possesses a natural closure which allows for providing a more intuitive ``raison d'\^etre'' of the non-Markovian parameterization $B=-M_{11}^{2,\lambda} (t,\omega)A^2$. 
Indeed, by integrating the $B$-equation of system \eqref{AB_syst}, a closed form equation can be naturally derived for describing the evolution of $B$, given by the following integro-differential equation with path-dependent kernel 
\be\label{Eq_closurex}
\frac{\d A}{\d t}=\lambda A -\bigg(\overbrace{\int_{-\infty}^t e^{\kappa (t-s) + \sigma (W_t(\omega)-W_s(\omega))} A^2(s) \d s}^{I(t)} \bigg) A.
\ee
Thus the quality of approximation of the integral term in \eqref{Eq_closurex} by the non-Markovian parameterization 
$B=-M_{11}^{2,\lambda} (t,\omega)A^2$ is conditioning the efficiency of the closure \eqref{Eq_closure2}. 

It is worthwhile noting that this non-Markovian parameterization can be obtained by integration of a backward-forward system  \cite[Chap.~4]{CLW15_vol2},
which for this example, consists of approximating $A(s)$ (for $s\leq t$) --- in the integral term of \eqref{Eq_closurex} and only for this term ---  by 
the solution of the following backward linear equation
\bea\label{Eq_ansatz}
&\frac{\d A}{\d s}^{(1)}=\lambda A^{(1)} (s), \; s<t\\
&A^{(1)}(t)=A(t).
\eea
Thus, by approximating $A(s)$ in the integral term of \eqref{Eq_closurex}, with 
\be
A(s)=e^{\lambda (s-t)} A(t) +o(A(t)), 
\ee
we obtain the following equation (up to order 3),
\bes
\frac{\d A}{\d t}=\lambda A -\bigg(\overbrace{\int_{-\infty}^t e^{(\kappa-2\lambda) (t-s) + \sigma (W_t(\omega)-W_s(\omega))} \d s}^{J(t)}\bigg)  A^3, 
\ees
from which we recover the closure \eqref{Eq_closure2} (see also \eqref{M_absyst}),  by using the change of variable $s'=s-t$ in the integral term.  The reader may wonder why approximating the exact closure \eqref{Eq_closurex} by  \eqref{Eq_closure2}, on this particular example. 
The reasons are multiple, especially in terms of insights gleaned by this operation. First, it shows the relevance of the non-Markovian parameterization resulting from application of Theorem \ref{Thm_CM_approx},  for certain parameter regimes like that corresponding of  Fig.~\ref{Fig_intro2} (for which $J(t)\approx I(t)$). Second, it shows the relevance --- for stochastic systems --- of approximating a fully non-Markovian closure like  \eqref{Eq_closurex} by a non-Markovian closure \eqref{Eq_closure2} that involves only the past of the noise, the latter allowing for an efficient simulation of the path-dependent coefficient $M_{11}^{2,\lambda} (t,\omega)$ in \eqref{M_absyst} by an SDE decoupled from the $A$-equation (i.e.~by simulating \eqref{MEq}). Instead,  the closure \eqref{Eq_closurex} would require to resolve the corresponding time integral $I(t)$ at each time step to compute $A(t)$, an operation 
that is costly numerically. This decoupling from the $A$-equation to simulate the $M$-term and thus $J(t)$ approximating $I(t)$ may be seen as a small advantage here due to the small size of system \eqref{AB_syst}, but such a feature is highly beneficial numerically when the size of the original system gets much larger, as for fluid models; see Sec.~\ref{RBC_sec} below.

As much as striking might be the ansatz  \eqref{Eq_ansatz}, the non-Markovian closure \eqref{Eq_closure2} depending only on the past of the noise, provides thus, under certain circumstances, a good approximation of $A(t)$ obtained from the closed form equation \eqref{Eq_closurex}  depending furthermore on the past-history of $A(t)$. 
The general error estimates of Secns.~\ref{Sec_modeling_error} and \ref{Sec_err_estimates} below allow us to clarify these circumstances, for this particular example and beyond. 
In particular, they show that when $\lambda$ is sufficiently small and the noise's intensity $\sigma$, scales as $\sqrt{\lambda}$, the non-Markovian closure \eqref{Eq_closure2}  (depending on the past of the noise) provides a good approximation of  the exact non-Markovian closure  \eqref{Eq_closurex} (depending not only on the past of the noise but also the past of $A$).

\section{Modeling error estimates from reduced equations}\label{Sec_modeling_error}
In this section we apply the approximation formulas of 
 Theorem \ref{Thm_CM_approx}  to derive general error estimates, establishing {\mkr further} legitimacy 
 of the reduced equations that can be derived by application of  Theorem \ref{Thm_CM_approx}, near the onset of instability. 
 
To simplify the presentation of the ideas, we focus on the following class of SPDE
\be \label{Eq_Quadratic}
\mathrm{d} u = \big(L_\lambda u  + B(u,u) \big) \mathrm{d} t + \sigma u \circ \d W_t,
 \ee
 and place ourselves in a pitchfork bifurcation scenario when $\sigma=0$.

\subsection{Assumptions and pitchfork scenario} \label{Sec_main_conditions}

 Recall that, $L_\lambda$ denotes a linear operator on the Hilbert space  $\cH$ such that $L_\lambda=-A + P_\lambda$, where  $A$ is the generator of a contraction semigroup in $\cH$, 
 while $P_\lambda$ is a low-order perturbation of $A$, i.e.~$P_\lambda: {\mkr V} \rightarrow \cH$ is a family of bounded linear operators depending continuously on the  parameter $\lambda$, with $D(A) \subset V \subset \cH$. Recall also that $D(A)$ is compactly and densely embedded in $H$.

Denote the eigenelements of $L_\lambda$ by $\{(\beta_k(\lambda), \boldsymbol{e}_k) \; : \; k \in \mathbb{N}\}$. 
To derive our error estimates, we assume

\bh \label{Cond_A1}
The linear operator $L_\lambda$ is self-adjoint.
\eh

\bh \label{Cond_A2}
The principle of exchange of stabilities (PES) condition is satisfied with $m=1$ in \eqref{splitting}, namely there exists $\lambda_c$ such that:
\begin{equation} \label{PES_pitch}
\begin{aligned}
& \beta_1(\lambda)
\begin{cases} <0 & \mbox{if } \lambda < \lambda_c, \\ =0 & \mbox{if } \lambda = \lambda_c,\\ >0 & \mbox{if } \lambda > \lambda_c,
\end{cases} &&  \\
& \beta_j(\lambda_c) < 0, \quad \textrm{for all } j\ge 2,
\end{aligned}
\end{equation}
and all the eigenvalues $\beta_j(\lambda)$, for $j\geq 2$, remain negative for $\lambda$ in some interval $[\lambda_c, \lambda^*]$, with $\lambda^*>\lambda_c$.
\eh

\bh \label{Cond_A3}
The nonlinear term $B$ is a continuous bilinear map from $V \times V$ to $H$. 
Furthermore the following conditions hold: 
\be  \label{Eq_energy_conservation} 
\langle B(u,u), u \rangle = 0, \;  \textrm{ for any } u  \textrm{ in }  V,
 \ee
and there exists $N\geq 1$ for which  
\be \label{cond_interaction_table}
\langle B(\boldsymbol{e}_1,\boldsymbol{e}_1), \boldsymbol{e}_n \rangle = 0, \; n > N, 
\ee
while for at least one $2 \leq n \leq N$, $\langle B(\boldsymbol{e}_1,\boldsymbol{e}_1), \boldsymbol{e}_n \rangle \neq 0.$ 
\eh

Note that if $\langle B(\boldsymbol{e}_1,\boldsymbol{e}_1), \boldsymbol{e}_n \rangle = 0$ for all $n \ge 2$, then the manifold $\Phi_\lambda$ will be identically zero and one needs then to consider higher-order approximation formulas; see \cite{CLW15_vol2}.

\subsection{Preparatory estimates for the pitchfork scenario} \label{Sec_apriori_bounds}

 The parameterization of the modes $\boldsymbol{e}_2$ to $\boldsymbol{e}_N$ is given by 
\be \label{Eq_Phi_Burgers}
\Phi(X,t,\omega) = \sum_{n=2}^N \Phi_{n}(X,t, \omega)    \boldsymbol{e}_n, 
\ee
with 
\be \label{Eq_Phi_Burgers_comp}
 \Phi_{n}(X, t, \omega) = B_{11}^n M_{11}^{n,\lambda} (t,\omega) X^2,
\ee
where 
\be
B_{ij}^n = \langle B(\boldsymbol{e}_i,\boldsymbol{e}_j), \boldsymbol{e}_n \rangle.
\ee
Here $X$ denotes the amplitude of $\boldsymbol{e}_1$.  In what follows we will omit the dependence on $\omega$, reserving this explicit dependence only when necessary.

Due to  \eqref{Eq_energy_conservation}, we have $\langle B(\boldsymbol{e}_1,\boldsymbol{e}_1), \boldsymbol{e}_1 \rangle = 0$ and the invariant manifold (IM) reduced equation,  up to order $O(X^3)$, takes the following form 
\be \label{Eq_reduced}
\mathrm{d} X = \Big(\epsilon X  + X \mathcal{F}_{\c}(\Phi(X,t))  \Big) \mathrm{d} t + \sigma X \circ \d W_t,
\ee
with 
\be \label{Def_Fc}
 \mathcal{F}_{\c}(v)=\left \langle \big(\hspace{-.2ex} B(\boldsymbol{e}_1, v) \hspace{-.3ex} +\hspace{-.3ex}  B(v,\boldsymbol{e}_1)\big), \boldsymbol{e}_1 \right \rangle,  \quad \forall v \in V_\s.
\ee   
Note that 
\bea\label{Cubic_term}
X\mathcal{F}_{\c}(\Phi(X,t)) &=  X\bigg(\sum_{n=2}^N \Phi_n(X,t) \big( B_{1n}^1 + B_{n1}^1\big)\bigg) \\
& =\bigg(\sum_{n=2}^N B_{11}^n M_{11}^{n,\lambda} (t) \big( B_{1n}^1 + B_{n1}^1\big)\bigg)  X^3 .
\eea

We assume then
\bhp  \label{Cond_P}
The coefficients in \eqref{Cubic_term}, satisfy for $\lambda$ in $[\lambda_c, \lambda^*]$, the dissipation condition 
\be \label{Cond_supercritical}
\sum_{n=2}^N B_{11}^n M_{11}^{n,\lambda} (t)  \big( B_{1n}^1 + B_{n1}^1\big) < 0, \; \text{ for all } t.
\ee
\ehp
 Condition \ref{Cond_P} is a sufficient condition for the reduced equation \eqref{Eq_reduced} to experience a stochastic supercritical pitchfork bifurcation as $\epsilon$ changes sign; see Sec.~\ref{Sec_stoch_pitch} below and Appendix \ref{Sec_1DSDE_pitchfork}.
 
Since $M_{11}^{n,\lambda}$ is always positive, \eqref{Cond_supercritical} can be ensured if we have  for all $2 \leq n\leq N$,
\be \label{Cond_supercritical2}
B_{11}^n \big( B_{1n}^1 + B_{n1}^1\big) \le 0,  
\ee
and at least one of the inequalities in \eqref{Cond_supercritical2} is a strict one. 
In case $\sigma=0$, such a condition is a sufficient condition for the deterministic reduced equation Eq.~\eqref{Eq_reduced} to experience a pitchfork bifurcation \cite{crawford1991introduction,MW05}. 

We turn now to the derivation of an equation satisfied by $\Phi$ which is key in our estimates presented below. 
From \eqref{Eq_Phi_Burgers_comp}, we obtain by using the Stratonovich form of
It\^o's formula \cite{Kunita90}:
\be \label{Eq_dPhi}
\d \Phi_{n} = B_{11}^n X^2 \circ \d M_{11}^{n,\lambda} + 2 B_{11}^n M_{11}^{n,\lambda} X \circ  \d X. 
\ee
Since here $(\boldsymbol{k},\beta_c(\lambda))=2 \epsilon$ and $k=2$ in \eqref{M_eq}, we have that  $M_{11}^{n,\lambda}$ satisfies
\be \label{M_eq_in_Burgers}
\mathrm{d} M  = \left(1 -  \big(2 \epsilon - \beta_n(\lambda) \big) M \right )\mathrm{d} t  - \sigma M \circ \mathrm{d} W_t.
\ee
Using \eqref{M_eq_in_Burgers} and \eqref{Eq_reduced} in \eqref{Eq_dPhi}, we get after simplification that 
\bea \label{Eq_dPhi2}
\d \Phi_{n} = \big( \beta_n(\lambda) \Phi_{n} + B_{11}^n X^2 & + 2\Phi_{n} \mathcal{F}_c(\Phi(X,t)) \big) \d t + \sigma \Phi_{n} \circ \d W_t.
\eea
Then, by \eqref{Eq_Phi_Burgers} and \eqref{Eq_dPhi2}, the equation satisfied by $\Phi$ is:
\bea \label{Eq_dPhi3}
\d \Phi = \Big( L^{\s}_\lambda \Phi + X^2 B_{11}^{\s} + 2\Phi \mathcal{F}_{\c}(&\Phi(X,t))  \Big) \d t + \sigma \Phi \circ \d W_t,
\eea
where $L^{\s}_\lambda = \Pi_{\s} L_\lambda$ and we have used
\be \label{Cond_finite_interaction}
\sum_{n=2}^N B_{11}^n X^2 \boldsymbol{e}_n = X^2 \Pi_{\s} B(\boldsymbol{e}_1,\boldsymbol{e}_1) ,
\ee
which holds thanks to \eqref{cond_interaction_table}, and we have introduced 
\be \label{Def_B_11s}
B_{11}^{\s}  = \Pi_{\s} B(\boldsymbol{e}_1,\boldsymbol{e}_1).
\ee

To summarize, we get from \eqref{Eq_reduced} and \eqref{Eq_dPhi3} that $(X, \Phi)$ satisfy the following system:
\bea \label{Eq_X_Phi}
\mathrm{d} X &= \Big(\epsilon X  + X \mathcal{F}_{\c}(\Phi(X,t))  \Big) \mathrm{d} t + \sigma X \circ \d W_t, \\
\d \Phi & = \Big( L^{\s}_\lambda \Phi + X^2 B_{11}^{\s} + 2\Phi  \mathcal{F}_{\c}(\Phi(X,t)) \Big) \d t + \sigma \Phi \circ \d W_t.
\eea

We have then the following lemma, whose proof is provided in Appendix \ref{Sec_Apriori_est}.

\bl[A priori estimates for $X$ and $\Phi$]  \label{Lem_apriori_est}
Assume that Condition \ref{Cond_P} holds. 
Consider the system \eqref{Eq_X_Phi} for which the initial condition $(X(0), \Phi(0,\omega))$
satisfies, for any $\omega$ in $\Omega$, 
\bea \label{Eq_Phi_IC}
X(0)&=X_0,\\
\Phi(0,\omega)&=\sum_{n=2}^N \Big(B_{11}^n M_{11}^{n,\lambda}(0,\omega)X_0^2\Big)  \boldsymbol{e}_n,
\eea
where $N$ is the integer that appears in \eqref{cond_interaction_table}. 

Assume furthermore that Conditions \ref{Cond_A1}, \ref{Cond_A2}, and \ref{Cond_A3} hold, with $\beta_1(\lambda) = \epsilon$ in \ref{Cond_A2}. Assume finally that $\sigma = \sqrt{\epsilon}$ and that 
\be \label{Eq_X_IC}
|X_0| \sim \sqrt{\epsilon}. 
\ee
 Then, for any given $T>0$ and any $\chi$ in $(0,1)$, there exists a constant $C>0$ depending on $\chi$ but independent of $\epsilon$, for which the following estimate holds:   
\be \label{A_priori_goal}
 \mathbb{P}\left( \sup_{[0, T/\epsilon]} |X(t)| \leq C \sqrt{\epsilon}  \;\; \textrm{and} \sup_{[0, T/\epsilon]} \|\Phi(t) \| \leq C \, \epsilon \right) \geq 1- \chi.
\ee
\el

\needspace{2ex}
\br\label{Rem_Apriori} 
\hspace*{.1em} 
\bi
\item[(i)] The proof of Lemma \ref{Lem_apriori_est} given in Appendix \ref{Sec_Apriori_est} relies on Lemmas \ref{Lemma_Mn_bound} and \ref{Lemma_BM_bound}, themselves proved in their respective appendices. These latter lemmas show that the constant $C$ in \eqref{A_priori_goal} grows as $\chi$ decreases, i.e.~the larger the probability with which  \eqref{A_priori_goal} holds,  the closer one needs to push $\epsilon$ towards $0$ in order for the corresponding a priori bounds about $ |X(t)|$ and  $\|\Phi(t) \|$ to be below a target value.

\item[(ii)] Note that a solution $(X(t),\Phi(t))$ to the system \eqref{Eq_X_Phi} emanating from an initial condition satisfying \eqref{Eq_Phi_IC}, implies  that 
$\Phi(X,t)$ is given by the parameterization \eqref{Eq_Phi_Burgers}, for all $t>0$. As it can be observed in the proof of Lemma \ref{Lem_apriori_est}, it is the dynamic interpretation of $\Phi(X,t)$ (i.e.~the equation it solves in  system \eqref{Eq_X_Phi}) that allows for deriving estimates about $\Phi(X,t)$, while $X$ is a solution to the reduced equation \eqref{Eq_reduced}.
\ei
\er

\subsection{Residual error estimates in large probability} \label{Sec_residual_estimates}

We first rewrite the original SPDE \eqref{Eq_Quadratic} into a coupled system for $u_{\c}(t) = \Pi_{\c} u(t) = x(t) \boldsymbol{e}_1$ and $u_{\s} = \Pi_{\s} u$: 
\begin{subequations} \label{Eq_Quadratic_v2}
\begin{align}
&\hspace{-.5ex}\mathrm{d} x  = \big(\epsilon x  + x \mathcal{F}_{\c}(u_\s) \hspace{-.1ex}+ \hspace{-.1ex} \Pi_{1}B(u_\s, u_\s)  \big) \mathrm{d} t \hspace{-.1ex}+\hspace{-.1ex} \sigma x \!\circ\! \d W_t,  \label{Eq_Quadratic_v2_a}\\
&\hspace{-.5ex}\mathrm{d} u_{\s} \! =\! \big(L_\lambda^\s u_{\s}  + x^2 B_{11}^{\s} + x\mathcal{F}_\s(u_\s)+\Pi_\s B(u_\s, u_\s)  \big) \mathrm{d} t   + \sigma u_{\s} \circ \d W_t, \label{Eq_Quadratic_v2_b}
\end{align}
\end{subequations}
where $F_{\c}$ is defined in \eqref{Def_Fc} and  
\bea\label{Eq_Fs}
& \Pi_{1}B(u_\s, u_\s)=\langle B(u_\s, u_\s), \boldsymbol{e}_1\rangle, \\
&\mathcal{F}_\s(v)=\Pi_\s (B(\boldsymbol{e}_1,v) +B(v,\boldsymbol{e}_1)), \quad \forall v  \in V_{\s}.
\eea
Note that the term $\Pi_{1} \big( x^2 B( \boldsymbol{e}_1, \boldsymbol{e}_1) \big)$ does not appear in the equation for $x$  since it equals 
$x^2 \langle B( \boldsymbol{e}_1, \boldsymbol{e}_1), \boldsymbol{e}_1 \rangle$, which is zero thanks to \eqref{Eq_energy_conservation}.

For any arbitrary stochastic processes $X(t)\boldsymbol{e}_1$ and $Y(t)$ that are adapted to the underlying stochastic basis and evolve respectively in $\cH_{\c}$ and  in $V_\s$,
we define the residual $\mathcal{R}(t) = (\mathcal{R}_1(t),  \mathcal{R}_2(t))$, where
\bea \label{Def_residual}
\mathcal{R}_1(t) & = X(t) - e^{\epsilon t + \sigma W_t} X(0) - \int_0^t  e^{\epsilon (t-s) + \sigma (W_t - W_s)}  E_1(s) \d s, \\
\mathcal{R}_2(t) & = Y(t) - e^{t L_{\lambda}^{\s} + \sigma W_t} Y(0) - \int_0^t  e^{L_\lambda^{\s} (t-s) + \sigma (W_t - W_s)} E_2(s) \d s, \; \text{with}\\
 E_1(s) &= X(s)  \mathcal{F}_\c (Y(s))+  \Pi_1 B(Y(s), Y(s)), \\
 E_2(s) &= X^2(s) B_{11}^\s + X(s) \mathcal{F}_\s (Y(s))+ \Pi_\s B(Y(s), Y(s)).
\eea
This residual measures the modeling error made by the process $(X(t),Y(t))$  in satisfying the system \eqref{Eq_Quadratic_v2}, and thus the modeling error made by $X(t) \boldsymbol{e}_1+Y(t)$ in satisfying the SPDE \eqref{Eq_Quadratic}.

Now pick up any solution $(X(t), \Phi(t))$ of the surrogate system \eqref{Eq_X_Phi}, for which the initial condition satisfies the relation \eqref{Eq_Phi_Burgers_comp} at $t=0$. Note that such a solution is well defined for all time $t\ge 0$, since under this constraint on the initial data, the system \eqref{Eq_X_Phi} is equivalent to the reduced equation \eqref{Eq_reduced}. 

We have the following bounds for the residual: 

\bt \label{Thm_residual}
Consider the solution  $(X(t), \Phi(t))$ to Eq.~\eqref{Eq_X_Phi} and its residual error defined in \eqref{Def_residual}. Assume that the assumptions of Lemma~\ref{Lem_apriori_est} are satisfied. Then, for any $T>0$ and any $\chi$ in $(0,1)$, there exists a constant $C>0$ independent of $\epsilon$ for which 
\be  \label{Eq_residual_goal1}
\mathbb{P}\left( \sup_{t \in [0, T/\epsilon]} |\mathcal{R}_1(t)| \leq C \epsilon \right) \geq 1-\chi,
\ee
and
\be  \label{Eq_residual_goal2}
\mathbb{P}\left( \sup_{t \in [0, T/\epsilon]}  \|\mathcal{R}_2(t) \| \leq  C \epsilon^{3/2}\right) \geq 1-\chi.
\ee
\et

\bp
Let  $(X(t), \Phi(t))$ be a solution to \eqref{Eq_X_Phi}. After integration and simplifications, we obtain that \eqref{Def_residual} becomes
 \bea
& \mathcal{R}_1(t)  = - \int_0^t  e^{\epsilon (t-s) + \sigma (W_t - W_s)}  \Pi_{1} B(\Phi(s), \Phi(s)) \d s, \\ 
& \mathcal{R}_2(t)  =  \int_0^t   e^{L_\lambda^{\s} (t-s) + \sigma (W_t - W_s)}   E_3(s) \d s,
 \eea
 with 
 \be \label{Eq_E3}
 E_3 = 2\Phi  \mathcal{F}_\c (\Phi)  -  X \mathcal{F}_\s (\Phi)   - \Pi_{s} B(\Phi, \Phi).
 \ee
Then, the desired results follow directly from the a priori estimates \eqref{A_priori_goal} together with the probabilistic estimate about the Brownian motion given in Lemma~\ref{Lemma_BM_bound}. Indeed, by letting $\Omega^*$ and $\gamma$ be the same as given in the proof of Lemma~\ref{Lem_apriori_est} (see \eqref{Eq_apriori_est2}), we have proved that $\|\Phi(t,\omega)\| \le C \epsilon$ for $t$ in $[0, T/\epsilon]$ and $\omega$ in $\Omega^*$; cf.~\eqref{Eq_Phi_est6}. Since $\Phi(t,\omega)$ as defined in \eqref{Eq_Phi_Burgers} takes value in a finite dimensional subspace of $V$ due to the assumption \eqref{cond_interaction_table}, the two norms $\|\cdot\|_V$ and $\|\cdot\|$ are equivalent when restricted to this subspace. We get then $\|\Phi(t,\omega)\|_V \le C \epsilon$ by redefining $C$. Thus,
\bea
|\mathcal{R}_1(t,\omega)| & \le C_B \int_0^t  e^{\epsilon (t-s) + 2 \gamma}  \|\Phi(s,\omega)\|_V^2 \d s  \\
& \le C_B C^2  e^{2 \gamma}  \int_0^t  e^{\epsilon (t-s)} \d s \, \epsilon^2 \\
& = C_B C^2  e^{2 \gamma}   (e^{\epsilon t} - 1) \epsilon \\
& \le C_B C^2  e^{2 \gamma}   (e^T - 1) \epsilon,   \quad \forall t \in [0, T/\epsilon], \omega \in \Omega^*.
\eea
The estimate \eqref{Eq_residual_goal1} follows. 

Finally, recall  that for such $t$ and $\omega$, we have that $|X(t,\omega)| \le C\sqrt{\epsilon}$ for some $C >0$ independent of $\epsilon$; see \eqref{Eq_apriori_est3}. Using the above estimates for $X$ and $\Phi$, we get then
\bea
\|\mathcal{R}_2(t)\|  & \le  \int_0^t   e^{\beta_2(\lambda) (t-s) + 2 \gamma}   E_3(s) \d s \\
& \le \widetilde{C} (\epsilon^{3/2} + \epsilon^2),  \quad \forall t \in [0, T/\epsilon], \omega \in \Omega^*,
\eea
where $\widetilde{C}>0$ is another constant independent of $\epsilon$. We have thus derived the estimate \eqref{Eq_residual_goal2}.
\ep

\br\label{Rmk_inhom}
When more general spatially inhomogeneous multiplicative noise such as given in Remark \ref{Rmk_Generalization_App_formulas} is considered, we can still derive residual estimates in the style of \eqref{Eq_residual_goal1} and \eqref{Eq_residual_goal2} under the same conditions as required in Theorem~\ref{Thm_residual} with $\sigma_n  = c_n \sqrt{\epsilon}$ for some $c_n \ge 0$, $n = 1, \ldots, N$. The residual $\mathcal{R}_2$ needs to be adapted accordingly by projecting onto each stable mode $n=2,\cdots, N$ as well. 

The derivation of the associated residual estimates follows the same lines of arguments presented in the proof of Theorem~\ref{Thm_residual} once the a priori estimate for $X(t)$ and $\Phi(t)$ given by \eqref{A_priori_goal} is established. The derivation of this latter a priori estimate follows the same steps as presented in Appendix \ref{Sec_Apriori_est}, but by relying on the following transformation rather than the one given in \eqref{Eq_random_transform}: 
\be
U = e^{-2 \epsilon t - 2 \sigma_1 W_t} X^2, \quad \Psi_n = e^{-2 \beta_n t - 2 \sigma_n W_t} \Phi_n^2.
\ee
\er

\section{Low- and high-mode error estimates} \label{Sec_err_estimates}
In this section, we go beyond the residual error estimates presented above to derive error estimates between the solutions of the finite-dimensional surrogate system \eqref{Eq_X_Phi} and those of the original SPDE \eqref{Eq_Quadratic}. To do so, it requires a priori bounds for solutions to the SPDE \eqref{Eq_Quadratic}. The main difficulty compared to those conducted for solutions of \eqref{Eq_X_Phi} (cf.~Lemma \ref{Lem_apriori_est}) lies in the infinite-dimensional nature of $u_\s$, whereas $\Phi$ is of finite-dimensional range due to \eqref{cond_interaction_table}. As a consequence, the loss of regularity via the nonlinearity $B$ is now effective which requires new estimates. 

To cope with these difficulties one needs extra assumptions on the linear part $L_\lambda$. To the previous assumptions on this operator, we 
require in particular an extra condition  to be satisfied  (Condition~\ref{Cond_L2}) for the linear dissipative effects to suitably counterbalance the loss of regularity due to the nonlinear terms in $B$. These conditions are described below. 

\bhl \label{Cond_L1}
Let $L_\lambda = -A + P_\lambda$ satisfy the Conditions \ref{Cond_A1} and \ref{Cond_A2} of Sec.~\ref{Sec_modeling_error}, with  $A$ given by 
\be\label{assumption_L}
A =  -\Delta \mbox{diag}(\boldsymbol{\nu}) 
\ee
where $\boldsymbol{\nu}$ is a $d$-dimensional vector with positive entries, denoted by $\nu_j$, and $\mbox{diag}(\boldsymbol{\nu})$ denotes the $d\times d$ diagonal matrix with $\boldsymbol{\nu}$ as the diagonal entries. 
We assume furthermore that the boundary conditions and functional setting are such that 
for any $u$ in $D(A)$,  $\langle A u, u \rangle \geq 0$ and $\sqrt{\langle A u, u \rangle}$ defines (after performing integration by parts) a norm $ \|u\|_{\bm{\nu}}$ on $V$ equivalent to $\| u \|_{V}$.  
\ehl

\bhl \label{Cond_L2}
There exist $\eta$ in $(0,1)$ and $\delta > 0$ such that for all $\lambda$ in the interval $[\lambda_c, \lambda^*]$ given in Condition~\ref{Cond_A2}, the following damping relations hold for all $u_{\s}$ in $\mathcal{W}_\s=D(A)\cap H_\s $: 
\begin{subequations} \label{Eq_P_control}
\begin{align}
& \langle -\eta A u_{\s} +  P_\lambda u_{\s}, u_{\s} \rangle \le - \delta \|u_{\s}\|^2, \label{Eq_P_control_a}  \\
&  \langle - \eta A u_{\s} + P_\lambda u_{\s}, A u_{\s} \rangle \le - \delta \| u_{\s}\|^2_{\boldsymbol{\nu}}.  \label{Eq_P_control_b}  
\end{align}
\end{subequations}
\ehl

For this class of linear operators, the error estimates derived hereafter allow us then, from the reduced equations based on Theorem \ref{Thm_CM_approx}, to conclude to a pitchfork bifurcation scenario in large probability. 
Our results apply to a broad class of SPDEs such as given by \eqref{Eq_Quadratic} (allowing even for inhomogeneous noise), as long as the noise intensity $\sigma$ and the eigenvalue's magnitude of the mildly unstable mode, scale accordingly, while the  linear and nonlinear terms obey standard (energy-preserving)  assumptions encountered in fluid problems  (Condition~\ref{Cond_A3}).  Applications to Rayleigh-B\'enard convection are detailed in Sec.~\ref{RBC_sec} below.

\br \label{Rmk_operator_A}
 The form of $A$ assumed in Condition~\ref{Cond_L1} is chosen to simplify the derivation of the error estimates in Theorem~\ref{Thm_error_estimates}. 
More general differential operators $A$ than \eqref{assumption_L} could have been considered for which the conclusions of Theorem~\ref{Thm_error_estimates} would still hold. In that respect,  an analogue of Conditions~\ref{Cond_L1} and \ref{Cond_L2} would still hold for a differential operator $A$ in divergence form, that (possibly) includes higher-order derivatives and that 
satisfies a strong ellipticity condition (see \cite{Tem97} and \cite[Sec.~7.2]{Pazy83}).  
Given a positive integer $m$, such an operator writes for any $u$ in the Sobolev space $(H^{2m}(\mathcal{D}))^d$ (over a smooth bounded domain $\mathcal{D}$)
\be \label{Def_A_general}
A u = \sum_{\substack{\alpha, \beta \in \mathbb{N}^p \\ [\alpha] = [\beta] = m}} (-1)^m D^\beta (a_{\alpha \beta}(x) D^\alpha u), 
\ee
where $\alpha = (\alpha_1, \ldots, \alpha_p) \in \mathbb{N}^p$, $[\alpha] = \alpha_1 + \cdots + \alpha_p$, the coefficients $a_{\alpha \beta}$ are matrix-valued mappings which are sufficiently smooth into the space of $d \times d$ symmetric matrices, and $a_{\alpha \beta} = a_{\beta \alpha}$.  If the lower  bound of $\langle Au, u\rangle$ is controlled as follows (for any $u$ in $D(A)$)
\be
\langle A u , u \rangle \ge \kappa \sum_{\alpha \in \mathbb{N}^p, [\alpha] =  m} \|D^\alpha u \|^2, \; \text{ for some } \; \kappa >0,
\ee
then with $\|u\|_{V} =\sum_{\alpha \in \mathbb{N}^p, [\alpha] =  m} \|D^\alpha u \|^2$, analogues of Conditions~\ref{Cond_L1} and \ref{Cond_L2} could be formulated for which 
the conclusions of  Theorem~\ref{Thm_error_estimates} would still hold.  In Sec.~\ref{RBC_sec} below, it is shown that Conditions~\ref{Cond_L1} and \ref{Cond_L2} as formulated here, i.e.~for the operator $A$ defined by \eqref{assumption_L}, and with the linear damping estimates  
\eqref{Eq_P_control_a}-\eqref{Eq_P_control_b}, are sufficient to apply Theorem~\ref{Thm_error_estimates} to a standard Rayleigh-B\'enard convection problem subject to multiplicative noise. 
\er

\subsection{A priori bounds on SPDE solutions}\label{Sec_dicuss_bounds}
Roughly speaking, we consider the following a priori bounds on the SDPE solution $(x(t),u_\s(t))$
\bea\label{Goal_123}
\text{(I)}:\; \sup_{t \in [0, T/\epsilon]} |x(t)| \le C \sqrt{\epsilon},\\
\textrm{(II)}:\; \sup_{t \in [0, T/\epsilon]} \|u_{\s}(t) \| \le C \epsilon,\\
\textrm{(III)}:\; \sup_{t \in [0, T/\epsilon]} \|u_{\s}(t) \|_V \le C \epsilon,
\eea
and are concerned whether these estimates hold in large probability. These a priori bounds are key to derive our main Theorem about error estimates; see Theorem \ref{Thm_error_estimates} and its proof.  We recall that the SPDE solutions are here understood in the sense recalled in Sec.~\ref{Sec_SEE}.  In that respect, to derive a priori bounds like \eqref{Goal_123}, one first work with the transformed equation Eq.~\eqref{REE}, followed by the inverse transformation to infer back the desired bounds about the SPDE solutions. We refer to Appendix \ref{Sec_Apriori_est_SPDE} for more details and provide here the main elements and ideas to derive  \text{(I)}, \text{(II)}, and \text{(III)} in large probability. 

In that perspective, observe that standard energy estimates on the solutions to \eqref{Eq_Quadratic_v2}  lead to 
\be
x^2(t,\omega) + \|u_{\s}(t,\omega)\|^2 \le e^{2\epsilon t + 2 \sigma W_t(\omega)}\|u(0,\omega)\|^2.
\ee
Thus, if one assumes $|u(0)| \sim \sqrt{\epsilon} $ one would obtain that, in large probability, $|x(t)| \le C \sqrt{\epsilon}$ and $\|u_{\s}(t) \| \le C \sqrt{\epsilon}$  by conducting similar estimates as in Step 1 of the proof of Lemma \ref{Lem_apriori_est} (and using \eqref{Eq_energy_conservation} for controlling the nonlinear term); see Appendix \ref{Sec_Apriori_est}.

To derive error estimates between $u_\s(t)$ and its parameterization $\Phi(t)$, this is however insufficient, as the bound on $\|u_{\s}(t) \| $ should scale as $\epsilon$ and not $\sqrt{\epsilon}$ since $\|\Phi (t)\|$ scales like $\epsilon$ due  \eqref{A_priori_goal}. 
The main difficulty  compared to Step 2 of the proof of Lemma \ref{Lem_apriori_est}, lies in the infinite-dimensional nature of $u_\s$, whereas $\Phi$ is of finite-dimensional range due to \eqref{cond_interaction_table}. As a consequence, the loss of regularity 
via the nonlinearity $B$ is now effective which requires new estimates. 

For instance, when estimating $\|u_\s\|$, by taking the inner product of  the $u_\s$-equation with  $u_\s$ in \eqref{Eq_Quadratic_v2}, 
this loss of regularity causes for instance the appearance in the RHS of terms depending on $ \|u_{\s}\|_V$ (and not only on $\|u_\s\|$) since e.g.~$\langle B(\boldsymbol{e}_1,  u_{\s}), u_{\s} \rangle \leq C_B  \|\boldsymbol{e}_1\|_V \|u_{\s}\|_V \|u_{\s}\|$.

Nevertheless, such difficulties can be handled through a careful exploitation of the damping effects brought by the linear part  
allowing for bypassing to estimate $ \|u_{\s}\|_V$ {\mkr (at this stage)}, and leading us to prove that, in large probability, and for $t$ in  $[0, T/\epsilon]$,
\be\label{Key_est1}
\|u_{\s}(t) \|^2\leq e^{2\gamma}  \|u_{\s}(0)\|^2 + \frac{ C e^{4\gamma}}{\delta}\epsilon^2,
\ee 
where $\gamma$ is as given in Lemma \ref{Lemma_BM_bound},  $\delta$ is given in \eqref{Eq_P_control}, and $C>0$ is a generic constant independent of $\epsilon$. The proof of this key estimate is given in Appendix~\ref{Sec_Apriori_est_SPDE_key1}.

As a consequence, if $|x(0)| \sim \sqrt{\epsilon}$ and $ \|u_{\s}(0)\| \sim \epsilon,$ one obtains that for any $T>0$ and any $\chi$ in $(0,1)$, 
\be
\mathbb{P}\Big(\textrm{(I) and (II)} \Big) \geq 1-\chi. 
\ee

{\mkr Nevertheless}, since $B: V \times V \rightarrow \cH$ (loss of regularity), to derive a full set of a priori estimates about the SPDE solution, {\mkr one cannot escape for estimating} $ \|u_{\s}\|_V$. {\mkr To do so, we} take the inner product of the $u_\s$-equation with $A u_\s$ in \eqref{Eq_Quadratic_v2}, {\mkr and after conducting the relevant estimates (see Appendix~\ref{Sec_Apriori_est_SPDE_key2}),} we arrive at, in large probability, and for $t$ in  $[0, T/\epsilon]$,
\be \label{Key_est2}
\|u_{\s} (t)\|_V^2 \le C \|u_{\s}(0)\|_V^2 + C  \epsilon^2  \int_0^t \hspace{-.05cm} e^{-  \delta(t-s) + 2 \sqrt{\epsilon}(W_t - W_s)} \d s, 
\ee
namely $\|u_{\s} (t)\|_V^2  \leq C \epsilon^2, $ provided that $\|u_{\s} (0,\omega)\|_V \sim \epsilon$ and $0\le \epsilon \le \epsilon^*$ for some $\epsilon^*>0$ sufficiently small.

Thus, if 
\be
|x(0)| \sim \sqrt{\epsilon},\; \; \|u_{\s}(0)\| \sim \epsilon, \; \text{ and } \|u_{\s} (0)\|_V \sim \epsilon.
\ee
we have that for any $T>0$ and any $\chi$ in $(0,1)$, 
\be \label{A_priori_SPDE_v3}
\mathbb{P}\Big(\textrm{(I), (II), and (III)} \Big) \geq 1 - \chi.
\ee

\subsection{Fluid problems subject to fluctuations: Low- and high-mode error estimates}
From these {\mkr probabilistic estimates we prove next our main error estimate theorem}. {\mkr This theorem gives the} error estimates {\mkr made on} the low- and high-mode dynamics, {\mkr when they are respectively approximated by the solution $X(t)$ to the reduced equation Eq.~\eqref{Eq_reduced}, and the parameterization 
$\Phi(X(t),t)$ given by  \eqref{Eq_Phi_Burgers}-\eqref{Eq_Phi_Burgers_comp}.} 
\bt\label{Thm_error_estimates}
 Assume that Conditions \ref{Cond_L1}, \ref{Cond_L2}, and \ref{Cond_A3} are satisfied. Denote $\epsilon = \beta_1(\lambda) >0$. Let  $u(t) = x(t) \boldsymbol{e}_1+u_\s(t)$ be a solution to the SPDE \eqref{Eq_Quadratic}.

 Assume furthermore that Condition~\ref{Cond_P} is satisfied, that $\sigma = \sqrt{\epsilon}$ 
and that
\bea\label{Eq_cond_init}
& X(0)=x(0)  \sim \sqrt{\epsilon},\\
 &  \|u_{\s}(0)\| \sim \epsilon, \; \|u_{\s} (0)\|_V \sim \epsilon, \\
 & \; \|u_{\s}(0) - \Phi(X,0)\| \sim \epsilon^{3/2}.
\eea
Consider the errors
\bea
&a(t) = x(t) - X(t), \; \textrm{(low-mode error)}\\
&b(t) = u_{\s}(t) - \Phi(X,t)  \; \textrm{(high-mode error)},
\eea
in which $X(t)$ solves the reduced equation Eq.~\eqref{Eq_reduced}, and the parameterization 
$\Phi$ is given by  \eqref{Eq_Phi_Burgers}-\eqref{Eq_Phi_Burgers_comp}.

Then, for any $T>0$ and any $\chi$ in $(0,1)$, there exist $\epsilon^{*} > 0$ and $C>0$ such that for any $\epsilon$ in $[0, \epsilon^{*}]$, the following error estimate holds:
\be \label{Err_est_goal}
\mathbb{P}\Big( \sup_{[0, T/\epsilon]} |a(t)| \le C \epsilon \; \text{and}\; \sup_{[0, T/\epsilon]} \|b(t)\| \le C \epsilon^{3/2} \Big) \ge 1 - \chi.
\ee
\et

\bp
First, let us write down the system of equations satisfied by the low-mode and high-mode errors, $a$ and $b$. By recalling that $(X,\Phi)$ and $(x,u_\s)$  satisfy \eqref{Eq_X_Phi} and \eqref{Eq_Quadratic_v2}, respectively,  we infer that $(a,b)$ satisfies  
\bea \label{Eq_error}
\hspace{-.25cm}\mathrm{d} a  &= \Big(\epsilon a  +  x\mathcal{F}_{\c}(u_\s) + \Pi_{1}B(u_\s, u_\s)  - X \mathcal{F}_{\c} (\Phi(X,t)) \Big) \mathrm{d} t + \sigma a \circ \d W_t, \\
\hspace{-.25cm} \mathrm{d} b &= \Big(L_\lambda^\s b  + x^2 B_{11}^{\s}  + x\mathcal{F}_\s(u_\s)+\Pi_\s B(u_\s, u_\s)- X^2 B_{11}^{\s} - 2 \Phi \mathcal{F}_{\c} (\Phi(X,t)) \Big) \mathrm{d} t + \sigma b \circ \d W_t,
\eea
with $\mathcal{F}_{\c}$ given by \eqref{Def_Fc}, $B_{11}^{\s}$ by \eqref{Def_B_11s}, and $\Pi_{1}B$ and $\mathcal{F}_\s$ defined in \eqref{Eq_Fs}.

We detail below the main elements to derive \eqref{Err_est_goal} from \eqref{Eq_error}. 
We start with the estimates of the nonlinear terms controlling the amplitudes of $a(t)$ and $b(t)$, in their respective dynamical equations. 
For the $a$-equation, this term is given  by $ \Pi_1 B(u_s, u_{\s})+ I(t)$ with 
\be
I(t)= x \mathcal{F}_{\c} (u_{\s}) - X \mathcal{F}_{\c} (\Phi(X,t)),
\ee
while for the $b$-equation, it is given by $J(t)+ (x^2 - X^2) B_{11}^{\s}$ with 
\be
J(t)= x \mathcal{F}_{\s} (u_\s) + \Pi_{\s} B(u_{\s}, u_{\s}) - 2 \Phi \mathcal{F}_{\c} (\Phi(X,t)).
\ee 

Observe that
\bea
I(t) & = a \mathcal{F}_{\c} (u_{\s}) + X (\mathcal{F}_{\c} (u_{\s}) - \mathcal{F}_{\c} (\Phi(X,t))) \\
& = a \mathcal{F}_{\c} (u_{\s}) +  X \langle B(\boldsymbol{e}_1, b) + B(b,\boldsymbol{e}_1), \boldsymbol{e}_1 \rangle,
\eea
and that 
\bea\label{I_interm}
 a I(t)    & \le  |\mathcal{F}_{\c} (u_{\s})|\, |a|^2 + C |a| \, |X| \|b\|_V \\
& \le  C \|u_{\s}\|_{V} |a|^2 + C \mathfrak{C} |a|\, |X| \|b\|_{\boldsymbol{\nu}},
\eea
due to the definition of $\mathcal{F}_{\c}$, and the norm equivalence between $\|\cdot\|_{\boldsymbol{\nu}}$ and $\|b\|_V$ (Condition~\ref{Cond_L1}). In the sequel, $C > 0$ will denote a generic constant that is allowed to change in the course of the estimates.

Now let us take $\Omega^*$ be the subset of $\Omega$ over which the a priori bounds  \eqref{A_priori_SPDE_v3} and \eqref{A_priori_goal} about $(x,u_\s)$ and $(X,\Phi)$, hold. 
Then over such a subset of events and for any $t$ in $[0,T/\epsilon]$, we deduce from \eqref{I_interm}, that 
\be\label{Est_fora1}
a(t) I(t)  \le C \epsilon  a(t)^2 + \frac{(1-\eta)}{2} \|b(t)\|_{\boldsymbol{\nu}}^2,
\ee
by using Young's inequality, \eqref{A_priori_SPDE_v3} to control $\|u_{\s}\|_{V}$, and \eqref{A_priori_goal} to control $X$. 
The factor $(1-\eta)/2$ in \eqref{Est_fora1} is chosen through application of the Young's inequality, in order to be appropriately ``absorbed'' by the linear damping effects on the high-modes, namely
\bea\label{stable_linear_part_est}
\langle L_\lambda^\s b, b \rangle & = \langle L_\lambda b, b \rangle \\
&  =  \langle (-\eta A + P_\lambda)  b, b \rangle  - (1-\eta) \langle A b, b\rangle \\
& \le - \delta \|b\|^2 - (1-\eta) \| b\|_{\boldsymbol{\nu}}^2,
\eea
where the latter inequality is a consequence of Conditions~\ref{Cond_L1} and \ref{Cond_L2}. 

As a result, we have over $\Omega^*$ and for any $t$ in $[0,T/\epsilon]$,
\be\label{Est_M1}
a I   + \langle L_\lambda^\s b, b \rangle \leq C \epsilon  a^2- \delta \|b\|^2 -\frac{(1-\eta)}{2} \|b\|_{\boldsymbol{\nu}}^2  \leq C \epsilon  a^2- \delta \|b\|^2. 
\ee

Similarly, we have over $\Omega^*$ and for any $t$ in $[0,T/\epsilon]$,
\be\label{Est_fora2}
a \Pi_1 B(u_s, u_{\s})  \le C \epsilon^2 |a| \le C \epsilon^3 + \epsilon a^2.
\ee

The estimates \eqref{Est_M1} and \eqref{Est_fora2} are used to control $a(t)$. To control $b(t)$ one needs to estimate $\langle  J(t), b \rangle$ and $\langle  (x^2 - X^2) B_{11}^{\s}, b \rangle$.

Here again by exploiting the a priori bounds \eqref{A_priori_goal}, \eqref{A_priori_SPDE_v3} and the Young's inequality, we infer that 
\be \label{nonlin_estimate1_eqn_b}
\langle  (x^2 - X^2) B_{11}^{\s}, b \rangle \le C \sqrt{\epsilon} |a| \|b\|  \le C \epsilon a^2 + \frac{\delta}{2} \|b\|^2, 
\ee
and
\be \label{nonlin_estimate2_eqn_b}
\langle  J(t), b \rangle \le  C \epsilon^{3/2} \|b\| \le C \epsilon^3  + \frac{\delta}{2} \|b\|^2,
\ee
still over over $\Omega^*$ and for any $t$ in $[0,T/\epsilon]$. 

From these estimates, and denoting by $M=a I+a \Pi_1B(u_s, u_{\s}) $, we observe, using \eqref{Est_M1}, that 
\be
M+ \langle  (x^2 - X^2) B_{11}^{\s}, b \rangle+  \langle  J, b \rangle  + \langle L_\lambda^\s b, b \rangle +  \epsilon a^2 \leq C_1 \epsilon^3  +C_2 \epsilon a^2, 
\ee 
where $C_1$ and $C_2$ are two positive constants. 

Following similar steps as in Appendix \ref{Sec_Apriori_est_SPDE_key1} (in particular after transformation as in \eqref{Eq_random_transform2}), we arrive then at
\bea
|a(t)|^2 + \|b(t)\|^2 & \le e^{C_2 \epsilon t + 2 \sqrt{\epsilon} W_t} (|a(0)|^2 + \|b(0)\|^2) + C_1 \epsilon^3 \int_0^t e^{C_2 \epsilon (t-s) + 2 \sqrt{\epsilon} (W_t - W_s)} \d s \\
& \le e^{C_2 T + 2 \gamma} (|a(0)|^2 + \|b(0)\|^2) + C_1 \epsilon^2,
\eea
where $\gamma$ is given by Lemma~\ref{Lemma_BM_bound}. 

 Note that $a(0) = 0$. So if we assume $\|b(0)\| \sim \epsilon$, then, 
\be \label{Error_est_v2}
|a(t)|^2 + \|b(t)\|^2 \le C \epsilon^2, \; \; \omega \in \Omega^*, \; t \in [0, T/\epsilon].
\ee
From this estimate, one can furthermore improve the estimate of  $\|b\|$.
To do so, note that by using the Young's inequality differently than for obtaining \eqref{nonlin_estimate1_eqn_b}, leads to 
\be \label{Est_nln_termb2}
\langle  (x^2 - X^2) B_{11}^{\s}, b \rangle \le C \sqrt{\epsilon} |a| \|b\|  \le C \epsilon a^2 + \frac{\delta}{4} \|b\|^2. 
\ee
Similarly we can arrange the constants through Young's inequality to get $\langle  J(t), b \rangle \le  C \epsilon^{3/2} \|b\| \le C \epsilon^3  + \delta\|b\|^2 /4.$

Now since $|a(t)|^2 \le C \epsilon^2$ due to \eqref{Error_est_v2}, we have
 \be \label{nonlin_estimate3_eqn_b}
\langle  (x^2 - X^2) B_{11}^{\s}, b \rangle \le C \epsilon^3 + \frac{\delta}{4} \|b\|^2. 
\ee
Then 
\be  \label{nonlin_estimate4}
\langle  (x^2 - X^2) B_{11}^{\s}, b \rangle+  \langle  J, b \rangle  + \langle L_\lambda^\s b, b \rangle \leq C \epsilon^3 -\frac{\delta}{2} \|b\|^2,
\ee 
which, by application of Gronwall's inequality (still following similar steps as in Appendix \ref{Sec_Apriori_est_SPDE_key1}), gives 
\be \label{Error_est_v3a}
\|b(t)\|^2 \le  e^{- \delta t + 2 \gamma} \|b(0)\|^2 + C \epsilon^3, \quad  \omega \in \Omega^*, \; t \in [0, T/\epsilon].
\ee
Thus, assuming $\|b(0)\| \sim \epsilon^{3/2}$, leads finally to 
\be \label{Error_est_v3b}
\|b(t)\|^2 \le C \epsilon^3, \quad  \omega \in \Omega^*,\; t \in [0, T/\epsilon].
\ee
The desired estimate \eqref{Err_est_goal} follows from \eqref{Error_est_v2} and \eqref{Error_est_v3b}. The proof is complete.

\ep

\needspace{2ex}
\br \label{err_esti_rmk}

\hspace*{.1em}  

\bi \setlength\itemsep{0.5em}

\item[i)] In the condition \eqref{Eq_cond_init}, if we drop the requirement $\|u_{\s}(0) - \Phi(X,0)\| \sim \epsilon^{3/2}$, the estimate \eqref{Error_est_v2} still holds while \eqref{Error_est_v3b} will be true after skipping a transient time of order $|\ln(\epsilon)|/\delta$, with $\delta$ given by \eqref{Eq_P_control_a}.

Indeed, since it is assumed that $X(0)  \sim \sqrt{\epsilon}$, then $\Phi(X,0) \sim \epsilon$. This together with the condition $\|u_{\s}(0)\| \sim \epsilon$ leads to $\|b(0)\| = \|u_{\s}(0) - \Phi(X,0)\| \sim \epsilon$. As a result, \eqref{Error_est_v2} still holds. 

Note also that since $\|b(0)\| \sim \epsilon$, then $e^{- \delta t + 2 \gamma} \|b(0)\|^2 \sim \epsilon^3$ for all $t \ge |\ln(\epsilon)|/\delta$. This together with \eqref{Error_est_v3a} leads to 
\be \label{Error_est_v3c}
\|b(t,\omega)\|^2 \le C \epsilon^3, \quad \forall \omega \in \Omega^*, t \in [|\ln(\epsilon)|/\delta, T/\epsilon].
\ee

\item[ii)] Note also that if we drop the requirement $\|u_{\s}(0) - \Phi(X,0)\| \sim \epsilon^{3/2}$ in the condition \eqref{Eq_cond_init}, then as a direct reinterpretation of \eqref{Error_est_v2}, the following relaxed version of the error estimate \eqref{Err_est_goal} holds without skipping any transient dynamics: 
\be \label{Err_est_relaxed}
\mathbb{P}\Big( \sup_{[0, T/\epsilon]} |a(t)| \le C \epsilon \; \text{and}\; \sup_{[0, T/\epsilon]} \|b(t)\| \le C \epsilon \Big) \ge 1 - \chi.
\ee 

\item[iii)] A close inspection of the proofs about the  error estimates shows that the requirement $\sigma = \sqrt{\epsilon}$ could be relaxed to the cases $ 0<\sigma \leq  \sqrt{\epsilon}$.

\item[iv)] Note that $\delta < |\beta_2(\lambda)|$, as can be seen by replacing $b$ in \eqref{stable_linear_part_est} with $\bm{e}_2$ and using $\langle L_\lambda \bm{e}_2,  \bm{e}_2 \rangle = \beta_2(\lambda) \|\bm{e}_2\|^2$. We can actually reduce the transient time $|\ln(\epsilon)|/\delta$ in \eqref{Error_est_v3c} to $|\ln(\epsilon)| / |\beta_2(\lambda)|$. This is because the estimate \eqref{Error_est_v3a} still holds by replacing the exponent $-\delta t$ therein by $-|\beta_2(\lambda)|t$. Indeed, we just need to use $|\beta_2(\lambda)|$ in place of $\delta$ when applying the Young's inequality that leads to the two estimates \eqref{Est_nln_termb2} and \eqref{nonlin_estimate3_eqn_b}. These together with $\langle L_\lambda^\s b, b \rangle \le \beta_2(\lambda) \|b\|^2$ lead to \eqref{nonlin_estimate4} with $\delta$ therein replaced by $|\beta_2(\lambda)|$. 
\ei

\er

\section{Applications to Rayleigh-B\'enard convection}\label{RBC_sec}

\subsection{The stochastic Rayleigh-B\'enard model and its mathematical formulation} We consider the following non-dimensionalized Boussinesq equations driven by a linear multiplicative noise:
\begin{equation} \label{Eq_RBC}
\begin{aligned}
&\mathrm{d} \vecu=\big( \Delta \vecu- \nabla p+ \sqrt{\Ra} \theta \vec{\bm{k}} - (\vecu\cdot\nabla )\vecu \big) \mathrm{d}t + \sigma \vecu \circ \mathrm{d} W_t, \\
&\mathrm{d} \theta = \big( \Delta \theta + \sqrt{\Ra} w -(\vecu\cdot\nabla )\theta \big) \mathrm{d}t  + \sigma  \theta \circ \mathrm{d} W_t,\\
&\text{div\,}\vecu=0.
\end{aligned}
\end{equation}
Here, the unknown functions are the velocity field $\vecu=(u, w)$, the pressure $p$, and the temperature fluctuation $\theta$. $\vec{\bm{k}} = (0,1)$ is the vertical unit vector, $\sigma>0$ is the noise amplitude. The nondimensional parameter $\Ra$ is the Rayleigh number, which serves as the control parameter for the bifurcation. We consider the case that the fluid is confined in a $2D$ nondimensional rectangular domain $\D = (0, L) \times (0, 1)$, where $L > 0$ represents the aspect ratio between the width and the height of the domain. To simplify the presentation, the above equations \eqref{Eq_RBC} are presented for the case that the Prandtl number, $\Pr$, is taken to be one. In the following, we use $(x,z)$ to denote the coordinates, $x$ for the horizontal direction and $z$ for the vertical direction.

Various physically sound boundary conditions can be handled in the theoretic setting; see \cite[Sec.~4.1.3]{MW14}. To fix ideas, we consider free-slip boundary condition for $\vecu$, and Dirichlet boundary condition for $T$  on the top and the bottom boundaries and Neumann boundary condition on the lateral boundaries:
\begin{equation} \label{RBC_bdry condition}
\begin{aligned} 
&u=0,\ \frac{\partial w}{\partial x}=\frac{\partial \theta}{\partial
x}=0  &&  \text{at}\ x=0,L,\\
&w=\theta=0,\ \frac{\partial u}{\partial z}=0  &&  \text{at}\ z=0,1.
\end{aligned}
\end{equation}

For the problem \eqref{Eq_RBC}-\eqref{RBC_bdry condition}, we set the spaces
\bea \label{RBC_spaces} 
& H = \{(\vecu,\theta) \in (L^2(\D ))^3 \  | \ \text{div\,} \vecu=0,\vecu\cdot \bm{n}|_{\partial\D }=0\},\\
& V = \big \{(\vecu,\theta) \in (H^1(\D ))^3 \cap H \ \big |\   \theta=0 \text{ at } z=0, 1\big \},\\
& H_1 = \big \{(\vecu,\theta)\in (H^2(\D ))^3 \ \big| \ \text{div\,} \vecu=0 \text{ and  \eqref{RBC_bdry condition}  holds} \big\},
\eea
where $\bm{n}$ is the unit outward normal vector to $\partial \D$.

Let $L_{\Ra}: H_1 \rightarrow H$ be defined for any $\bm{\psi} =(\vecu,\theta)$ in $ H_1$  by
\bea \label{RBC_lin_operators}
& L_{\Ra}= -A+P_{\Ra}, \text{ with }\\
& A\bm{\psi} = (- \mathbb{L} \Delta \vecu, - \Delta \theta), \\
& P_\Ra \bm{\psi} = ( \sqrt{\Ra} \, \mathbb{L}  (\theta \vec{\bm{k}}), \sqrt{\Ra} \,w),
\eea
where $\mathbb{L}$ denotes the Leray projection defined on $ (L^2(\D ))^2$, which projects each element into the divergence-free subspace of $ (L^2(\D ))^2$.  We define also $B: V \times V \rightarrow H$ as
\be \label{RBC_nonlin}
B(\bm{\psi}_1,\bm{\psi}_2)= (- \mathbb{L} ((\vecu_1\cdot\nabla )\vecu_2),-(\vecu_1 \cdot\nabla )\theta_2),
\ee
for all $\bm{\psi}_1 =(\vecu_1,\theta_1)$ and $\bm{\psi}_2 =(\vecu_2,\theta_2)$ in $V$.

We denote by $\|\cdot\|$ the norm on $H$ defined for any $\bm{\psi}=(u,w,\theta)$ by 
\bes
 \| \boldsymbol{\psi} \|^2= |u|_{L^2(\mathcal{D})}^2 +| w|_{L^2(\mathcal{D})}^2+| \theta|_{L^2(\mathcal{D})}^2,
\ees
and the norm on $V$, the norm $\|\cdot\|_V$ defined by
\bes
 \| \boldsymbol{\psi} \|_V^2=\| \nabla  \boldsymbol{\psi} \| + \| \boldsymbol{\psi}\|.
\ees

Then the problem \eqref{Eq_RBC}-\eqref{RBC_bdry condition} can be written as
\begin{equation} \label{Eq_RBC_v2}
\begin{aligned} &\mathrm{d} \bm{\psi} =\big( L_{\Ra}\bm{\psi} +B(\bm{\psi}, \bm{\psi}) \big)\mathrm{d}t  + \diffusion \bm{\psi} \circ \mathrm{d} W_t.
\end{aligned}
\end{equation}

In the following, we assume for the non-dimensionalized spatial domain $\mathcal{D}= [0, L] \times [0,1]$ that $L$ is chosen such that the deterministic analogue of \eqref{Eq_RBC} admits a pitchfork bifurcation as $\Ra$ crosses its critical value $\Ra_c$ from below. As pointed out in Sec.~\ref{Sec_verifying_A1-A3} when verifying Condition~\ref{Cond_A2}, there is exactly one critical wave vector for all values of $L$ except those from a subset of $\mathbb{R}^+$ with measure zero; cf.~\eqref{Exception_set} for the latter exceptional set. So the pitchfork bifurcation scenario is a generic situation for the deterministic case.

\subsection{Verification of Conditions~\ref{Cond_A1}-\ref{Cond_A3} and \ref{Cond_P}} \label{Sec_verifying_A1-A3}  

In this section, we provide details on verifying the Conditions~\ref{Cond_A1}-\ref{Cond_A3} and \ref{Cond_P}.

\medskip
{\noindent \bf Verification of Condition~\ref{Cond_A1}.} Since the linear operator $L_{\Ra}$ defined by \eqref{RBC_lin_operators} involves the Leray projection $\mathbb{L}$, it is not immediately clear that $L_{\Ra}$ is self-adjoint. However, this property will be verified once the eigenfunctions of $L_{\Ra}$ are computed. 

We consider thus the eigenvalue problem $L_{\Ra} \bm{\psi} = \beta \bm{\psi}$ in the space $\cH_1$. That is
\begin{equation} \label{RBC_eign}
\begin{aligned}
& \mathbb{L} \Delta \vecu +  \sqrt{\Ra}\, \mathbb{L} (\theta \vec{\bm{k}}) =\beta \vecu,\\
&\Delta \theta +  \sqrt{\Ra}\, w =\beta \theta,\\
&\text{div\,}\vecu=0.
\end{aligned}
\end{equation}
Recall that the eigenvalues and eigenfunctions for \eqref{RBC_eign} are given as follows (see e.g.~\cite[Section 4.1.6]{MW14}): 
\begin{itemize}
\item[] {\hspace{-2.5em}  {\bf Group one}}:
\begin{align} \label{RBC_beta_0k}
&\beta_{0k} =-k^2\pi^2, \quad \bm{e}_{0k}=\sqrt{2/L} (0,0,\sin ( k \pi z) ), \; \; k \in \mathbb{N}. 
\end{align}
\item[] {\hspace{-2.5em}  {\bf Group two}}: The eigenvalues are given, for all $j, k$  in $\mathbb{N}$, by 
\bea \label{RBC_beta_jk}
\beta_{jk}^{\pm}(\Ra) = - \gamma_{jk}^2 \pm \sqrt{\frac{\Ra \alpha_j^2}{\gamma_{jk}^2}},
\eea
with 
\be
\alpha_j = j\pi/L, \;  \gamma_{jk}=\sqrt{\alpha_j^2+k^2\pi^2};
\ee
and the corresponding eigenfunctions are given by 
\bea \label{eigen vector}
\bm{e}^{\pm}_{jk}=\Bigl( & a_{jk}^{\pm} \sin ( \alpha_j x) \cos ( k \pi z), b_{jk}^{\pm} \cos (\alpha_j  x) \sin (k \pi z),  c_{jk}^{\pm} \cos ( \alpha_j  x) \sin (k \pi z) \Bigr ),
\eea
where 
\be
a_{jk}^{\pm} =  - k \pi N^{\pm}_{jk}, \quad b_{jk}^{\pm} = \alpha_j N^{\pm}_{jk}, \quad c_{jk}^{\pm} = \frac{\sqrt{\Ra} \alpha_j}{\gamma_{jk}^2+\beta_{jk}^{\pm}(\Ra)} N^{\pm}_{jk}, 
\ee
with 
\be
N^{\pm}_{jk} = \left(\frac{4 (\gamma_{jk}^2 + \beta_{jk}^{\pm}(\Ra))^2}{L(\gamma_{jk}^2(\gamma_{jk}^2 + \beta_{jk}^{\pm}(\Ra))^2 + \Ra \alpha_j^2)} \right)^{1/2}. 
\ee
\end{itemize}

Note that in Group two above, the coefficient $\frac{1}{\gamma^2_{jk}+\beta^{\pm}_{jk}(\Ra)}$ in the third component of $\bm{e}^{\pm}_{jk}$ is well defined for $j,k$ in $\mathbb{N}$ and for all $\Ra > 0$ because from \eqref{RBC_beta_jk}, we know that $\beta_{jk}^{\pm}(\Ra) \neq -\gamma^2_{jk}$. The constant $N^{\pm}_{jk}$, included in each of the three components of $\bm{e}^{\pm}_{jk}$, is a normalization constant ensuring $\|\bm{e}^{\pm}_{jk}\| = 1$.

Now, for any eigenfunctions $\bm{\psi}_1$ and $\bm{\psi}_2$ listed above, it follows from a direct calculation that the following identity holds 
\be
\langle L_{\Ra} \bm{\psi}_1,\bm{\psi}_2 \rangle = \langle \bm{\psi}_1, L_{\Ra} \bm{\psi}_2 \rangle. 
\ee
Thus, $L_{\Ra}$ is self-adjoint.

\medskip
{\noindent \bf Verification of Condition~\ref{Cond_A2}.} Note that only eigenvalues from the second group, i.e.,  $\beta_{jk}^{\pm}(\Ra)$ given by \eqref{RBC_beta_jk} can change signs when $\Ra$ varies. To identify the critical value of $\Ra$, we equate $\beta_{jk}^{\pm}(\Ra)$ to zero to obtain 
\be
\Ra = \gamma_{jk}^6 / \alpha_j^2.
\ee 
The critical Rayleigh number $\Ra_c$ is thus given by:
\bea \label{Ra_c}
\Ra_c & = \min_{j,k\in \mathbb{N}} \frac{\gamma_{jk}^6}{\alpha_j^2} =  \min_{j,k\in \mathbb{N}} \frac{(k^2\pi^2 + j^2 \pi^2/L^2)^3}{j^2\pi^2/L^2} \\
& =  \min_{j\in \mathbb{N}} \frac{\pi^4 (1 + j^2/L^2)^3}{j^2/L^2}.
\eea
One can readily check that the minimum is achieved at either $(j,k) = (\lfloor L/\sqrt{2} \rfloor,1)$ or $(\lceil L/\sqrt{2} \rceil,1)$, or both, depending on the value of $L$, where $\lfloor L/\sqrt{2} \rfloor$ denotes the largest integer below $L/\sqrt{2}$, and  $\lceil L/\sqrt{2}\, \rceil$ denotes the smallest integer above $L/\sqrt{2}$. The cases with the minimum achieved at both of the two indices $(j_{1},1) = (\lfloor L/\sqrt{2} \rfloor,1)$ and $(j_{2},1) = (\lceil L/\sqrt{2} \, \rceil,1)$ occur when $\gamma_{j_{1}, 1}^6/\alpha_{j_{1}}^2 = \gamma_{j_{2}, 1}^6 / \alpha_{j_{2}}^2$. That is when $L$ is chosen from the following subset of $\mathbb{R}^+$: 
\be \label{Exception_set}
\mathcal{S} = \left\{ L > 0 \; \Big | \;  \frac{(1 + (j_1)^2/L^2)^3}{(j_1)^2} = \frac{(1 + (j_2)^2/L^2)^3}{(j_2)^2} \text{ with } j_{1} = \lfloor L/\sqrt{2} \rfloor,  j_2 = \lceil L/\sqrt{2} \, \rceil, j_1 \neq j_2 \right\}.
\ee
Note that $\mathcal{S}$ has measure zero. In the following, we focus on the case of pitchfork bifurcation, and assume thus that $L$ takes values outside of $\mathcal{S}$. We denote the unique index that achieves the minimum value $\Ra_c$ by $(j_c,1)$.

Thus, Condition~\ref{Cond_A2} is verified with $\Ra$ playing the role of $\lambda$ and $\beta_1(\lambda)$ in \eqref{PES} taken to be $\beta_{j_c1}^{+}(\Ra)$ here. Although the eigenelements here are labeled using a double index associated with the corresponding wave vectors, they can apparently also be labeled using a single index to fit into the setting of \eqref{PES}.

\medskip
{\noindent \bf Verification of Condition~\ref{Cond_A3}.} Since the spatial domain is taken to be a (2D) rectangle, then $B$ defined by \eqref{RBC_nonlin} is a continuous bilinear map from $V \times V$ to $\cH$ with $V$ being the the subspace of $(H^1(\mathcal{D}))^3$ given in \eqref{RBC_spaces}.

Note also that $B$ satisfies that 
\be \label{Eq_energy_conservation_RBC}
\langle B(\bm{\psi}_1,\bm{\psi}_2), \bm{\psi}_3 \rangle =  - \langle B(\bm{\psi}_1,\bm{\psi}_3), \bm{\psi}_2 \rangle,
\ee
for all $\bm{\psi}_1$, $\bm{\psi}_2$, and $\bm{\psi}_3$ in $V$. The condition \eqref{Eq_energy_conservation} is a direct consequence of \eqref{Eq_energy_conservation_RBC}. 

It follows also from a direct calculation that 
\bea \label{Eq_interaction_table_RBC}
& \hspace{-0.7em} \langle B(\bm{e}^{+}_{j_c1},\bm{e}^{+}_{j_c1}), \bm{e}^{\pm}_{jk})  \rangle \!=\! 0,  \; \forall j, k \in \mathbb{N}, \\
&\hspace{-0.7em}  \langle B(\bm{e}^{+}_{j_c1},\bm{e}^{+}_{j_c1}), \bm{e}_{0k})  \rangle \!=\! \begin{cases}
{\displaystyle - \frac{\pi \sqrt{2 \Ra} (\gamma_{j_c1}^2 + \beta_{j_c1}^{+}(\Ra)) \alpha_{j_c}^2}{\sqrt{L} (\gamma^2_{j_c1}(\gamma_{j_c1}^2 + \beta_{j_c1}^{+}(\Ra))^2+\Ra \alpha_{j_c}^2)}} &  \text{ if } k = 2, \\
0 & \hspace{-0.2em} \text{ otherwise.} 
\end{cases}
\eea
The condition \eqref{cond_interaction_table} is thus verified with $N$ therein equals 2 here, after re-arranging the double index into a single one. 

\medskip
{\noindent \bf Verification of Condition~\ref{Cond_P}.} Recall that $N$ in \eqref{Cond_supercritical} is $2$ here (see again \eqref{Eq_interaction_table_RBC}). From \eqref{Eq_interaction_table_RBC}, we also know 
\be \label{Eq_B112}
B_{11}^2 =  \langle B(\bm{e}^{+}_{j_c1},\bm{e}^{+}_{j_c1}), \bm{e}_{02})  \rangle = - \frac{\pi \sqrt{2 \Ra} (\gamma_{j_c1}^2 + \beta_{j_c1}^{+}(\Ra)) \alpha_{j_c}^2}{\sqrt{L} (\gamma^2_{j_c1}(\gamma_{j_c1}^2 + \beta_{j_c1}^{+}(\Ra))^2+\Ra \alpha_{j_c}^2)} < 0.
\ee
Note also that
\bea \label{Eq_B121_coef}
\big( B_{12}^1 + B_{21}^1\big) & = \langle B(\bm{e}^{+}_{j_c1},\bm{e}_{02}), \bm{e}^{+}_{j_c1} \rangle + \langle B(\bm{e}_{02}, \bm{e}^{+}_{j_c1}), \bm{e}^{+}_{j_c1} \rangle \\
& = \frac{\pi \sqrt{2 \Ra} (\gamma_{j_c1}^2 + \beta_{j_c1}^{+}(\Ra)) \alpha_{j_c}^2}{\sqrt{L} (\gamma^2_{j_c1}(\gamma_{j_c1}^2 + \beta_{j_c1}^{+}(\Ra))^2+\Ra \alpha_{j_c}^2)}. 
\eea
Since $M_{11}^{2,\Ra}$ is always positive, we have then
\be
B_{11}^2 M_{11}^{2,\Ra}  \big( B_{12}^1 + B_{21}^1\big) < 0.
\ee
Condition~\ref{Cond_P} is thus verified. 

\subsection{Verification of Conditions~\ref{Cond_L1} and \ref{Cond_L2}} \label{Sec_verifying_L1-L2}  

Since it has already been checked above that the linear operator $L_{\Ra}$ defined by \eqref{RBC_lin_operators} satisfies  Conditions~\ref{Cond_A1} and \ref{Cond_A2}, to verify Condition~\ref{Cond_L1}, it remains to show that $\sqrt{\langle A \boldsymbol{\psi}, \boldsymbol{\psi} \rangle}$ induces a norm on $V$ equivalent to $\| \boldsymbol{\psi} \|_{V}$. For this purpose, we note that for $A$ defined in \eqref{RBC_lin_operators}, one has, due to the boundary conditions, that for any $\boldsymbol{\psi}=(u,w,\theta)$ in $D(A) = H_1$:
\bes
\langle A \boldsymbol{\psi}, \boldsymbol{\psi} \rangle = |\nabla u|_{L^2(\mathcal{D})}^2 +| \nabla w|_{L^2(\mathcal{D})}^2+| \nabla \theta|_{L^2(\mathcal{D})}^2  = \| \nabla \boldsymbol{\psi} \|^2,
\ees
Note that
\bes
\| \nabla  \boldsymbol{\psi}\| \leq  \|\boldsymbol{\psi}\|_V \leq (1+C) \|  \boldsymbol{\psi}\|,
\ees
where $C>0$ is a generic constant in the Poincar\'e inequality. 
Thus, $\| \nabla \boldsymbol{\psi} \|$ defines a norm equivalent to  $\|\boldsymbol{\psi}\|_V$.

Condition~\ref{Cond_L2} is a consequence of the following lemma. 

\bl \label{Lem_RBC_cond_L}
There exist an $\eta$ in $(0,1)$, $\delta > 0$ and $\Ra^* > \Ra_c$, for which the condition \eqref{Eq_P_control} holds for the linear operator $L_{\Ra}$ defined in \eqref{RBC_lin_operators} for all $\Ra$ in $[\Ra_c, \Ra^*]$.
\el

See Appendix~\ref{Sec_proof_Lemma_RBC} for a proof.  

\subsection{Stochastic pitchfork bifurcation in large probability}\label{Sec_stoch_pitch}

Due to \eqref{Eq_interaction_table_RBC}, the parameterization $\Phi$ defined in \eqref{Eq_Phi_Burgers} reads here as
\be \label{Eq_Phi_RBC}
\Phi(X,t,\omega) = \Phi_{2}(X,t, \omega) \boldsymbol{e}_{02}, \quad \forall X \in \mathbb{R},
\ee
with 
\be \label{Eq_Phi_RBC_comp}
 \Phi_{2}(X, t, \omega) = B_{11}^{2} M_{11}^{2,\Ra} (t,\omega) X^2,
\ee
and $B_{11}^{2}$ given by \eqref{Eq_B112}.

The abstract reduced equation \eqref{Eq_reduced} for the stochastic RBC problem takes then the following explicit form 
\be\label{Eq_RBC_reduced}
\d X= (\beta_{j_c1}^{+}(\Ra) X-\alpha(\Ra) M_{11}^{2,\Ra} (t,\omega) X^3 )\d t +\sigma X \circ \d W_t,
\ee
with 
\be
\alpha(\Ra) = - B_{11}^2 \big( B_{12}^1 + B_{21}^1\big)=\frac{2\pi^2 \Ra (\gamma_{j_c1}^2 + \beta_{j_c1}^{+}(\Ra))^2 \alpha_{j_c}^4}{L (\gamma^2_{j_c1}(\gamma_{j_c1}^2 + \beta_{j_c1}^{+}(\Ra))^2+\Ra \alpha_{j_c}^2)^2}.
\ee
We call \eqref{Eq_RBC_reduced}, the non-Markovian RBC reduced equation.

\begin{SCfigure} 
\includegraphics[width=0.6\linewidth,height=.42\textwidth]{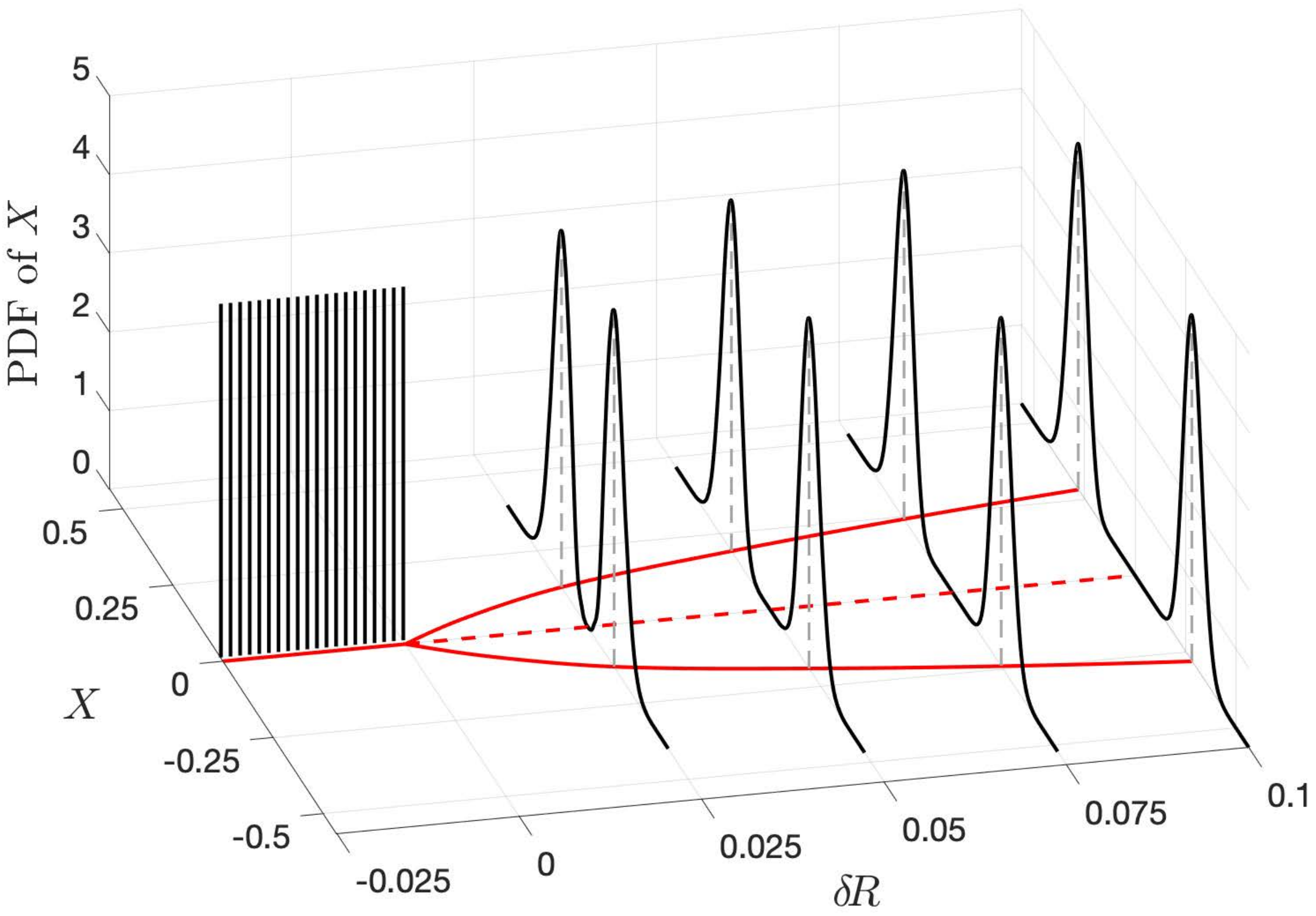}
\caption{\footnotesize {\bf The stochastic supercritical pitchfork bifurcation for the RBC reduced equation \eqref{Eq_RBC_reduced}.}
When $\deltaR = \Ra - \Ra_c \le 0$, the probability density function (PDF) is a Dirac delta function since $X=0$ is globally stable. When $\deltaR > 0$, two locally stable random steady states $\pm X_{\Ra}$ emerge and $X=0$ loses its stability. The most probable values of $\pm X_{\Ra}$ as $\deltaR$ varies form the solid red curves in the $(\deltaR,X)$-plane, and the unstable trivial steady state for $\deltaR > 0$ is marked by the dashed red line segment. The most probable states scale like $\sqrt{\deltaR}$, which is consistent with the theoretic results presented in Lemma~\ref{Lem_SS_estimation}, by noting that $\epsilon = \beta_{j_c1}^{+}(\Ra)$ in Lemma~\ref{Lem_SS_estimation} scales linearly with $\deltaR$; cf.~\eqref{Eq_eps_Ra_relation} below. For the parameters, we have set $L=3$, leading to $\Ra_c \approx 660.52$ and the critical wave number $j_c = 2$; and the noise amplitude is set to $\sigma = 0.01$.}\label{Fig_pitchfork_pdf}
\end{SCfigure}

Note that the RBC reduced equation \eqref{Eq_RBC_reduced} fits into the non-Markovian normal form of a supercritical pitchfork bifurcation as given by Appendix~\ref{Sec_1DSDE_pitchfork}. In particular, the trivial steady state $X=0$ is globally stable when the Rayleigh number $\Ra$ is below the critical value $\Ra_c$ given by \eqref{Ra_c}. It becomes unstable when $\Ra > \Ra_c$, and two locally stable random steady states emerge.
These bifurcated random equilibria are given by $\pm X_{\Ra}$, where 
\be \label{Eq_stationary_RBC_reduced}
X_{\Ra}(t,\omega) = \frac{1}{\sqrt{2 \alpha(\Ra) \int_{-\infty}^t M_{11}^{2,\Ra} (s, \omega) \exp(f_{\Ra}(t,s,\omega)) \d s}}, \quad \Ra \ge \Ra_c,
\ee
with  
\be  
f_{\Ra}(t,s,\omega) = -2 \beta_{j_c1}^{+}(\Ra) (t-s) - 2 \sigma (W_t(\omega) - W_s(\omega)).
\ee
This stochastic pitchfork bifurcation is shown in Fig.~\ref{Fig_pitchfork_pdf} at the level of the probability density function (PDF) of the (random) steady states for a particular choice of the domain size $L = 3$, with $\sigma = 0.01$. 

In the following, we provide an interpretation of the general error estimate results given by Theorem~\ref{Thm_error_estimates} for the stochastic RBC problem considered here. By doing so, we establish a rigorous link between the pitchfork bifurcation for the non-Markovian RBC reduced equation \eqref{Eq_RBC_reduced} and the dynamics of the full stochastic RBC problem \eqref{Eq_RBC_v2} when the noise amplitude parameter $\sigma$ scales like $\sqrt{\deltaR}$ with $\deltaR = \Ra - \Ra_c$. For this purpose, we introduce
\begin{subequations} \label{Def_SS_RBC}
\begin{align}
& A_{\Ra}^+(t)=X_{\Ra}(t)\boldsymbol{e}_{j_c1}^+ + \Phi(X_{\Ra}(t),t), \\
& A_{\Ra}^-(t)=-X_{\Ra}(t)\boldsymbol{e}^+_{j_c1} + \Phi(-X_{\Ra}(t),t),
\end{align}
\end{subequations} 
where $\boldsymbol{e}^{+}_{j_c1}$ is the eigenmode that becomes linearly unstable,
\bes 
\boldsymbol{e}^{+}_{j_c1}=\Bigl( - \pi \sin ( \alpha_{j_c} x) \cos (  \pi z), \alpha_{j_c} \cos (\alpha_{j_c}  x) \sin (\pi z),  \frac{\sqrt{\Ra} \alpha_{j_c}}{\gamma_{j_c 1}^2+\beta_{j_c1}^{+}(\Ra)} \cos ( \alpha_{j_c}  x) \sin ( \pi z) \Bigr ),
\ees
and the manifold function $\Phi$ is given by \eqref{Eq_Phi_RBC}. Then, the following theorem holds.

\bt \label{RBC_thm}
Consider the stochastic RBC problem \eqref{Eq_RBC_v2}. Assume that the domain aspect ratio parameter $L$ is chosen outside of the measure zero set $\mathcal{S}$ given by \eqref{Exception_set}, so that there is exactly one mode that becomes unstable as $\Ra$ crosses the first critical value $\Ra_c$ defined by \eqref{Ra_c}. Denote $\deltaR = \Ra - \Ra_c$. Assume that $\sigma = (\beta_{j_c1}^{+}(\Ra))^{1/2}$ and that
\be \label{RBC_IC_assumption}
X(0)= \langle \boldsymbol{\psi}(0), \boldsymbol{e}^+_{j_c1} \rangle \sim \sqrt{\deltaR}, \; \|\boldsymbol{\psi}_{\s}(0)\| \sim \deltaR,  \;\textrm{and} \;   \|\boldsymbol{\psi}_{\s}(0)\|_V \sim \deltaR.
\ee  
Then, if $X(0)>0$, the solution  $\boldsymbol{\psi} (t)$ of \eqref{Eq_RBC_v2} emanating from $\boldsymbol{\psi}(0)$
satisfies that for any $\chi$ in $(0,1)$ and $T > 0$, there exist $\Ra^*>\Ra_c$ and $C > 0$ such that the following estimate holds for all $\Ra$ in $[\Ra_c, \Ra^*]$:
\be \label{RBC_error_est}
\mathbb{P}\Big(\sup_{[0, \eta T/\deltaR]} \|\boldsymbol{\psi} (t) - A_{\Ra}^+(t) \|  \leq C \deltaR \Big) \ge 1 - \chi,
\ee
where $A_{\Ra}^+$ is defined in \eqref{Def_SS_RBC} and $\eta > 0$ is a constant depending only on the aspect ratio $L$. The same estimate holds with  $A_{\Ra}^-(t)$ replacing $A_{\Ra}^+(t)$, if $X(0)<0$. 
\et

\bp

The desired result is a direct consequence of Theorem~\ref{Thm_error_estimates} and Remark~\ref{err_esti_rmk}-(ii). Indeed, the Conditions \ref{Cond_L1}, \ref{Cond_L2}, \ref{Cond_A3}, and \ref{Cond_P} required in Theorem~\ref{Thm_error_estimates} are verified for the RBC problem considered here in Sections~\ref{Sec_verifying_A1-A3} and \ref{Sec_verifying_L1-L2}. To derive \eqref{RBC_error_est}, we are only left to show that 
the assumption \eqref{RBC_IC_assumption} on the initial condition is equivalent to the analogues given in \eqref{Eq_cond_init} stated in terms of the parameter $\epsilon = \beta_{j_c1}^{+}(\Ra)$ instead.\footnote{Note that the condition $\|u_{\s}(0) - \Phi(X,0)\| \sim \epsilon^{3/2}$ required in \eqref{Eq_cond_init} is not needed here as pointed out in Remark~\ref{err_esti_rmk}-(ii).} Namely, we just need to show that $\deltaR$ is proportional to $\epsilon$. The parameter $\eta$ in the estimate \eqref{RBC_error_est} is related to the associated proportionality constant when converting from $T/\epsilon$. We are thus only left with verifying the following scaling relation:
\be \label{Eq_eps_Ra_relation}
\beta_{j_c1}^{+}(\Ra) = \frac{(j_c\pi)^2/L^2}{ (\pi^2 + (j_c\pi)^2/L^2)^2}(\Ra - \Ra_c) + O(|\Ra - \Ra_c|^2),
\ee 
which can be derived through direct algebraic operations. Indeed, by introducing $\delta_{j_c} = \alpha_{j_c}^2 / \gamma_{j_c1}^2$, we have for $\Ra$ close to $\Ra_c$ that 
\bea
\sqrt{\delta_{j_c} \Ra} & = \sqrt{\delta_{j_c} \Ra_c + \delta_{j_c}(\Ra - \Ra_c)} = \sqrt{\delta_{j_c} \Ra_c} \sqrt{1 + (\Ra - \Ra_c)/\Ra_c} \\
& = \sqrt{\delta_{j_c} \Ra_c} \left(1 + \frac{(\Ra - \Ra_c)}{2\Ra_c} + O(|\Ra - \Ra_c|^2) \right).
\eea
With $\delta_{j_c}$ defined above, we can rewrite $\beta_{j_c1}^{+}(\Ra)$ given by \eqref{RBC_beta_jk} as $ \beta_{j_c1}^{+}(\Ra) = - \gamma_{j_c1}^2 + \sqrt{\delta_{j_c} \Ra}$, leading in turn to
\be
\beta_{j_c1}^{+}(\Ra) = - \gamma_{j_c1}^2 + \sqrt{\delta_{j_c} \Ra_c} +  \frac{\sqrt{\delta_{j_c}}(\Ra - \Ra_c)}{2 \sqrt{\Ra_c}} + O(|\Ra - \Ra_c|^2).
\ee
Since  $- \gamma_{j_c1}^2 + \sqrt{\delta_{j_c} \Ra_c} = \beta_{j_c1}^{+}(\Ra_c) = 0$, we get 
\be \label{Eq_eps_Ra_relation1}
\beta_{j_c1}^{+}(\Ra) = \frac{\sqrt{\delta_{j_c}}(\Ra - \Ra_c)}{2\sqrt{\Ra_c}} + O(|\Ra - \Ra_c|^2).
\ee
Recall also from \eqref{Ra_c} that $\Ra_c = \gamma_{j_c 1}^6 / \alpha_{j_c}^2$. Using this identity in \eqref{Eq_eps_Ra_relation1}, we get 
$\beta_{j_c1}^{+}(\Ra) = \frac{\alpha_{j_c}^2}{2\gamma_{j_c 1}^4}(\Ra - \Ra_c) + O(|\Ra - \Ra_c|^2)$, and \eqref{Eq_eps_Ra_relation} follows.
\ep

The above theorem shows that the SPDE solution $\boldsymbol{\psi}(t)$ of \eqref{Eq_RBC_v2} is, in large probability, within an $C \deltaR$-cones centered around either $A_{\Ra}^{+}$ or $A_{\Ra}^{-}$ defined in \eqref{Def_SS_RBC} depending on the sign of the projected initial data $X(0) = \langle \boldsymbol{\psi}(0), \boldsymbol{e}^+_{j_c1} \rangle$, for $t$ in $[0, \eta T/\deltaR]$. Recall also from Lemma~\ref{Lem_SS_estimation} (and the scaling relation \eqref{Eq_eps_Ra_relation}) that the magnitude of the bifurcated random steady states $\pm X_{\Ra}$ for the reduced equation \eqref{Eq_RBC_reduced} (and hence the corresponding lifted states $A_{\Ra}^{\pm}$ in $\cH$) are of order $\sqrt{\deltaR}$ with large probability for all $t > 0$. 

\begin{SCfigure} 
\includegraphics[width=0.62\linewidth, height=.5\textwidth]{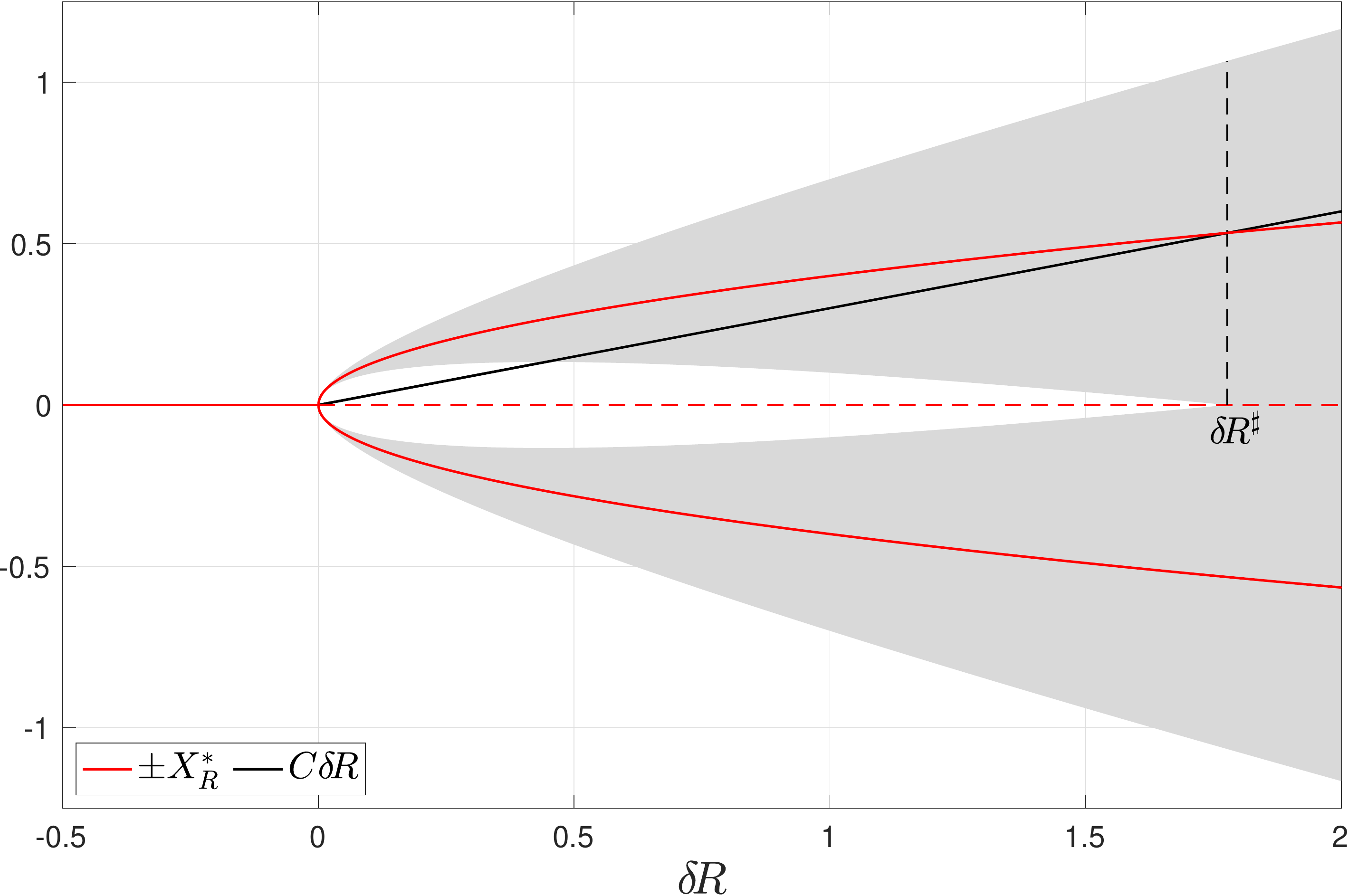}
\caption{\footnotesize  {\bf A schematic to illustrate the stochastic pitchfork transition scenario for the stochastic RBC problem \eqref{Eq_RBC_v2}.}
The red curves indicate the most probable values, denoted by $\pm X^*_{\Ra}$, of the random steady states $\pm X_{\Ra}$ for the reduced equation \eqref{Eq_RBC_reduced}. Recall that $X_{\Ra}$ scales like $\sqrt{\deltaR}$ with large probability; see again Lemma~\ref{Lem_SS_estimation} and \eqref{Eq_eps_Ra_relation}. The gray zones indicate the $C \deltaR$-cones centered at $\pm X^*_{\Ra}$ with the radius $C \deltaR$ being the upper bound appearing in the error estimate \eqref{RBC_error_est}. The vertical dashed line marks out the threshold $\deltaR^\sharp$ at which $X^*_{\Ra}$ intersects with the line segment $y = C \deltaR$. For $\deltaR$ in $(0, \deltaR^\sharp)$, the two $C \deltaR$-cones are disjoint from each other.
Theorem \ref{RBC_thm} ensures that for $\deltaR \le \min \{\deltaR^\sharp, \deltaR^*\}$, the (projected) SPDE solution from an initial condition satisfying \eqref{RBC_IC_assumption} is confined in one of the two cones with large probability and over a large time interval of order $1/\deltaR$. Here, $\deltaR^* = \Ra^* - \Ra_c$ with $\Ra^*$ being as specified in Theorem~\ref{RBC_thm}.}\label{Fig_Bif_diag}
\end{SCfigure}

As a result, there exists a $\deltaR^\sharp > 0$ such that for any $\deltaR$ in $(0, \deltaR^\sharp)$, the two $C \deltaR$-cones, one centered around $A_{\Ra}^{+}$ and the other around $A_{\Ra}^{-}$, do not intersect with large probability; moreover, the SPDE solutions emanating from initial data that satisfy the condition \eqref{RBC_IC_assumption} stay within one of these two cones with large probability and over a long time interval of order $1/\deltaR$. See Fig.~\ref{Fig_Bif_diag} for a schematic that illustrates this dynamical property in the projected unstable subspace. This scenario is thus consistent with the idea of a stochastic pitchfork transition for the SPDE dynamics that takes place when $\Ra$ crosses the threshold $\Ra_c$ from below, albeit subject to large probability and finite time intervals requirements.  

\section{Concluding remarks} \label{Sec_conclusion}

Based on a dynamical reformulation  of the leading-order approximation formulas for stochastic invariant manifolds, 
we derived rigorous energy estimates between the solutions of the corresponding reduced system  \eqref{Eq_X_Phi} and those of the original SPDE \eqref{Eq_Quadratic}. Such estimates as summarised in Theorem \ref{Thm_error_estimates} hold with large probability and are obtained for the case of a stochastic pitchfork bifurcation scenario. The key assumption \ref{Cond_L2}  that enables these estimates exploits dissipation effects brought by the stable modes of the linear operator in order to suitably counterbalance the loss of regularity due to the nonlinear terms. 

Although seeking a $\mathbb{P}$-almost sure description of the bifurcation for a given SPDE is still out of reach for most problems arising from applications\footnote{Recall that for stochastic ODEs already the question of how to describe a stochastic bifurcation, is not completely settled; see e.g.~\cite[Chap.~9]{Arnold98} and \cite{CF98}.}, we showed in this article that the dynamical properties of the SPDE solutions are actually captured by our reduced systems with large probability and over long time intervals, as long as the noise's intensity and the eigenvalue's magnitude of the mildly unstable mode, scale accordingly. While this type of characterizations has been investigated in the literature before, such as via the amplitude equation approach \cite{Blomker07}, these approaches do not seem to be directly applicable to problems considered here which are with quadratic nonlinearities and subject to multiplicative noise. While the focus here is on the pitchfork bifurcation, our approach can be suitably adapted to handle stochastic disturbances of fluids experiencing more general bifurcations such as $S^1$-attractor bifurcations \cite{MW05}, Hopf, and double-Hopf bifurcations. 

We have also shown that the underlying analytic formulas of our reduced systems near the onset of instability for SPDE models involve, unlike for PDEs,  non-Markovian coefficients which depend explicitly on the noise path (and its ``past"), the model's coefficients, the noise intensity, as well as the model's eigenelements. In that respect, the memory functions captured by these coefficients are not set arbitrarily, but derived in a consistent manner, so that in particular the characteristics of the memory function is self-consistently determined by the intensity of the random force from the model's equation.   

As aready pointed out in \cite{CLW15_vol2}, higher-order approximations of the genuine stochastic center manifold,  $h^\lambda_\omega(X,t)$, involve  also 
memory coefficients albeit of more complicated functional dependence on the noise path. 
The underlying high-order non-Markovian parameterizations have also analytic expressions that can be derived by solving the appropriate backward-forward (BF) systems  \cite[Chaps.~4 and 7.3]{CLW15_vol2}. The resulting random coefficients solve then new auxiliary SDEs that can also be used for deriving the error estimates as done for the leading-order approximation. 
 Such high-order non-Markovian parameterizations should lead to improved error estimates, as numerically shown in  \cite[Chap.~7]{CLW15_vol2} in a stochastic Burgers example; see also \cite{CL15}.

Finally and more generally, we mention that the usage of such backward-forward systems allows also for moving away from criticality, by adopting and adapting the variational approach and theory of optimal parameterizing manifolds (OPMs) of \cite{CLM16_Lorenz9D,CLM19_closure,Chekroun2021c} to the stochastic context. The new parameterizations obtained this way can be viewed as optimized homotopic deformations of those valid near criticality recalled in Theorem \ref{Thm_CM_approx} above. In the OPM approach, the optimization step is performed using data of the fully resolved problem (in a parsimonious way) by minimizing a least-square metric measuring typically the parameterization defect,  
and where the parameters to optimize are reduced to the backward integration times $\tau$ of the BF systems involved; see  \cite[Sec.~4.3]{CLM19_closure}. In the stochastic context, it means that the $M_n$-coefficients in \eqref{Phi_n_general}  are replaced away from criticality by integrals of the form $\int_{-\tau}^0 e^{g_n(\lambda) s + \sigma W_s(\omega)}\, \mathrm{d}s$, that are optimized (in $\tau$) per mode $\boldsymbol{e}_n$. Details about such an OPM approach for the stochastic context will be communicated elsewhere.

\section*{Acknowledgments}
This work has been supported by the Office of Naval Research (ONR) Multidisciplinary University Research Initiative (MURI) grant N00014-20-1-2023, and by the National Science Foundation grant DMS-2108856. This study was also supported by the European Research Council (ERC) under the European Union's Horizon 2020 research and innovation program [Grant Agreement No.~810370], the Ben May Center grant for theoretical and/or computational research and by the Israeli Council for Higher Education (CHE) via the Weizmann Data Science Research Center.

\appendix

\section{$M$-term: First moments, non-Gaussian statistics, and large probability estimates}\label{Appendix_A}

\subsection{First moments} \label{Sec_Mn_first_moments}
We consider the random variable of the form
\be \label{Mn recall2}
M(\omega) = \int_{-\infty}^0 e^{g(\lambda) s + \sigma W_s(\omega)}\, \mathrm{d}s,
\ee
in which $\sigma$ is allowed to take negative values and $g(\lambda)$ is a function of the control parameter $\lambda$. Such random variables (with negative $\sigma$) arise for instance 
when the stable modes are forced stochastically but not the unstable ones; see Eq.~\eqref{Eq_Mn_formula_general}.

The following lemma provides exact formulas for the expectation, variance and autocorrelation function, of such a general $M$-term.
\bl  \label{lem:Mn}  
Let $g>0$ and 
\be \label{eq:g and sigma}
 \qquad  \sigma_*= \sqrt{2 g}, \qquad \sigma_{\#} = \sqrt{g}.
\ee

Then $M$-term defined by \eqref{Mn recall2} a wide-sense stationary random process provided that $|\sigma| < \sigma_{\#}$, and 
furthermore
\bi
\item[(i)] The expectation of $M$ exists if and only if $|\sigma| < \sigma_{*}$, and is given by
\be \label{expectation Mn}
\mathbb{E}(M) = \frac{2}{2g - \sigma^2}, \quad  |\sigma |< \sigma_{*}.
\ee 

\item[(ii)] The variance of $M$ exists if and only if $|\sigma |< \sigma_{\#}$, and is given  by
\be\label{variance Mn}
\mathrm{Var}(M)  = \frac{2 \sigma^2}{(2g - \sigma^2 )^2(g -  \sigma^2)}.
\ee
\ei
Finally, the autocorrelation $R(t)$ of the stochastic process
\be
M(t,\omega) = \int_{-\infty}^0 e^{g s + \sigma (W_{s+t}(\omega) - W_t(\omega))}\, \d s,
\ee
exists if and only if  $|\sigma| < \sigma_{\#}$, and is given by
\be \label{eq:autocorrelation}
R(t)= \exp\Bigl( -\Bigl ( g - \frac{\sigma^2}{2} \Bigr ) |t|\Bigr), \quad t \in \mathbb{R}.
\ee

\el

This lemma results from direct application of the  Fubini Theorem, the independent increment property of the Wiener process, and the fact that $\mathbb{E}(e^{\sigma W_t(\cdot)}) = e^{\sigma^2|t|/2}$ for any $t\in \mathbb{R}$, as expectation of the geometric Brownian motion generated by $\d S_t = \frac{\sigma^2}{2} S_t \d t + \sigma S_t \d W_t$. See the Appendix of \cite{CLW15_vol2} for a proof.

\subsection{Non-Gaussian statistics and log-normal approximations}\label{Sec_lognormal}
We provide here a simple analytic understanding regarding the non-Gaussian statistics followed by an arbitrary $M$-term, in the case of real eigenvalues. Recall that an $M$-term is obtained as the stationary solution of the scalar SDE, 
\bea \label{M_eq_v3} 
\mathrm{d} M= \left(1 - g M \right )\mathrm{d} t -   \sigma M \circ \mathrm{d} W_t,
\eea
where $g$ denotes a distance to resonance, namely $g=(\boldsymbol{k},\boldsymbol{\beta}_c(\lambda)) - \beta_n(\lambda)$ for some $n$.  

Rescale \eqref{M_eq_v3} by using 
\be
W_t(\omega) = \frac{1}{\sqrt{g}}\widetilde{W}_{g t}(\omega), \quad \widetilde{t} = g t, \quad \widetilde{M} = g M.
\ee
Then  \eqref{M_eq_v3} becomes
\be \label{Eq_Mn_v2}
\d \widetilde{M} = (1 - \widetilde{M}) \d \widetilde{t} - \frac{\sigma}{\sqrt{g}} \widetilde{M} \circ \d \widetilde{W}_{\widetilde{t}}. 
\ee 

To have an explicit formula for the probability distribution function (PDF) of the random variable $\widetilde{M}$ is non-trivial. We illustrate nevertheless that this PDF is well approximated by 
the PDF of an analytic, log-normal distributed random variable $N$, provided that  $\sigma < \sqrt{g}$. 

\begin{figure}
\centering
\includegraphics[width=1\linewidth, height=.4\textwidth]{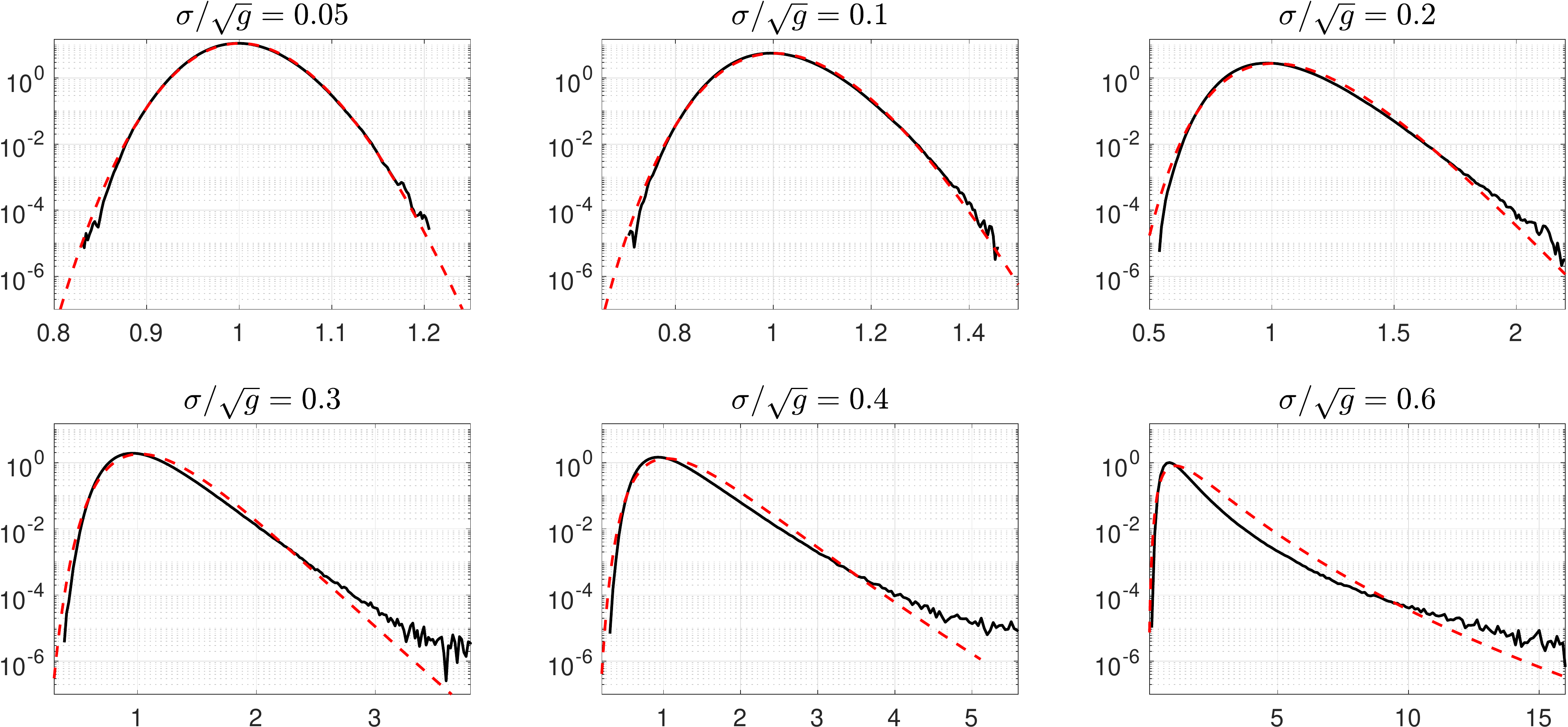}
\vspace{-1ex}
\caption{\footnotesize PDF of $\widetilde{M}$ (black) vs.~its lognormal approximation given by \eqref{Eq_PDF_LN} (red dash curve), for various values of $\sigma/\sqrt{g}$. A semi-log scale is used to plot these PDFs.  Note that $\sigma/\sqrt{g} < 1$ is required to ensure the mean and variance of $\widetilde{M}$ to exist; see Lemma \ref{lem:Mn}.}\label{Fig_LN_approx}
\end{figure}

Our analysis revealed that the parameters of this log-normal distribution are 
\bea
& \mu_{\sigma,g} = \ln(\kappa_{\sigma,g}) -\frac{1}{2} (\Sigma_{\sigma,g})^2,\\
& \Sigma_{\sigma,g} = \sqrt{\ln\left( 1+\left(\frac{ \eta_{\sigma,g}}{\kappa_{\sigma,g}}\right)^2 \right)}, 
\eea
where
\bea
& \kappa_{\sigma,g}= \frac{2}{2 - \sigma^2/g}, \\
& \eta_{\sigma,g} =  \sqrt{\frac{2 \sigma^2/g}{(2 - \sigma^2/g )^2(1 -  \sigma^2/g)}}.
\eea
The PDF of $N$  is then given by (\cite[Chap.~14]{Johnson_al94})
\be \label{Eq_PDF_LN}
p(N) = \frac{1}{N  \Sigma_{\sigma,g}\sqrt{2\pi}} \exp\left(-\frac{(\ln(N) - \mu_{\sigma,g})^2}{2\Sigma_{\sigma,g}^2}  \right),  \quad N > 0,
\ee
and turns out to approximate robustly the distribution of simulated $M$-term for different values of $\sigma/\sqrt{g}$ as shown in Fig.~\ref{Fig_LN_approx}.

\subsection{Probability estimates in the large } \label{Sec_proof_Mn_estimates0}

\bl \label{Lemma_Mn_bound}
We assume that Conditions~\ref{Cond_A1} and \ref{Cond_A2} as well as \eqref{cond_interaction_table} of Sec.~\ref{Sec_modeling_error} hold. Consider the $M_{11}^{n,\lambda}$-terms appearing in \eqref{Eq_Phi_Burgers_comp} for $n = 2, \cdots, N$ and $\lambda \in [\lambda_c, \lambda^*]$, where $N$ is given by \eqref{cond_interaction_table} and $\lambda^*$ is given in Condition~\ref{Cond_A2}. Then, for any $\chi$ in $(0, 1)$, there exists $\kappa > 0$ such that 
\be \label{Eq_Mn_bounds}
\mathbb{P}\left \{ 0 \le  M_{11}^{n,\lambda}(\omega)  < \kappa \; \big | \;  n = 2, \cdots, N, \lambda \in [\lambda_c, \lambda^*] \right \} \ge 1- \chi.
\ee
\el

\bp
The result is a direct consequence of the Chebyshev's inequality, which states that if $Y$ is a random variable and $p \ge 1$, then for any $\kappa >0$, it holds that 
\be \label{Chebyshev}
\mathbb{P}\{\omega: |Y| \ge \kappa \} \le \frac{1}{\kappa^p} \mathbb{E}(|Y|^p).
\ee
Since $M_{11}^{n,\lambda}(\omega) = \int_{-\infty}^0  e^{(2\epsilon - \beta_n(\lambda)) s + \sigma W_s(\omega)} \mathrm{d}s$ and the eigenvalues are all real-valued and are arranged in descending order, we have for all $\omega \in \Omega, \lambda \in[\lambda_c, \lambda^*]$, 
\be \label{Eq_Mn_order1}
0\le M_{11}^{N,\lambda}(0,\omega) \le M_{11}^{N-1,\lambda}(0,\omega) \le \cdots \le M_{11}^{3,\lambda}(0,\omega) \le M_{11}^{2,\lambda}(0,\omega).
\ee
Denote 
\be
\lambda^{\sharp} = \mathop{\mathrm{argmax}}_{\lambda \in [\lambda_c, \lambda^*]} \beta_2(\lambda).
\ee
We have then
\be \label{Eq_Mn_order2}
M_{11}^{2,\lambda}(\omega) \le  M_{11}^{2,\lambda^\sharp}(\omega), \quad \forall  \omega \in \Omega, \lambda \in[\lambda_c, \lambda^*].
\ee

Recall that due to Lemma~\ref{lem:Mn}
\be
\mathbb{E}(M_{11}^{2,\lambda^\sharp}) = \frac{2}{2(2\epsilon-\beta_2(\lambda^\sharp)) - \sigma^2},
\ee 
which holds for all $\sigma  < \sqrt{2(2\epsilon-\beta_2(\lambda^\sharp))}$. Since $\sigma = \sqrt{\epsilon}$ and $\beta_2(\lambda^\sharp) < 0$ thanks to Condition~\ref{Cond_A2}, this desired condition $\sigma  < \sqrt{2(2\epsilon-\beta_2(\lambda^\sharp))}$ is satisfied. Thus, for any $\chi$ in $(0,1)$, thanks to \eqref{Eq_Mn_order1} and \eqref{Eq_Mn_order2}, we just need to apply \eqref{Chebyshev} with $Y = M_{11}^{2,\lambda^\sharp}$, $p=1$, and 
\be \label{Eq_kappa}
\kappa = \frac{1}{|\beta_2(\lambda^\sharp)| \chi }, 
\ee
to ensure \eqref{Eq_Mn_bounds}.
\ep

\bl \label{Lemma_BM_bound}
For any $\epsilon > 0$, $T>0$ and $\chi$ in $(0, 1]$, let $\gamma = \sqrt{-2T\ln(1-\chi))}$, then 
\be \label{Eq_BM_bound}
\mathbb{P}\left \{ \sqrt{\epsilon} \sup_{0\le t \le T/\epsilon}|W_t(\omega)| \le \gamma \right \} \ge 1 - \chi.
\ee
\el

\bp
Recall that for any $\gamma>0$ and $T^*>0$, it holds that
\be 
\mathbb{P}\{\sup_{0 \le t \le T^*}|W_t(\omega)| \ge \gamma\} \le 2 \mathbb{P}\{|W_{T^*}(\omega)|\ge \gamma\};
\ee 
see e.g. Billingsley \cite{Bil86} page 529. So for any $\sigma > 0$, we have
\be \label{est. W}
\mathbb{P}\bigg\{ \sigma  \sup_{0 \le  t \le T^*}|W_t(\omega)|  \ge \gamma \bigg \} \le 2 \mathbb{P}\left\{ \sigma  |W_{T^*}(\omega)| \ge \gamma \right \}  \le  \exp \left(-\frac{\gamma^2}{2\sigma^2 T^*} \right), 
\ee
where in the last inequality above we used the facts that $W_{T^*}(\omega)$ follows a normal distribution with mean zero and variance equal to $T^*$ and that 
\be
\int_a^\infty e^{-x^2} \d x \le \frac{\sqrt{\pi}}{2} e^{-a^2}. 
\ee  
Hence, with $T^* = T/\epsilon$ and $\sigma = \sqrt{\epsilon}$, we have 
\be
\mathbb{P}\left \{ \sqrt{\epsilon} \sup_{0\le t \le T/\epsilon}|W_t(\omega)| \ge \gamma \right \} \le \exp \left(-\frac{\gamma^2}{2 T} \right).
\ee
Then, by choosing $\gamma = \sqrt{-2T\ln(1-\chi))}$, we get the desired estimate \eqref{Eq_BM_bound}. 
\ep

\section{A priori estimates in large probability: Surrogate system~\eqref{Eq_X_Phi}}\label{Sec_Apriori_est}

\bp[Proof of Lemma \ref{Lem_apriori_est}]

The desired estimate \eqref{A_priori_goal} follows from (i) basic probabilistic estimates about  the $M$-terms and the Brownian motion (see Lemmas \ref{Lemma_Mn_bound} and \ref{Lemma_BM_bound} above), and (ii) elementary energy estimates based on the dissipation condition \eqref{Cond_supercritical}. 

To perform (and simplify) our estimates we first transform the SDE system \eqref{Eq_X_Phi} into differential equations with random coefficients by using $Z(t,\omega)$, the Ornstein-Uhlenbeck process, stationary solution of the scalar Langevin equation \eqref{Eq_Langevin}. 

More exactly, we perform the change of variables
\be \label{Eq_random_transform}
U = e^{-Z(t,\omega)} X, \quad \Psi = e^{-Z(t,\omega)} \Phi.
\ee
The system \eqref{Eq_X_Phi} is then transformed into 
\begin{subequations} \label{Eq_X_Phi_v2}
\begin{align}
\frac{\mathrm{d} U}{\d t} &= \epsilon U + Z(t,\omega) U + e^{Z(t,\omega)}  U \mathcal{F}_\c( \Psi ),  \label{Eq_X_Phi_v2_a} \\
\frac{\d \Psi}{\d t} & = L^{\s}_\lambda \Psi +  Z(t,\omega) \Psi + e^{Z(t,\omega)} B^\s_{11}U^2 + 2 e^{Z(t,\omega)} \Psi \mathcal{F}_\c( \Psi ).   \label{Eq_X_Phi_v2_b}
\end{align}
\end{subequations}
The estimation of $X$ and $\Phi$ based on \eqref{Eq_X_Phi_v2} are organized into three steps below. 

\medskip
\noindent{\bf Step 1}: Energy estimates for $X$. Multiplying \eqref{Eq_X_Phi_v2_a} by $U$, we get 
\bea \label{Eq_X_Phi_v3}
\frac{1}{2} \frac{\d }{\d t} U^2 = \epsilon U^2  + Z(t,\omega) U^2 + e^{Z(t,\omega)}  U^2\mathcal{F}_\c( \Psi ). 
\eea

Observe that by \eqref{Eq_Phi_Burgers} and \eqref{Def_Fc} we have
\bea \label{Eq_Fc_est}
\mathcal{F}_\c( \Psi ) & =  e^{-Z(t,\omega)}  \sum_{n=2}^N \Phi_n \big( B_{1n}^1 + B_{n1}^1\big) \\
& = e^{-Z(t,\omega)} X^2 \sum_{n=2}^N B_{11}^n M_{11}^{n,\lambda} \big( B_{1n}^1 + B_{n1}^1\big) \\
& \le 0. \qquad \text{(due to \eqref{Cond_supercritical})}
\eea
We get then from \eqref{Eq_X_Phi_v3} that 
\bea
\frac{1}{2} \frac{\d }{\d t} U^2 \le \epsilon U^2 + Z(t,\omega) U^2,
\eea
which leads to
\be
U^2(t,\omega) \le e^{\alpha(t,\omega)} U^2_0,
\ee
with $\alpha(t,\omega)=2\epsilon t + 2 \int_0^t Z(s,\omega) \d s.$

Transforming back to the $X$-variable by using \eqref{Eq_random_transform}, and by noting that 
\be \label{OU_identity}
Z(t,\omega) = Z(0,\omega) - \int_0^t Z(s, \omega) \d s + \sqrt{\epsilon} W_t(\omega),
\ee 
we obtain
\be \label{Eq_apriori_est1}
X^2(t,\omega)  \le e^{2\epsilon t + 2 \sqrt{\epsilon} W_t(\omega)} X^2_0. 
\ee

\medskip
\noindent{\bf Step 2}: Energy estimates for $\Phi$. Note that $\langle L^{\s}_\lambda \Psi,  \Psi \rangle \le \beta_2(\lambda) \|\Psi\|^2$. Taking the inner-product of \eqref{Eq_X_Phi_v2_b} with $\Psi$ and using again the fact that $\mathcal{F}_\c( \Psi ) \le 0$ (cf.~\eqref{Eq_Fc_est}), we get 
 \be  \label{Eq_Phi_est1}
 \frac{1}{2}\frac{\d \|\Psi\|^2}{\d t} \le  \beta_2(\lambda) \|\Psi\|^2 +  Z(t,\omega) \|\Psi\|^2 + e^{Z(t,\omega)}U^2 \langle B^\s_{11}, \Psi \rangle.
\ee
Now, note that
\be
I \stackrel{\textrm{def}}{=} e^{Z(t,\omega)}U^2 \langle B^\s_{11}, \Psi \rangle = e^{Z(t,\omega)}U^2 \langle B(\boldsymbol{e}_1,\boldsymbol{e}_1), \Psi \rangle  \le e^{Z(t,\omega)} U^2\|B(\boldsymbol{e}_1,\boldsymbol{e}_1)\| \|\Psi\|,
\ee
which leads to 
\be
I \le  \frac{1}{2 |\beta_2(\lambda)|} e^{2Z(t,\omega)} \|B(\boldsymbol{e}_1,\boldsymbol{e}_1)\|^2 U^4  + \frac{1}{2}  |\beta_2(\lambda)|  \|\Psi\|^2.
\ee
Using this last inequality in \eqref{Eq_Phi_est1}, we arrive at
\be  \label{Eq_Phi_est2}
 \frac{1}{2}\frac{\d \|\Psi\|^2}{\d t} \le  \frac{1}{2} \beta_2(\lambda) \|\Psi\|^2 +  Z(t,\omega) \|\Psi\|^2 + \frac{1}{2 |\beta_2(\lambda)|} e^{2Z(t,\omega)} \|B(\boldsymbol{e}_1,\boldsymbol{e}_1)\|^2 U^4,
\ee
which leads, after integration, to
\be
 \|\Psi(t,\omega)\|^2 \le  e^{\beta_2(\lambda)t  + 2 \int_0^t  Z(s, \omega) \d s} \|\Psi (0,\omega)\|^2 +  \frac{1}{|\beta_2(\lambda)|} \|B(\boldsymbol{e}_1,\boldsymbol{e}_1)\|^2  J(t,\omega),
 \ee
 with
 \be
J(t,\omega)= \int_0^t e^{\beta_2(\lambda) (t-s)+ 2\int_s^t  Z(\tau, \omega) \d \tau}  e^{2 Z(s, \omega)}  U^4(s,\omega) \d s. 
 \ee
 
By transforming back to the $(X, \Phi)$-variable using \eqref{Eq_random_transform} and taking into account the identity \eqref{OU_identity}, we obtain
 \bea \label{Eq_Phi_est3}
\|\Phi(t,\omega)\|^2 & \le e^{\beta_2(\lambda) t + 2\sqrt{\epsilon} W_t(\omega)} \|\Phi(0,\omega)\|^2 + \frac{1}{|\beta_2(\lambda)|} \|B(\boldsymbol{e}_1,\boldsymbol{e}_1)\|^2 \widetilde{J}(t,\omega),
 \eea
 with 
 \be \label{Def_J_tilde}
 \widetilde{J}(t,\omega)=\int_0^t e^{\beta_2(\lambda) (t-s) + 2\sqrt{\epsilon} (W_t(\omega) - W_s(\omega))} X^4(s,\omega) \d s. 
 \ee

 \medskip
\noindent{\bf Step 3}: Probabilistic estimates for $X$ and $\Phi$. Based on \eqref{Eq_apriori_est1} and \eqref{Eq_Phi_est3}, to derive \eqref{A_priori_goal}, we just need to obtain suitable bounds for $W_t$ and the random coefficients $M_{11}^{n,\lambda}(0,\omega)$'s appearing in the expression of $\Phi(0,\omega)$; cf.~\eqref{Eq_Phi_IC}.  More precisely, for $T$ and $\chi$ as given in the statement of Lemma~\ref{Lem_apriori_est}, we show that there exists a subset $\Omega^*$ of $\Omega$ with $\mathbb{P}(\Omega^*) \ge 1 - \chi$, over which $M_{11}^{n,\lambda}(0, \omega)$ is bounded for all $n = 1, \cdots, N$, and $\sqrt{\epsilon}W_t(\omega)$ is also bounded for all $\omega$ in $\Omega^*$ and $t$ in $[0, T/\epsilon]$. The existence of such a subset $\Omega^*$ is guaranteed by Lemma~\ref{Lemma_Mn_bound} and Lemma~\ref{Lemma_BM_bound}. 

Indeed, applying Lemma~\ref{Lemma_Mn_bound} with $\chi$ therein replaced by $\chi/2$, there exists a constant $\kappa > 0$ and a subset $\Omega_1$ with $\mathbb{P}(\Omega_1) \ge 1 - \chi/2$, for which  
\be  
0 \le  M_{11}^{n,\lambda}(0,\omega)  < \kappa, \;  \omega \in \Omega_1, n = 2, \cdots, N, \lambda \in [\lambda_c, \lambda^*].
\ee 
Applying Lemma~\ref{Lemma_BM_bound} with $\chi$ therein replaced by $\chi/2$, there exists a constant $\gamma > 0$ and a subset $\Omega_2$ with $\mathbb{P}(\Omega_2) \ge 1 - \chi/2$, for which  
\be
\sqrt{\epsilon} \sup_{0 \le t \le T/\epsilon}|W_t(\omega)| \le \gamma, \; \omega \in \Omega_2. 
\ee

Now, take $\Omega^* = \Omega_1 \cap \Omega_2$. We have $\mathbb{P}(\Omega^*) \ge 1 - \chi$, and 
\begin{subequations} \label{Eq_apriori_est2}
\begin{align}
& \sqrt{\epsilon} \sup_{0 \le t \le T/\epsilon}|W_t(\omega)| \le \gamma,\;  \omega \in \Omega^*, \label{Eq_apriori_est2_a} \\
& 0 \le M_{11}^{n,\lambda}(0,\omega)  < \kappa, \;  \omega \in \Omega_1, n = 2, \cdots, N, \lambda \in [\lambda_c, \lambda^*].\label{Eq_apriori_est2_b}
\end{align}
\end{subequations}

Using \eqref{Eq_apriori_est2_a} and \eqref{Eq_X_IC}, we conclude from \eqref{Eq_apriori_est1} that 
 \be \label{Eq_apriori_est3}
 X^2(t,\omega) \le (C_1)^2 e^{2 T + 2 \gamma} \epsilon, \quad t \in [0, T/\epsilon],  \; \omega  \in \Omega^*,
 \ee
where $C_1 > 0$ is a constant such that $|X_0| \leq C_1 \sqrt{\epsilon}$ (Due to \eqref{Eq_X_IC}). 

For $\Phi$, the terms on the RHS of \eqref{Eq_Phi_est3} can be estimated as follows. Using the definition of $\Phi(0,\omega)$ in \eqref{Eq_Phi_IC}, the assumption on $X_0$ in \eqref{Eq_X_IC}, and the probabilistic bound \eqref{Eq_apriori_est2_b} about the $M_{11}^{n,\lambda}$-terms, we get
\be \label{Eq_Phi_est4}
e^{\beta_2(\lambda) t + 2\sqrt{\epsilon} W_t(\omega)} \|\Phi(0,\omega)\|^2 \le e^{2\gamma} \|\Phi(0,\omega)\|^2 \le \kappa^2 (C_1)^4 e^{2\gamma}  \left(\sum_{n=2}^N (B_{11}^n)^2\right) \epsilon^2, 
\ee
 for any $t$ in $[0, T/\epsilon]$, and  $\omega$ in $\Omega^*$, where $C_1$ is the same as in \eqref{Eq_apriori_est3}.

For $\widetilde{J}$ defined by \eqref{Def_J_tilde}, using \eqref{Eq_apriori_est2_a} and \eqref{Eq_apriori_est3}, we obtain that for all $\omega$ in $\Omega^*$ and $t$ in $[0, T/\epsilon]$,
\be \label{Eq_Phi_est5}
\widetilde{J}(t,\omega)\le \int_0^t e^{\beta_2(\lambda) (t-s) + 4 \gamma} \d s \,
(C_1)^4 e^{4 T + 4 \gamma} \epsilon^2 \le \frac{1}{|\beta_2(\lambda)|} (C_1)^4 e^{4 T + 8 \gamma} \epsilon^2.
\ee
Now, by using \eqref{Eq_Phi_est4} and \eqref{Eq_Phi_est5}  in \eqref{Eq_Phi_est3}, we arrive, for any $t$ in $[0, T/\epsilon]$, and  $\omega$ in $\Omega^*$, at
\be \label{Eq_Phi_est6}
\|\Phi(t,\omega)\| \le C_2\, \epsilon, 
\ee
in which $C_2$ is a constant independent of $\epsilon$. The desired estimate \eqref{A_priori_goal} follows now from \eqref{Eq_apriori_est3} and \eqref{Eq_Phi_est6} by taking $C = \max\{C_1 e^{T + \gamma}, C_2\}$.
\ep


\section{High-mode estimates in large probability} \label{Sec_Apriori_est_SPDE}

We provide in this section the proofs of the two key estimates \eqref{Key_est1} and \eqref{Key_est2}, regarding 
the high-mode component, $u_\s$,  of the solution, $u$, to Eq.~\eqref{Eq_Quadratic}. 
\subsection{High-mode estimate  \eqref{Key_est1}}\label{Sec_Apriori_est_SPDE_key1}

First note that,  since the nonlinearity $B$ in Eq.~\eqref{Eq_Quadratic} is bilinear, the transformed equation \eqref{REE} for the variable $v = e^{-Z(t,\omega)} u$ associated with Eq.~\eqref{Eq_Quadratic} becomes here 
\begin{equation} \label{REE_bilinear} 
\frac{\mathrm{d} v}{\mathrm{d} t} = L_\lambda v + Z(t,\omega) v + e^{Z(t,\omega)} B(v),
\end{equation}
with $Z(t,\omega)$ denoting again the OU process, stationary solution of the scalar Langevin equation \eqref{Eq_Langevin}. 

Taking inner product of \eqref{REE_bilinear} with $v$ and using the energy conservation property of $B$ (cf.~\eqref{Eq_energy_conservation}), we obtain
\be
\frac{1}{2}\frac{\d \|v\|^2}{\d t} = \langle L_\lambda v, v\rangle  +  Z(t,\omega) \|v_{\s}\|^2 \le \epsilon \|v_{\s}\|^2 + Z(t,\omega)\|v_{\s}\|^2. 
\ee
Integrating this last equation, then transforming back to the $u$-variable and taking into consideration the identity \eqref{OU_identity}, we obtain
\be \label{Eq_apriori_est1_SPDE}
x^2(t,\omega) + \|u_{\s}(t,\omega)\|^2 \le e^{2\epsilon t + 2 \sqrt{\epsilon} W_t(\omega)} \|u(0,\omega)\|^2.
\ee
Recall also that, thanks to Lemma~\ref{Lemma_BM_bound}, for any $\chi$ in $(0,1)$, there exists a subset $\Omega^* \subset \Omega$ with $\mathbb{P}(\Omega^*) \ge 1 - \chi$ and a constant $\gamma > 0$, such that 
\be  \label{Eq_est_BM_recall}
\sqrt{\epsilon} \sup_{0 \le t \le T/\epsilon}|W_t(\omega)| \le \gamma,\;  \omega \in \Omega^*;
\ee
see \eqref{Eq_apriori_est2_a}. Thus, if $|u(0)| \sim \sqrt{\epsilon}$, then thanks to \eqref{Eq_apriori_est1_SPDE} and \eqref{Eq_est_BM_recall}, there exists a constant $C>0$ independent of $\epsilon$ such that 
\be \label{A_priori_SPDE_v1}
|x(t,\omega)|^2 + \|u_{\s}(t,\omega) \|^2 \le C \epsilon, \;\; \forall t \in [0, T/\epsilon], \; \omega \in \Omega^*.
\ee

By introducing 
\be \label{Eq_random_transform2}
y = e^{-Z(t,\omega)} x, \quad v_{\s} = e^{-Z(t,\omega)} u_{\s},
\ee
we get from \eqref{Eq_Quadratic_v2_b} the following equation for $v_{\s}$:
\be \label{Eq_us_b}
\frac{\d v_{\s} }{\d t} = L^{\s}_\lambda v_{\s}  +  Z(t,\omega)v_{\s} + e^{Z(t,\omega)} \big( y^2  B^{\s}_{11} + y \mathcal{F}_{\s}(v_{\s}) + \Pi_{
\s} B(v_{\s},v_{\s}) \big).   
\ee
Noting that $\langle B(v_{\s},v_{\s}), v_{\s} \rangle = 0$ and that $\langle L^{\s}_\lambda v_{\s}, v_{\s} \rangle = \langle L_\lambda v_{\s}, v_{\s} \rangle$, we get then
\be \label{Eq_energy_est_us}
\frac{1}{2}\frac{\d \|v_{\s}\|^2 }{\d t} = \langle L_\lambda v_{\s}, v_{\s}\rangle  +  Z(t,\omega) \|v_{\s}\|^2  + e^{Z(t,\omega)} \Big \langle y^2  B^{\s}_{11} + y \mathcal{F}_{\s}(v_{\s}), v_{\s} \Big \rangle. 
\ee

For $\langle L_\lambda v_{\s}, v_{\s}\rangle$, first note that by Condition~\ref{Cond_L1}, we have 
\be \label{Eq_div_thm}
\langle A u, u \rangle =  \|u\|_{\boldsymbol{\nu}}^2, \quad u \in \mathcal{W}=D(A).
\ee
Now, let $\eta$ be as given by Condition~\ref{Cond_L2}. Since $L_\lambda = - A + P_\lambda$,  we get, thanks to \eqref{Eq_P_control_a} and \eqref{Eq_div_thm}
\bea \label{Eq_damping_term_est}
\langle L_\lambda v_{\s}, v_{\s}\rangle &= 
\langle (-\eta A + P_\lambda) v_{\s}, v_{\s}\rangle - (1-\eta) \langle A v_{\s}, v_{\s}\rangle \\
& \le - \delta \|v_{\s}\|^2 - (1-\eta) \|v_{\s}\|_{\boldsymbol{\nu}}^2.
\eea

Recall that $B^{\s}_{11} = \Pi_{\s} B(\boldsymbol{e}_1,  \boldsymbol{e}_1)$; see \eqref{Def_B_11s}. We get
\bea \label{Eq_B_term_est1}
e^{Z(t,\omega)}  \langle y^2 B^{\s}_{11}, v_{\s} \rangle  & \le e^{Z(t,\omega)}   C_B  y^2 \|\boldsymbol{e}_1\|_V^2 \|v_{\s}\| \\
& \le \frac{e^{2Z(t,\omega)} (C_B)^2 \|\boldsymbol{e}_1\|_V^4}{2 \delta} y^4  +  \frac{\delta}{2}\|v_{\s}\|^2,
\eea
where $C_B > 0$ denotes the smallest constant for which $\|B(u,v)\| \le C_B \|u\|_V \|v\|_V$ for all $u$ and $v$ in $V$.

Recall also that $\mathcal{F}_\s(v) =\Pi_\s (B(\boldsymbol{e}_1,v) +B(v,\boldsymbol{e}_1))$; see \eqref{Eq_Fs}. We get then  
\bea \label{Eq_B_term_est2a}
e^{Z(t,\omega)} \langle y  \mathcal{F}_\s(v_{\s}), v_{\s} \rangle & \le 2e^{Z(t,\omega)} C_B |y| \|\boldsymbol{e}_1\|_V \|v_{\s}\|_V \|v_{\s}\| \\
& \le 2e^{Z(t,\omega)} \mathfrak{C} C_B |y| \|\boldsymbol{e}_1\|_V \|v_{\s}\|_{\boldsymbol{\nu}} \|v_{\s}\|,
\eea
where $\mathfrak{C} > 0$ denotes the smallest constant for which $ \|u\|_V \le \mathfrak{C} \|u\|_{\boldsymbol{\nu}}$ for all $u$ in $V$. Such a constant exists due to the assumption that $\|\cdot \|_{\boldsymbol{\nu}}$ defines an equivalent norm to $\|\cdot\|_V$ on $V$;  see Condition~\ref{Cond_L1}. We get then 
\be \label{Eq_B_term_est2}
e^{Z(t,\omega)}   \langle y \mathcal{F}_\s(v_{\s}), v_{\s} \rangle \le \frac{1-\eta}{2}  \|v_{\s}\|^2_{\boldsymbol{\nu}} + \frac{2e^{2Z(t,\omega)} (\mathfrak{C} C_B \|\boldsymbol{e}_1\|_V)^2}{(1-\eta)} |y|^2  \|v_{\s}\|^2. 
\ee

Using \eqref{Eq_damping_term_est}, \eqref{Eq_B_term_est1}, and \eqref{Eq_B_term_est2} in \eqref{Eq_energy_est_us}, we obtain
\bea \label{Eq_energy_est_us_v2}
\frac{1}{2}\frac{\d \|v_{\s}\|^2 }{\d t}  & \le -\frac{1}{2} \delta \|v_{\s}\|^2 - \frac{(1-\eta)}{2} \|v_{\s}\|^2_{\boldsymbol{\nu}} +  Z(t,\omega) \|v_{\s}\|^2  \\
& \quad + \frac{e^{2Z(t,\omega)} (C_B)^2 \|\boldsymbol{e}_1\|_V^4}{2\delta} y^4 + \frac{2 e^{2Z(t,\omega)} (\mathfrak{C} C_B \|\boldsymbol{e}_1\|_V)^2}{(1-\eta)} |y|^2  \|v_{\s}\|^2. 
\eea

By introducing 
\be
C_1 = \frac{(C_B)^2 \|\boldsymbol{e}_1\|_V^4}{2 \delta}, \quad C_2 =  \frac{2(\mathfrak{C} C_B \|\boldsymbol{e}_1\|_V)^2}{(1-\eta)}, 
\ee
and restricting to the subset $\Omega^*$, we have thanks to \eqref{A_priori_SPDE_v1} and \eqref{Eq_random_transform2} that 
\be \label{Eq_energy_est_us_v4}
C_1 e^{2Z(t,\omega)}y^4 +  C_2 e^{2Z(t,\omega)}  |y|^2  \|v_{\s}\|^2 \le C e^{-2Z(t,\omega)} \epsilon^2, 
\ee
for $\omega \in \Omega^*$, $t \in [0, T/\epsilon]$, and some generic constant $C>0$ independent of $\epsilon$.

Applying Gronwall's inequality to \eqref{Eq_energy_est_us_v2} while taking into consideration \eqref{Eq_energy_est_us_v4}, we get 
\be \label{Eq_energy_est_us_v5}
\|v_{\s}(t,\omega) \|^2  \le e^{- \delta t + 2\int_0^t Z(s, \omega) \d s} \|v_{\s}(0,\omega)\|^2 + C \epsilon^2 \int_0^t 
e^{- \delta (t-s) +  2\int_s^t Z(\tau, \omega)  \d \tau - 2  Z(s, \omega)} \d s,
\ee
which holds for $\omega \in \Omega^*$ and $t \in [0, T/\epsilon]$. 

Converting back to $u_{\s}$ and using \eqref{Eq_est_BM_recall}, we get 
\bea \label{Eq_energy_est_us_v6}
\|u_{\s}(t,\omega)\|^2  & \le e^{- \delta t + 2 \sqrt{\epsilon} W_t(\omega)} \|u_{\s}(0,\omega)\|^2 + C \epsilon^2 \int_0^t 
e^{- \delta (t-s) + 2 \sqrt{\epsilon} (W_t(\omega) - W_s(\omega))} \d s, \\
& \le e^{- \delta t + 2 \gamma} \|u_{\s}(0,\omega)\|^2  + C  e^{4 \gamma} \epsilon^2 \int_0^t 
e^{- \delta (t-s)} \d s\\
& \le e^{2\gamma}  \|u_{\s}(0,\omega)\|^2 + \frac{C e^{4 \gamma}}{\delta}\epsilon^2, 
\eea
which holds for $\omega \in \Omega^*$ and $t \in [0, T/\epsilon]$. The desired estimate \eqref{Key_est1} is thus derived.

\subsection{High-mode estimate \eqref{Key_est2}}\label{Sec_Apriori_est_SPDE_key2}

In the following, to simplify the presentation, the estimates presented below are articulated at a formal level. The equations/inequations involving $\d \|v_{\s}\|^2_{\boldsymbol{\nu}}/\d t$ can be made rigorous by working with the corresponding integral formulation allowed by the regularity of $v$ recalled in \eqref{Eq_regularity}.

Back to \eqref{Eq_us_b}, multiplying both sides by $A v_{\s}$ and noting that $\langle v_{\s}, A v_{\s} \rangle = \|v_{\s}\|^2_{\boldsymbol{\nu}}$ (cf.~\eqref{Eq_div_thm}), we get
\bea \label{Eq_us_c}
\frac{1}{2}\frac{\d \|v_{\s}\|^2_{\boldsymbol{\nu}}}{\d t} & =  \langle L^{\s}_\lambda v_{\s}, A v_{\s}\rangle  +  Z(t,\omega) \|v_{\s}\|^2_{\boldsymbol{\nu}} \\
& \quad  + e^{Z(t,\omega)}  \left \langle \big( y^2  B^{\s}_{11} + y \mathcal{F}_{\s}(v_{\s}) +\Pi_{\s} B(v_{\s},v_{\s}) \big), A v_{\s} \right \rangle.   
\eea

For $\langle L^{\s}_\lambda v_{\s}, A v_{\s}\rangle$, using the assumption \eqref{Eq_P_control_b}, we get
\bea \label{Eq_damping_term_est2}
\langle L^{\s}_\lambda v_{\s}, A v_{\s}\rangle &= \langle L_\lambda v_{\s}, A v_{\s}\rangle \\
& =  \langle (-\eta A + P_\lambda) v_{\s}, A v_{\s}\rangle - (1-\eta) \langle A v_{\s}, A v_{\s}\rangle \\
& \le - \delta \|v_{\s}\|_{\boldsymbol{\nu}}^2  - (1-\eta) \| A v_{\s}\|^2.
\eea

We have also
\bea
e^{Z(t,\omega)} \langle y^2  B^{\s}_{11}, A v_{\s} \rangle &\le C_B e^{Z(t,\omega)} y^2 \|\boldsymbol{e}_1\|_V^2 \|A v_{\s}\| \\
& \le \frac{1}{1-\eta}(C_B)^2 e^{2Z(t,\omega)} y^4 \|\boldsymbol{e}_1\|_V^4 + \frac{1-\eta}{4} \|A v_{\s}\|^2.
\eea
By using the same type of estimates as given in \eqref{Eq_B_term_est2a}--\eqref{Eq_B_term_est2}, we get 
\be
e^{Z(t,\omega)} \langle y \mathcal{F}_{\s}(v_{\s}), A v_{\s} \rangle  \le \frac{2 (\mathfrak{C} C_B)^2}{1-\eta}  e^{2Z(t,\omega)} y^2 \|\boldsymbol{e}_1\|_V^2 \|v_{\s}\|_{\boldsymbol{\nu}}^2 + \frac{1-\eta}{2} \|A v_{\s}\|^2.
\ee
We also have
\bea
e^{Z(t,\omega)} \langle  B(v_{\s},v_{\s}), A v_{\s} \rangle & \le C_B e^{Z(t,\omega)} \|v_{\s}\|_V^2 \|A v_{\s}\| \\
& \le  \frac{1}{1-\eta}(C_B)^2 e^{2Z(t,\omega)} \|v_{\s}\|_V^4 + \frac{1-\eta}{4} \|A v_{\s}\|^2 \\
& \le  \frac{1}{1-\eta}\mathfrak{C}^4 (C_B)^2 e^{2Z(t,\omega)} \|v_{\s}\|_{\boldsymbol{\nu}}^4 + \frac{1-\eta}{4} \|A v_{\s}\|^2.
\eea

Using the above estimates in \eqref{Eq_us_c}, we get 
\be \label{Eq_us_d}
\frac{1}{2}\frac{\d \|v_{\s}\|^2_{\boldsymbol{\nu}}}{\d t}  \le -\delta \|v_{\s}\|^2_{\boldsymbol{\nu}}  +  Z(t,\omega) \|v_{\s}\|^2_{\boldsymbol{\nu}} + C e^{2Z(t,\omega)} (y^4  + y^2 \|v_{\s}\|^2_{\boldsymbol{\nu}} + \|v_{\s}\|^4_{\boldsymbol{\nu}}), 
\ee
where $C>0$ is a constant independent of $\epsilon$.

Now, using the estimate of $x$ given in \eqref{A_priori_SPDE_v1} and noting that $y = e^{-Z(t,\omega)} x$ from \eqref{Eq_random_transform2}, we obtain
\be \label{Eq_us_e}
\frac{1}{2}\frac{\d \|v_{\s}\|^2_{\boldsymbol{\nu}}}{\d t}  \le -\delta \|v_{\s}\|^2_{\boldsymbol{\nu}}  +  Z(t,\omega) \|v_{\s}\|^2_{\boldsymbol{\nu}} + C  \big( e^{-2Z(t,\omega)} \epsilon^2 + \epsilon \|v_{\s}\|^2_{\boldsymbol{\nu}} + e^{2Z(t,\omega)} \|v_{\s}\|^4_{\boldsymbol{\nu}} \big).
\ee
By choosing $\epsilon$ sufficiently small (achieved by choosing $\lambda^*$ sufficiently close to $\lambda_c$), we can ensure that 
\be
C \epsilon \|v_{\s}\|^2_{\boldsymbol{\nu}}  \le \frac{\delta}{2} \|v_{\s}\|^2_{\boldsymbol{\nu}}.
\ee
We get then 
\be \label{Eq_us_f}
\frac{1}{2}\frac{\d \|v_{\s}\|^2_{\boldsymbol{\nu}}}{\d t}  \le - \frac{\delta}{2} \|v_{\s}\|^2_{\boldsymbol{\nu}}  +  Z(t,\omega) \|v_{\s}\|^2_{\boldsymbol{\nu}} 
+ C e^{-2Z(t,\omega)} \epsilon^2 + C e^{2Z(t,\omega)} \|v_{\s}\|^4_{\boldsymbol{\nu}}.
\ee
For the remaining part of the proof, the generic constant $C>0$ (independent of $\epsilon$) is allowed to change in the course of the estimates.

Now, define
\be \label{Eq_transform_v}
q(t,\omega) = e^{\delta t - 2 \int_0^t Z(s, \omega) \d s}  \|v_{\s}\|^2_{\boldsymbol{\nu}},
\ee
we get from \eqref{Eq_us_f} that
\be \label{Eq_us_g}
\frac{\d q}{\d t} \le  C e^{\delta t - 2 \int_0^t Z(s, \omega) \d s -2Z(t,\omega)} \epsilon^2 + C \|u_{\s}\|^2_{\boldsymbol{\nu}}\, q(t,\omega).
\ee

To proceed further, we need an estimate on $\int_0^t \|u_{\s}\|^2_{\boldsymbol{\nu}} \d s$. For this purpose, 
first note that if $\|u_{\s}(0)\| \sim \epsilon$, then it follows from \eqref{Eq_energy_est_us_v6} that 
\be
\|u_{\s}(t,\omega) \| \le C \epsilon, \;\; \forall t \in [0, T/\epsilon], \; \omega \in \Omega^*.
\ee
Recall also from \eqref{A_priori_SPDE_v1} that $|x(t,\omega)| \le C \sqrt{\epsilon}$ for $\omega \in \Omega^*$ and $t \in [0, T/\epsilon]$. 
Now, by integrating \eqref{Eq_energy_est_us_v2} and using the above estimates for $|x(t,\omega)|$ and $\|u_{\s}(t,\omega)\|$, we get 
\bea \label{Eq_energy_est_us_v7}
\|v_{\s}(t,\omega)\|^2 & +  (1 - \eta) \int_0^t e^{2 \int_s^t Z(s, \omega) \d s} \|v_{\s}(s,\omega)\|^2_{\boldsymbol{\nu}} \d s \\
& \le  
e^{2\int_0^t Z(s, \omega) \d s} \|v_{\s}(0,\omega)\|^2  + C\epsilon^2 \int_0^t 
e^{2\int_s^t Z(\tau, \omega)  \d \tau - 2 Z(s, \omega)} \d s,
\eea
which holds for $\omega$ in $\Omega^*$ and $t$ in $[0, T/\epsilon]$. It follows that 
\bea
\int_0^t e^{2 \sqrt{\epsilon} (W_t(\omega)-W_s(\omega))} \|u_{\s}(s,\omega)\|^2_{\boldsymbol{\nu}}  \d s  & \le   e^{2 \sqrt{\epsilon} W_t(\omega)} \|u_{\s}(0,\omega)\|^2 + C \epsilon^2 \int_0^t 
e^{2\sqrt{\epsilon} (W_t(\omega)-W_s(\omega))} \d s  \\
& \le C \epsilon^2, \quad \forall \omega \in \Omega^*, t \in [0, T/\epsilon].
\eea
Using \eqref{Eq_est_BM_recall}, we get (after redefining $C$) that
\bea \label{Eq_int_u_est}
\int_0^t \|u_{\s}(s,\omega)\|^2_{\boldsymbol{\nu}} \d s  \le C \epsilon^2, \quad \forall \omega \in \Omega^*, t \in [0, T/\epsilon]. 
\eea

Now, let us denote $f(s,\omega) = \exp\big(\delta s - 2 \int_0^s Z(\tau, \omega) \d \tau -2Z(s, \omega)\big)$. Applying Gronwall's inequality to \eqref{Eq_us_g} and using \eqref{Eq_int_u_est}, we get 
\bea \label{Eq_us_h}
q(t,\omega) & \le  \exp\left(C \int_0^t  \|u_{\s}(s,\omega)\|^2_{\boldsymbol{\nu}} \d s\right) q(0,\omega) \\
& \quad + C\epsilon^2  \int_0^t \exp\left(C \int_s^t \|u_{\s}(\tau,\omega) \|^2_{\boldsymbol{\nu}}\d \tau\right) f(s,\omega) \d s\\
& \le C q(0,\omega) + C  \epsilon^2 \int_0^t f(s,\omega) \d s,
\eea
which holds for $\omega$  in $\Omega^*$ and $t$ in $[0, T/\epsilon]$. 

Now, converting back to $\|u_{\s}\|_V^2$ using \eqref{Eq_transform_v} and \eqref{Eq_random_transform2}, we get from \eqref{Eq_us_h} that 
\be \label{Eq_us_i}
\|u_{\s}(t,\omega)\|^2_{\boldsymbol{\nu}} \le C \|u_{\s}(0,\omega)\|^2_{\boldsymbol{\nu}} + C  \epsilon^2  \int_0^t e^{-\delta (t-s) + 2 \sqrt{\epsilon}(W_t(\omega) - W_s(\omega))} \d s,
\ee 
which holds for $\omega \in \Omega^*$ and $t \in [0, T/\epsilon]$. The desired estimate \eqref{Key_est2} follows then from \eqref{Eq_us_i} since $\|\cdot \|_{\boldsymbol{\nu}}$ defines an equivalent norm to $\|\cdot\|_V$ on $V$; see Condition~\ref{Cond_L1}.



\section{The non-Markovian normal form of a pitchfork bifurcation} \label{Sec_1DSDE_pitchfork}
We study here the following non-Markovian normal form of a supercritical pitchfork bifurcation
\be \label{Eq_reduced_generic}
\d X =( \epsilon X - \alpha M_\sigma(t, \omega) X^3)\d t + \sigma X \circ \d W_t,
\ee
where $\epsilon$ is the bifurcation parameter, $\alpha$ and $\sigma$ are positive constants, and $M_\sigma$ is the stationary process given by 
\be \label{Mt}
M_\sigma(t,\omega)  = \int_{-\infty}^0  e^{g s + \sigma (W_{s+t}(\omega) - W_t(\omega))} \mathrm{d}s, \; \omega \in \Omega, \; g>0. 
\ee

Note that since, $M_\sigma(t,\omega)$ is always positive, the above equation has a supercritical stochastic pitchfork bifurcation as $\epsilon$ crosses $0$ from below. In fact, by introducing a change of variables $\widetilde{X} = X^{-2}$, one can compute explicitly the bifurcated random steady states for \eqref{Eq_reduced_generic}. They are given by $\pm a_\epsilon$, where 
\be \label{Eq_stationary_generic}
a_\epsilon(t,\omega) = \frac{1}{\sqrt{2 \alpha \int_{-\infty}^t M_\sigma(s, \omega) \exp(f_\epsilon(t,s,\omega)) \d s}}, \quad \epsilon \ge 0,
\ee
with  
\be \label{Eq_forcing_in_steadystate}
f_\epsilon(t,s,\omega) = -2 \epsilon (t-s) - 2 \sigma (W_t(\omega) - W_s(\omega)).
\ee
In the following, we show that when $\sigma = \sqrt{\epsilon}$, then $a_\epsilon(t)$ is on the order of $\sqrt{\epsilon}$ with large probability for all $t\ge 0$. 

\bl \label{Lem_SS_estimation}
Consider the stationary solution $a_\epsilon$ for \eqref{Eq_reduced_generic} given by \eqref{Eq_stationary_generic}. Assume that $\sigma = \sqrt{\epsilon}$. Then, for any $\chi$ in $(0,1)$, there exist positive constants $c_1$ and $c_2$ independent of $\epsilon$ such that 
\be \label{steadystate_bound_goal}
\mathbb{P} \left(c_1\sqrt{\epsilon} \le a_\epsilon(t) \le c_2\sqrt{\epsilon} \right) \ge 1 - \chi, \quad \forall \epsilon \in [0, 2 g], t \ge 0.
\ee
\el

\bp
When $\epsilon = 0$, \eqref{steadystate_bound_goal} holds trivially since $a_\epsilon = 0$ since the integral term in the RHS of \eqref{Eq_stationary_generic} becomes $+\infty$ when $\epsilon = 0$. We assume thus $\epsilon > 0$ in the calculation below. 

Let us first recall that if a stochastic process $X(t, \omega)$ has continuous sample paths for almost every $\omega$ and satisfies  
\be \label{asymptotic X_t}
\lim_{t\rightarrow \infty}\frac{X(t, \omega)}{t} = 0 \text{ for  almost all } \omega,
\ee
then for any $\epsilon_1$ in $(0,1)$ and $\epsilon_2 > 0$, there exist $t_0 > 0$ and $\Omega_{\epsilon_1, \epsilon_2} \subset \Omega$, such that
\bea \label{large deviation X_t}
& \mathbb{P}(\Omega_{\epsilon_1, \epsilon_2}) \ge 1 - \epsilon_1,  \\
& \frac{|X(t, \omega)|}{t} < \epsilon_2, \quad \Forall t \ge t_0,  \; \omega \in \Omega_{\epsilon_1, \epsilon_2}.
\eea
The above result follows from a straightforward application of the Egoroff's Theorem by considering random variables of the form $X(t_n, \omega)/t_n$ for each non-zero rational number $t_n$.

Now, let us consider the following scaling
\be \label{BM_scaling}
\widetilde{t} = \epsilon t, \quad \widetilde{W}_{\widetilde{t}}(\omega) = \sqrt{\epsilon} W_t(\omega).
\ee

Since $\widetilde{W}$ is a Brownian motion, we can apply the above recalled general result to ensure that for any $\chi$ in $(0,1)$, there exist $T^* > 0$, such that
\bea \label{large deviation BM_scaling}
\mathbb{P}\left \{ |\widetilde{W}_{\widetilde{t}}(\omega)|  < \frac{1}{4} |\widetilde{t}| \right\}  \ge 1 - \chi/2, \quad \Forall \widetilde{t} \le -T^*. 
\eea
We obtain then using \eqref{BM_scaling} that 
\be \label{large deviation BM_scaling2}
\mathbb{P}\left \{ \sqrt{\epsilon} |W_t(\omega)| < \frac{1}{4} \epsilon |t| \right\} \ge 1 - \chi/2, \quad \Forall t \le -T^*/\epsilon. 
\ee
From Lemma~\ref{Lemma_BM_bound} (applied to Brownian motion defined on $(-\infty, 0]$), we also know that by choosing $\gamma = \sqrt{-2T^*\ln(\chi/2)}$, we get 
\be \label{Eq_BM_bound_finitetime}
\mathbb{P}\left \{ \sqrt{\epsilon} \sup_{-T^*/\epsilon \le t \le 0} |W_t(\omega)| \le \gamma \right \} \ge 1 - \chi/2.
\ee
Now, from \eqref{large deviation BM_scaling2} and \eqref{Eq_BM_bound_finitetime}, we obtain a subset $\Omega^* \subset \Omega$ with $\mathbb{P}(\Omega^*) \ge 1 - \chi$ such that for all $\omega$ in $\Omega^*$ it holds that 
\bea  \label{Eq_BM_bounds_alltime}
& \sqrt{\epsilon} |W_t(\omega)| \le \gamma, \quad \Forall t \in [-T^*/\epsilon, 0], \\
& \sqrt{\epsilon} |W_t(\omega)| < \frac{1}{4} \epsilon |t|, \quad \Forall t \le -T^*/\epsilon. 
\eea
Note that from \eqref{Mt} with $\sigma = \sqrt{\epsilon}$, we have 
\be \label{Mt2}
M_\sigma(t,\omega) = e^{-g t - \sqrt{\epsilon} W_{t}(\omega)} \int_{-\infty}^t  e^{g s' +  \sqrt{\epsilon} W_{s'}(\omega)} \mathrm{d}s'.
\ee
Using \eqref{Eq_BM_bounds_alltime} and \eqref{Mt2}, one readily obtain for any $\epsilon$ in $[0, 4 g)$ and $\omega$ in $\Omega^*$ that
\bea \label{Mt_est1}
M_\sigma (t,\omega) & \le \left(\frac{e^{2\gamma}}{g} + \frac{4e^{\gamma + (T^*/4)}}{4g - \epsilon}\right), \quad t \in [-T^*/\epsilon, 0],
\eea
and 
\be\label{Mt_est2}
M_\sigma (t,\omega) \le \frac{4 e^{-\epsilon t/2}}{4 g - \epsilon}, \quad t \in (-\infty, -T^*/\epsilon].
\ee

Now, we are in position to estimate the random steady state $a_\epsilon$ at $t=0$. Using \eqref{Eq_forcing_in_steadystate}, and \eqref{Mt_est1}--\eqref{Mt_est2}, we get
\be
\int_{-\infty}^0 M_\sigma (s, \omega) \exp(f_\epsilon(0,s,\omega)) \d s \le \frac{e^{2\gamma}}{2 \epsilon}\left(\frac{e^{2\gamma}}{g} + \frac{4 e^{\gamma + (T^*/4)}}{4g - \epsilon}\right) + \frac{4e^{-T^*}}{(4g-\epsilon) \epsilon}, \; \omega \in \Omega^*, \epsilon \in (0, 4g).
\ee
By restricting $\epsilon$ to, for instance $(0, 2g]$, we get then 
\be \label{int_Mt_est}
\int_{-\infty}^0 M_\sigma (s, \omega) \exp(f_\epsilon(0,s,\omega)) \d s \le \frac{e^{2\gamma}}{2 \epsilon}\left(\frac{e^{2\gamma}}{g} + \frac{2 e^{\gamma + (T^*/4)}}{g}\right) + \frac{2e^{-T^*}}{g\epsilon}.
\ee
Using this last estimate in \eqref{Eq_stationary_generic} with $t=0$ therein, we obtain 
\be \label{steadystate_bound1}
a_\epsilon(0,\omega) \ge c_1 \sqrt{\epsilon}, \quad  \omega \in \Omega^*, \epsilon \in (0, 2g],
\ee
with  
\be
c_1 = \left(\frac{g}{ \alpha e^{2\gamma} (e^{2\gamma}+ 2 e^{\gamma + (T^*/4)}) + 4\alpha e^{-T^*}}\right)^{1/2}.
\ee

Using again \eqref{Eq_BM_bounds_alltime}, we also have 
\be\label{Mt_est3}
 M_\sigma (t,\omega) \ge \frac{4 e^{\epsilon t/2}}{4 g + \epsilon}, \quad t \in (-\infty, -T^*/\epsilon], \; \omega \in \Omega^*,
\ee
and we can simply bound $M_\sigma (t,\omega)$ below by zero for all $t$ in $[-T^*/\epsilon, 0]$. This way, we get 
\be
\int_{-\infty}^0 M_\sigma (s, \omega) \exp(f_\epsilon(0,s,\omega)) \d s \ge  \frac{2e^{-3T^*}}{9 g \epsilon}, \quad \epsilon \in (0, 2 g], \omega \in \Omega^*,
\ee
leading thus to 
\be \label{steadystate_bound2} 
a_\epsilon(0,\omega) \le c_2 \sqrt{\epsilon}, \quad \omega \in \Omega^*, \epsilon \in (0, 2g],
\ee
with  $c_2 = \sqrt{9 g e^{3T^*}/(4\alpha)}.$ Since the constants $c_1$ and $c_2$ in respectively \eqref{steadystate_bound1} and \eqref{steadystate_bound2} depend only on $\chi$ (through $\gamma$ and $T^*$) but not on $\epsilon$, the estimate \eqref{steadystate_bound_goal} for $t=0$ follows. Since $a_\epsilon$ is stationary, the same estimation holds for all $t\ge0$. The proof is complete. 
\ep

\section{Proof of Lemma~\ref{Lem_RBC_cond_L}} \label{Sec_proof_Lemma_RBC}

\bp[Proof of Lemma~\ref{Lem_RBC_cond_L}]  

Recall that the linear operator $L_{\Ra} = - A + P_{\Ra}$ is defined in \eqref{RBC_lin_operators} and its eigenelements are provided in 
Sec.~\ref{Sec_verifying_A1-A3}. Note that $A$ leaves invariant each of the one-dimensional subspace spanned by any of the eigenvectors of $L_{\Ra}$. Thus, to verify \eqref{Eq_P_control}, it is sufficient to show that there exist $\eta$ and $\delta$, for which \eqref{Eq_P_control} holds when $u_{\s}$ is taken to be any stable eigenmode of $L_{\Ra}$.

We first verify \eqref{Eq_P_control} for $u_{\s} = \bm{e}_{0k}$ given in \eqref{RBC_beta_0k} for any $k$ in $\mathbb{N}$. Note that since $- A \bm{e}_{0k} = - k^2 \pi^2 \bm{e}_{0k} = L_\Ra \bm{e}_{0k}$, we get $P_{\Ra} \bm{e}_{0k} = (L_{\Ra} + A) \bm{e}_{0k} = 0$ for all $k$. As a result, for any $\eta$ in $(0,1)$, we have 
\bea \label{RBC_cond_L_summary1}
& \langle (- \eta A + P_{\Ra}) \bm{e}_{0k}, \bm{e}_{0k}\rangle = - \eta k^2 \pi^2 \| \bm{e}_{0k}\|^2 \le  - \eta \pi^2 \| \bm{e}_{0k}\|^2, \\
& \langle (-\eta A + P_{\Ra}) \bm{e}_{0k}, A \bm{e}_{0k}\rangle = - \eta k^2 \pi^2 \| \nabla \bm{e}_{0k}\|^2 \le - \eta \pi^2 \| \nabla \bm{e}_{0k}\|^2, 
\eea
which holds for all $\Ra$ and all $k$.

Now, for eigenfunctions in group two given by \eqref{eigen vector}, we have 
\bea \label{RBC_cond_L_est1}
\langle A \bm{e}^{\pm}_{jk}, \bm{e}^{\pm}_{jk}\rangle & = \gamma_{jk}^2 \| \bm{e}^{\pm}_{jk}\|^2, 
\eea
and
\bea \label{RBC_cond_L_est2}
\langle P_{\Ra} \bm{e}^{\pm}_{jk}, \bm{e}^{\pm}_{jk}\rangle & = \langle (L_\Ra + A) \bm{e}^{\pm}_{jk}, \bm{e}^{\pm}_{jk}\rangle = (\beta_{jk}^{\pm}(\Ra) + \gamma_{jk}^2) \| \bm{e}^{\pm}_{jk}\|^2  = \pm \sqrt{\Ra \alpha_j^2 / \gamma_{jk}^2}  \| \bm{e}^{\pm}_{jk}\|^2,
\eea
where the last equality follows from \eqref{RBC_beta_jk}.

By using \eqref{RBC_cond_L_est1} and \eqref{RBC_cond_L_est2}, we get 
\be \label{RBC_cond_L_est3}
\langle (- \eta A + P_\Ra) \bm{e}^{\pm}_{jk}, \bm{e}^{\pm}_{jk}\rangle  = \left(- \eta \gamma_{jk}^2 \pm \sqrt{\Ra \alpha_j^2 / \gamma_{jk}^2} \right)\| \bm{e}^{\pm}_{jk}\|^2.
\ee
Similarly, we have 
 \bea \label{RBC_cond_L_est3b}
\langle (- \eta A + P_\Ra) \bm{e}^{\pm}_{jk}, A \bm{e}^{\pm}_{jk}\rangle = \left(- \eta \gamma_{jk}^2 \pm \sqrt{\Ra \alpha_j^2 / \gamma_{jk}^2} \right) \|\nabla \bm{e}^{\pm}_{jk}\|^2.
\eea

To proceed further, we discuss the two cases $u_{\s} = {e}^{+}_{jk}$ and $u_{\s} = {e}^{-}_{jk}$ separately. For $u_{\s} = {e}^{-}_{jk}$, let $\eta$ in $(0,1)$ be arbitrarily fixed. Since 
\be
- \eta \gamma_{jk}^2 - \sqrt{\Ra \alpha_j^2 / \gamma_{jk}^2} \le - \eta \gamma_{jk}^2 \le -\eta \gamma_{11}^2,
\ee
we get
\bea \label{RBC_cond_L_summary2}
& \langle (- \eta A + P_\Ra) \bm{e}^{-}_{jk}, \bm{e}^{-}_{jk}\rangle \le - \eta \gamma_{11}^2 \| \bm{e}^{-}_{jk}\|^2, \\
& \langle (- \eta A + P_\Ra) \bm{e}^{-}_{jk},  A\bm{e}^{-}_{jk}\rangle \le - \eta \gamma_{11}^2 \| \nabla \bm{e}^{-}_{jk}\|^2,
\eea
which holds for all $\Ra$ and all indices $(j,k)$ in $\mathbb{N}^2$.  

For $u_{\s} = {e}^{+}_{jk}$, we aim to show that there exist an $\eta$ in $(0,1)$, $\Ra^* > \Ra$, and a $\mu > 0$ such that 
\be \label{RBC_gap_goal}
\max_{\substack{(j,k) \in \mathbb{N}^2 \\ (j,k) \neq (j_c, 1)}} \left(- \eta \gamma_{jk}^2 + \sqrt{\Ra \alpha_j^2 / \gamma_{jk}^2} \right) \le -\mu, \quad \forall \Ra \in [\Ra_c, \Ra^*].
\ee
First note that since $\gamma_{jk}^2$ is always positive and increases as $k$ increases, we have either
\be \label{RBC_gap_equiv}
\max_{\substack{(j,k) \in \mathbb{N}^2 \\ (j,k) \neq (j_c, 1)}} \left(- \eta \gamma_{jk}^2 + \sqrt{\Ra \alpha_j^2 / \gamma_{jk}^2} \right) = \max_{\substack{j \in \mathbb{N} \\ j \neq j_c}} \left( -\eta \gamma_{j1}^2 + \sqrt{\Ra \alpha_j^2 / \gamma_{j1}^2} \right),
\ee
or 
\be \label{RBC_gap_equiv2}
\max_{\substack{(j,k) \in \mathbb{N}^2 \\ (j,k) \neq (j_c, 1)}} \left(- \eta \gamma_{jk}^2 + \sqrt{\Ra \alpha_j^2 / \gamma_{jk}^2} \right) = \left( -\eta \gamma_{j_c2}^2 + \sqrt{\Ra \alpha_{j_c}^2 / \gamma_{j_c2}^2} \right).
\ee

To handle the case \eqref{RBC_gap_equiv}, let us introduce 
\be \label{RBC_gap_est0}
f_{j}(\eta,\Ra) = - \eta \gamma_{j1}^2 + \sqrt{\Ra \alpha_j^2 / \gamma_{j1}^2}. 
\ee
Since $\alpha_j^2 / \gamma_{j1}^2 \le 1$ for any $j$, and $\gamma_{j1}^2$ approaches $+\infty$ as $j$ increases, we see that for any fixed $\eta$ in $(0,1)$, we can make $f_{j}(\eta,\Ra)$ as negative as needed by taking $j$ sufficiently large. More precisely, for any $M>0$, and $\Ra_1 > \Ra_c$, there exists an index $J$ such that 
\be
f_{j}(\eta,\Ra) \le -M, 
\ee 
for all $j \ge J$ and  $\Ra$ in $[\Ra_c, \Ra_1]$. This implies that there exists $J^* > 0$ such that the maximum in \eqref{RBC_gap_equiv} is achieved at some $j \le J^*$ for all $\Ra$ in $[\Ra_c, \Ra_1]$, where $\Ra_1 > \Ra_c$ is arbitrarily fixed. We consider then the following set of functions
\be
f_j(\eta,\Ra) = - \eta \gamma_{j1}^2 + \sqrt{\Ra \alpha_j^2 / \gamma_{j1}^2},  \quad j \in \{1, \ldots, J^*\} \setminus \{j_c\}.
\ee
From now on, both cases  \eqref{RBC_gap_equiv} and \eqref{RBC_gap_equiv2} can be handled in the same way. We thus introduce also
\be
g(\eta, \Ra) = -\eta \gamma_{j_c2}^2 + \sqrt{\Ra \alpha_{j_c}^2 / \gamma_{j_c2}^2}. 
\ee

Note that $f_j(1,\Ra) = \beta_{j1}^+(\Ra)$ and $g(1,\Ra) = \beta^{+}_{j_c 2}(\Ra)$; cf.~\eqref{RBC_beta_jk}. Then, by the definition of $\Ra_c$, we know that $f_j(1, \Ra_c) < 0$ for all such $j$'s and $g(1,\Ra_c) < 0$ as well. Since $f_j$ and $g$ are continuous in both $\eta$ and $\Ra$, and we are considering only finitely many such functions, there exists $\eta$ in $(0,1)$, $\Ra^* > \Ra_c$ and $\mu>0$ for which 
\be
g(\eta, \Ra) < - \mu \; \text{ and }\;  f_j(\eta, \Ra) < - \mu, \quad j \in \{1, \ldots, J^*\} \setminus \{j_c\}, \; \Ra \in [\Ra_c, \Ra^*]. 
\ee
It follows then
\be \label{RBC_gap_est}
-\eta \gamma_{j_c2}^2 + \sqrt{\Ra \alpha_{j_c}^2 / \gamma_{j_c2}^2} < -\mu \; \text{ and }\;  \max_{\substack{j \in \mathbb{N} \\ j \neq j_c}} \left(- \eta \gamma_{j1}^2 + \sqrt{\frac{\Ra \alpha_j^2}{\gamma_{j1}^2}} \right) \le -\mu, \quad \forall \Ra \in [\Ra_c, \Ra^*].
\ee
The desired estimate \eqref{RBC_gap_goal} follows then from \eqref{RBC_gap_est}, taking into consideration \eqref{RBC_gap_equiv} and \eqref{RBC_gap_equiv2}. 

It follows then from \eqref{RBC_cond_L_est3}, \eqref{RBC_cond_L_est3b} and \eqref{RBC_gap_goal} that 
 \bea \label{RBC_cond_L_summary3}
& \langle (- \eta A + P_\Ra) \bm{e}^{+}_{jk}, \bm{e}^{+}_{jk}\rangle \le - \mu \| \bm{e}^{+}_{jk}\|^2, \\
& \langle (- \eta A + P_\Ra) \bm{e}^{+}_{jk}, A \bm{e}^{+}_{jk}\rangle \le - \mu \| \nabla \bm{e}^{+}_{jk}\|^2,
\eea
which holds for all $\Ra$ in $[\Ra_c, \Ra^*]$ and $(j,k) \neq (j_c, 1)$.  

To summarize, based on the estimates in \eqref{RBC_cond_L_summary1}, \eqref{RBC_cond_L_summary2}, and \eqref{RBC_cond_L_summary3}, by choosing $\eta$ in $(0,1)$, $\mu>0$ and $\Ra^* > \Ra_c$ for which \eqref{RBC_cond_L_summary3} holds, and letting $\delta = \min\{ \eta \pi^2, \mu, \eta \gamma_{11}^2\}$, the condition \eqref{Eq_P_control} holds for all $u_{\s}$ in $D(A)\cap H_\s$ and $\Ra$ in $[\Ra_c, \Ra^*]$. The proof is now complete. 
\ep

\bibliographystyle{amsalpha}
\input{Transitions_SPDEs_arXiv.bbl}

\end{document}

%% file: Transitions_SPDEs_arXiv.bbl
\newcommand{\etalchar}[1]{$^{#1}$}
\providecommand{\bysame}{\leavevmode\hbox to3em{\hrulefill}\thinspace}
\providecommand{\MR}{\relax\ifhmode\unskip\space\fi MR }
\providecommand{\MRhref}[2]{%
  \href{http://www.ams.org/mathscinet-getitem?mr=#1}{#2}
}
\providecommand{\href}[2]{#2}